\documentclass[11pt]{amsart}
\usepackage{amsxtra}
\usepackage{amssymb,color}
\usepackage{mathrsfs}
\usepackage{tikz}
\usetikzlibrary{positioning,calc}
\theoremstyle{plain}

\newcommand{\cleqn}{\setcounter{equation}{0}}
\newcommand{\clth}{\setcounter{theorem}{0}}
\newcommand {\sectionnew}[1]{\section{#1}\cleqn\clth}
\newcommand{\nn}{\hfill\nonumber}
\newtheorem{theorem}{Theorem}[section]
\newtheorem{lemma}[theorem]{Lemma}
\newtheorem{definition-theorem}[theorem]{Definition-Theorem}
\newtheorem{proposition}[theorem]{Proposition}
\newtheorem{corollary}[theorem]{Corollary}
\newtheorem{definition}[theorem]{Definition}
\newtheorem{example}[theorem]{Example}
\newtheorem{remark}[theorem]{Remark}
\newtheorem{conjecture}[theorem]{Conjecture}
\newtheorem{notation}[theorem]{Notation}

\newtheorem*{maintheorem*}{Main Theorem}
\newcommand \bth[1] { \begin{theorem}\label{t#1} }
\newcommand \ble[1] { \begin{lemma}\label{l#1} }

\newcommand \bpr[1] { \begin{proposition}\label{p#1} }
\newcommand \bco[1] { \begin{corollary}\label{c#1} }
\newcommand \bde[1] { \begin{definition}\label{d#1}\rm }
\newcommand \bex[1] { \begin{example}\label{e#1}\rm }
\newcommand \bre[1] { \begin{remark}\label{r#1}\rm }
\newcommand \bcj[1] { \begin{conjecture}\label{j#1}\rm }

\newcommand \bnota[1] { \begin{notation}\label{n#1}\rm }
\renewcommand {\eth} { \end{theorem} }
\newcommand {\ele} { \end{lemma} }

\newcommand {\epr} { \end{proposition} }
\newcommand {\eco} { \end{corollary} }
\newcommand {\ede} { \end{definition} }
\newcommand {\eex} { \end{example} }
\newcommand {\ere} { \end{remark} }
\newcommand {\ecj} { \end{conjecture} }

\newcommand {\enota} { \end{notation} }
\newcommand \thref[1]{Theorem \ref{t#1}}
\newcommand \leref[1]{Lemma \ref{l#1}}
\newcommand \prref[1]{Proposition \ref{p#1}}
\newcommand \coref[1]{Corollary \ref{c#1}}

\newcommand \deref[1]{Definition \ref{d#1}}
\newcommand \exref[1]{Example \ref{e#1}}
\newcommand \reref[1]{Remark \ref{r#1}}

\def \tt {t} 
\def \del {{\partial}}   

\def \Rset {{\mathbb R}}         
\def \Cset {{\mathbb C}}
\def \KK {{\mathbb K}}

\def \Zset {{\mathbb Z}}

\def \Qset {{\mathbb Q}}

\def \AA {{\mathcal{A}}} 

\def \BB {{\mathcal{B}}} 
\def \FF {{\mathcal{F}}}

\def \TT {{\mathcal{T}}}

\def \QQ {{\mathcal{Q}}}

\def \PP {{\mathcal{P}}}

\def \UU {{\mathcal{U}}}
\def \RR {{\mathcal{R}}}
\def \SS {{\mathcal{S}}}
\def \LL {{\mathcal{L}}}
\def \TT {{\mathcal{T}}} 

\def \rb {{\bf{r}}}

\def \ex {{\bf{ex}}}
\def \Scr {{\mathscr{S}}}
\def \mb {{\bf{m}}}
\def \inv {{\bf{inv}}}
\def \de {\delta}
\def \al {\alpha}
\def \be {\beta}
\def \vpi {\varpi}
\def \la {\lambda}

\def \om {\omega}
\def \Om {\Omega}
\def \ga {\gamma}
\def \de {\delta}
\def \Ga {\Gamma}

\def \sig {\sigma}

\def \vp {\varphi}
\def \ep {\epsilon}
\def \De {\Delta}
\def \rb { {\bf{r}} }

\def \mt  {\mapsto}
\def \lra {\longrightarrow}

\def \hra {\hookrightarrow}

\def \sy  {\ast}                         
\def \ci  {\circ}

\def \rcor {\rangle}
\def \lcor {\langle}

\def \ol {\overline}
\def \wt {\widetilde}
\def \wh {\widehat}

\def \id { {\mathrm{id}} }

\DeclareMathOperator \rk { {\mathrm{rk}} }
\def \Lie { {\mathrm{Lie \,}} }

\def \redd { {\mathrm{red}} }
\def \g  {\mathfrak{g}}   

\def \sl {\mathfrak{sl}}

\def \n  {\mathfrak{n}}

\def \b  {\mathfrak{b}}

\def \sl {\mathfrak{sl}}

\DeclareMathOperator \Span { {\mathrm{Span}} }
\DeclareMathOperator \Aut { {\mathrm{Aut}} }

\DeclareMathOperator \opp { {\mathrm{op}} }

\DeclareMathOperator \lt  { {\mathrm{lt}} }

\DeclareMathOperator \supp { {\mathrm{supp}} }
\DeclareMathOperator \msupp { {\mathrm{msupp}} }
\DeclareMathOperator \Fract { {\mathrm{Fract}} }

\newcommand\kx{\KK^*}

\newcommand\HH{{\mathcal{H}}}
\newcommand\xh{X(\HH)}

\newcommand \Znn {\Zset_{\ge 0}}

\def\lab{{\boldsymbol \lambda}}

\def\nub{{\boldsymbol \nu}}
\begin{document}
\title[The Berenstein--Zelevinsky conjecture]
{The Berenstein--Zelevinsky \\
quantum cluster algebra conjecture}
\author[K. R. Goodearl]{K. R. Goodearl}
\address{
Department of Mathematics \\
University of California\\
Santa Barbara, CA 93106 \\
U.S.A.
}
\email{goodearl@math.ucsb.edu}
\author[M. T. Yakimov]{M. T. Yakimov}
\thanks{The research of K.R.G. was partially supported by NSF grant DMS-0800948 and NSA grant H98230-14-1-0132,
and that of M.T.Y. by NSF grant DMS-1303036 and Louisiana Board of Regents grant Pfund-403.}
\address{
Department of Mathematics \\
Louisiana State University \\
Baton Rouge, LA 70803 \\
U.S.A.
}
\email{yakimov@math.lsu.edu}
\date{}
\begin{abstract} We prove the Berenstein--Zelevinsky conjecture that the quantized 
coordinate rings of the double Bruhat cells of all finite dimensional simple algebraic groups 
admit quantum cluster algebra structures with initial seeds as specified by \cite{BZ}.
We furthermore prove that the corresponding upper quantum cluster algebras 
coincide with the constructed quantum cluster algebras and exhibit a large number of 
explicit quantum seeds. Along the way a detailed study of the properties of quantum double Bruhat cells 
from the viewpoint of noncommutative UFDs is carried out and a quantum analog of the Fomin--Zelevinsky 
twist map is constructed and investigated for all double Bruhat cells. The results are valid 
over base fields of arbitrary characteristic and the deformation parameter is 
only assumed to be a non-root of unity.
\end{abstract}
\maketitle

\sectionnew{Introduction}
\label{intro}
Cluster Algebras are an axiomatic class of algebras which was introduced by Fomin and Zelevinsky \cite{FZ1}. In the last 10 years they have played a fundamental role 
in many different contexts in combinatorics, representation theory, geometry and mathematical physics \cite{Ke,Wi}.

One of the major applications of cluster algebras is to representation theory. Fomin and Zelevinsky \cite{FZ1} set up the program of 
(1) constructing cluster algebra structures on quantum algebras and coordinate rings of algebraic varieties, and (2) capturing
a large part of the (dual) canonical basis of Kashiwara and Lusztig via the algorithmically constructed
set of cluster monomials \cite{CKLP}. The first step 
has been carried out in several settings. Berenstein, Fomin and Zelevinsky constructed upper cluster algebra structures 
on the double Bruhat cells of all complex simple {algebraic} groups \cite{BFZ}. Geiss, Leclerc and Schr\"oer constructed cluster structures on the coordinate rings of Schubert cells for symmetric Kac--Moody algebras \cite{GLS2}. Scott and the last group of authors constructed cluster structures on the homogeneous coordinate rings of flag varieties \cite{GLS1,Sc}.  

Let $G$ be a connected, simply connected, complex simple {algebraic} group with a fixed pair of opposite Borel subgroups $B^\pm$ and a Weyl group $W$. The double Bruhat cells \cite{FZ} 
of $G$ are the subvarieties
$$
G^{u,w}:= B^+ u B^+ \cap B^- w B^-, \quad u, w \in W.
$$
Berenstein and Zelevinsky introduced the notion of a quantum cluster algebra {\cite[Sections 3-5]{BZ}} and conjectured that the quantized coordinate rings $R_q[G^{u,w}]$
of all double Bruhat cells admit {upper} quantum cluster algebras structures with an explicit initial seed which they constructed {\cite[Conjecture 10.10]{BZ}}. 
In a certain sense this is the most general situation and the cases above arise via reductions from it. The conjecture
has been proved only in very special cases. Geiss, Leclerc and Schr\"oer proved it in the case when $u=1$ and $G$ is a symmetric Kac--Moody group
\cite{GLS3}, and in a recent breakthrough Kang--Kashiwara--Kim--Oh \cite{KKKO} proved that all cluster monomials belong to the dual canonical basis, thus completing the Fomin--Zelevinsky program in this important case. {In the symmetric finite dimensional case, the latter theorem was also proved by different methods by Qin \cite{Q}.} 
The case of the conjecture for $u=1$ and general simple $G$  was established 
in \cite{BR,GYbig}. These cases are simpler and none of the existing methods allows for an extension to the situation when both $w,u \neq 1$, 
even for special $G$.

In this paper we prove the Berenstein--Zelevinsky conjecture in full generality and deduce a number of desired properties of the 
quantum cluster algebras in question. In order to formulate these results, we recall that for all fundamental weights $\vpi_i$
and pairs $u', w' \in W$, one defines the quantum minors $\De_{u' \vpi_i, w' \vpi_i} \in R_q[G^{u,w}]$. Fix reduced expressions
$$
w= s_{i_1} \ldots s_{i_N}, \quad u = s_{i'_1} \ldots s_{i'_M}
$$
and set
$$
w_{< k}^{-1} := s_{i_N} \ldots s_{i_{k+1}}, \; \; 1 \leq k \leq N,  \quad
u_{ \leq j}:= s_{i'_1} \ldots s_{i'_j}, \; \; 1 \leq j \leq M.
$$

\begin{maintheorem*}
    For all connected, simply connected, complex simple {algebraic} groups $G$, Weyl group elements $u,w \in W$, infinite base fields $\KK$ of arbitrary characteristic, 
and non-roots of unity $q \in \kx$ such that $\sqrt{q} \in \KK$, the quantized coordinate ring $R_q[G^{u,w}]$ of the double Bruhat cell $G^{u,w}$ 
is isomorphic to the quantum cluster algebra with an initial seed having 

{\rm(1)} cluster variables $\De_{\vpi_i, w^{-1} \vpi_i}$, $1 \leq i \leq r$, $\De_{\vpi_{i_k}, w^{-1}_{<k} \vpi_{i_k}}$, $1 \leq k \leq N$,
$\De_{u_{\leq j} \vpi_{i'_j}, \vpi_{i'_j}}$, $1 \leq j \leq M$,

{\rm(2)} exchange matrix equal to the Berenstein--Fomin--Zelevinsky matrix \eqref{BZ-exmatr} and 

{\rm(3)} set of mutable variables given \eqref{BZ-exch}. 

Furthermore, the corresponding upper quantum cluster algebra 
coincides with the quantum cluster algebra, and thus is isomorphic to $R_q[G^{u,w}]$.
\end{maintheorem*}

This result extends the Geiss--Leclerc--Schroer quantum cluster algebras \cite[Theorem 1.1]{GLS3} to (1) the non simply laced case, (2) when both Weyl group elements $w, u \neq 1$, and
(3) an arbitrary base field.  

In addition, we construct large sets of explicit seeds of these quantum cluster algebras. The seeds are described in
\thref{SwSu-qcluster} (a)--(b). (\prref{connections} establishes a needed isomorphism between two settings.) There is one such seed 
for each element of the following subset of the symmetric group $S_{M+N}$
\begin{equation}
\label{Xi-intro}
\Xi_{M+N} := \{ \sig \in S_{M+N} \mid 
\mbox{$\sig([1,k])$ is an interval for all $k \in [2,M+N]$} \}.
\end{equation}

The second property in 
the theorem is a very desired property of a cluster algebra and enters as {the EGM condition of} the Gross--Hacking--Keel--Kontsevich proof \cite{GHKK} of the Fock-Goncharov conjecture \cite{FG} on {theta} bases 
of cluster algebras. {More precisely}, in a companion paper we prove that the equality of cluster and upper cluster algebras for all 
double Bruhat cells also holds at the classical level 
(which is a somewhat easier result) and that one can construct (many) maximal green sequences in the sense of Keller \cite{Ke2}
using the explicit clusters in \thref{SwSu-qcluster} indexed by the set \eqref{Xi-intro}. Using these two results, the main theorem of \cite{GHKK} implies 
the validity of the Fock-Goncharov conjecture \cite{FG} for all Berenstein--Fomin--Zelevinsky cluster algebras
for double Bruhat cells $G^{u,w}$. We expect that the constructed large sets of explicit clusters of $R_q[G^{u,w}]$ will have many additional 
applications, just like Scott's clusters \cite{Sc} for Grassmannians coming from Postnikov diagrams which are used in a wide array of situations, such as
integrable systems \cite{KW}, mirror symmetry \cite{MR,RW}, scattering amplitudes \cite{Aetal}, local acyclicity \cite{MS} and many others.
 
We carry out the proof of the Main Theorem in two steps: 

(I) Using our results on quantum cluster algebra 
structures on quantum nilpotent algebras \cite{GYbig}, we construct a quantum cluster algebra structure on each $R_q[G^{u^{-1},w^{-1}}]$ and prove
that it has the 
additional property that it coincides with the corresponding upper quantum cluster algebra. 

(II) In the second step we construct a quantum analog of the Fomin--Zelevinsky twist map which
is an algebra antiisomorphism  $\zeta_{u,w} : R_q[G^{u^{-1}, w^{-1}}] \to R_q[G^{u,w}]$
and show that it transforms the quantum cluster algebra from step I to the one conjectured by Berenstein and Zelevinsky \cite{BZ}.
(Additional minor details require an interchange of the Weyl group elements $w$ and $u$, 
realized via another antiismorphism involving the antipode of $R_q[G]$.)

In the remaining part of the introduction we briefly indicate the nature of the proof. 
For every nilpotent Lie algebra $\n$, the universal enveloping algebra $U(\n)$ can be realized as an iterated skew polynomial extension 
\begin{equation}
\label{class}
U(\n) = \KK[x_1][x_2; \id, \de_2] \ldots [x_N; \id, \de_N]
\end{equation}
where all derivations $\de_k$ are locally nilpotent. Quantum nilpotent algebras are similar kinds of skew polynomial extensions 
\begin{equation}
\label{chain-intro}
R = \KK[x_1][x_2; \sig_2, \de_2] \ldots [x_N; \sig_N, \de_N]
\end{equation}
equipped with a rational action of a torus $\HH$ for which the automorphisms $\sig_k$ come from the $\HH$-action, 
the skew derivations $\de_k$ are locally nilpotent, and certain diagonal eigenvalues of the $\HH$-action are non-roots of unity (\deref{CGL}). 
The algebras \eqref{class} have all but the last property; their corresponding eigenvalues are all equal to 1. In \cite{GYbig} 
we proved that all quantum nilpotent algebras (under mild conditions) possess quantum cluster algebra structures where the initial cluster 
variables are sequences of the form $p_1, \ldots, p_N$ such that $p_k$ is the unique homogeneous prime element 
in $R_k$ that is not in $R_{k-1}$ (for the intermediate algebras $R_k$ generated by $x_1, \ldots, x_k$). These arguments are based on the 
theory of noncommutative UFDs.

To carry out step I, we first realize each quantum double Bruhat cell $R_q[G^{u,w}]$ in terms of a bicrossed product 
\begin{equation}
\label{bicross}
\UU^-[u]_{\opp} \bowtie \UU^+[w].
\end{equation}
The algebras $\UU^\pm[w]$ are the {\em{quantum Schubert cell algebras}} (also called {\em{quantum unipotent cells}}) introduced by 
De Concini--Kac--Procesi \cite{DKP} and Lusztig \cite{L}. They are deformations of the universal enveloping algebras
$\UU(\n_\pm \cap w(\n_\mp))$ where $\n_\pm$ are the unipotent radicals of the Borel subalgebras $\Lie (B_\pm)$. 
Here $(.)_{\opp}$ refers to the opposite algebra structure. Part I of the 
proof proceeds in the following steps:

(A) Proving that the algebras $\UU^-[u]_{\opp} \bowtie \UU^+[w]$ are quantum nilpotent algebras with respect to a particular ordering of the Lusztig root vectors of 
$\UU^-[u]$ and $\UU^+[w]$.

(B) Showing that for different pairs of Weyl group elements $w, u \in W$, the algebras \eqref{bicross} embed in each other 
appropriately and that the intermediate algebras $R_k$ in the extension \eqref{chain-intro} are algebras of the same type for particular 
subwords of $w$ and $u$.

(C) Classifying the prime elements of the algebras $\UU^-[u]_{\opp} \bowtie \UU^+[w]$ using techniques from the study of the spectra of quantum groups 
\cite{HLT,J,Y-sqg}.

The Fomin--Zelevinsky twist map \cite{FZ} is a birational isomorphism $\zeta_{u,w} : G^{u,w} \to G^{u^{-1}, w^{-1}}$ which, in addition to cluster algebras, 
plays a major role in many other constructions. {In this paper we construct a quantum version of the Fomin--Zelevinsky twist map for all (quantum) double Bruhat cells. 
The second part of 
the proof of the Main Theorem relies on proving a compatibility property of it} with respect to two different kinds of quantum cluster algebras 
(the ones that are the result of step I and the conjectured ones by Berenstein and Zelevinsky).
The original twist map is defined in terms of set theoretic factorizations in $G$, not in terms of maps on coordinate rings that are manifestly homomorphisms of algebras. This makes the construction of a quantum twist map rather 
delicate. For the construction of the map we employ the open embedding of the reduced double Bruhat cell
$$
G^{u,w}/H \hra (B^+ u B^+)/B^+ \times (B^- w B^-)/B^- \subset G/B^+ \times G/B^-
$$
obtained from the projections of $G$ to the flag varieties $G/B^\pm$ (here $H := B^+ \cap B^-$). We express the restriction of the action of the classical twist 
map to the pullbacks of the two projections in geometric terms and then quantize that picture. The fact that the quantum twist map is an 
algebra antiisomorphism $R_q[G^{u,w}] \cong R_q[G^{u^{-1}, w^{-1}}]$ is established by first showing that this happens at the level of
division rings of fractions and then going back to the original algebras. {After the completion of the paper, Kimura and Oya \cite{KO}
proved that a variant of the quantum twist map in the special case of one Weyl group element maps the dual canonical basis to itself.
It is an important problem to investigate the properties of the quantum twist map
for quantum double Bruhat cells with respect to dual canonical bases.}

The paper is organized as follows. Sect. \ref{qnilp} contains necessary background on quantum cluster algebras and a brief description of 
our construction of cluster structures on quantum nilpotent algebras \cite{GYbig}. Sect. \ref{qdoubB} contains background on quantum groups and the realization of the 
quantum double Bruhat cell algebras $R_q[G^{w,u}]$ in terms of bicrossed products $\UU^-[u]_{\opp} \bowtie \UU^+[w]$. Sect. \ref{qdoubB}--\ref{prime} carry out steps (A), (B) and (C) above. Sect. \ref{qS+S-} constructs a cluster algebra structure on each of the algebras $\UU^-[u]_{\opp} \bowtie \UU^+[w]$, proves that the 
corresponding upper quantum cluster algebra equals the quantum cluster algebra and constructs seeds associated to the elements of the set \eqref{Xi-intro}. Sect. \ref{twist} constructs the quantum twist map. Sect. \ref{App} and \ref{QcldoubleB} describe the connection between the cluster algebra structure on 
$\UU^-[w]_{\opp} \bowtie \UU^+[u]$ and the conjectured one on $R_q[G^{u,w}]$ and contain the proof of the Berenstein--Zelevinsky conjecture, 
\thref{BZc}. To establish this connection, a second isomorphism $R_q[G^{u,w}] \to R_q[G^{w^{-1},u^{-1}}]$ is used in Sect. \ref{App}, which is the product of the antipode and Cartan involution of $R_q[G]$. This one is simpler than the quantum twist map and only plays an auxiliary role. However, as a result of it, 
in the body of the paper the roles of $u$ and $w$ are interchanged, compared to the introduction.

Throughout the paper we work over an infinite field $\KK$ of arbitrary characteristic. For a $\KK$-algebra $A$ and elements $a_1, \ldots, a_k \in A$, we denote by 
$\KK \lcor a_1, \ldots, a_k \rcor$ the unital subalgebra of $A$ generated by $a_1, \ldots, a_k$. For two integers $m \leq n$, denote $[m,n] := \{m,  m+1, \ldots, n\}$.
\medskip

\noindent
{\bf Acknowledgement.} {We are grateful to the referee for making a number of valuable suggestions which helped us to improve the exposition.}
\sectionnew{Quantum cluster algebra structures on quantum nilpotent algebras}
\label{qnilp}
In this section we review quantum cluster algebras and noncommutative unique factorization 
domains. We describe the key parts of
our construction of quantum cluster algebra structures on all
symmetric quantum nilpotent algebras (alias CGL extensions) \cite{GYbig}, which plays a key role in the paper.
We also derive some additional general results to those in \cite{GYbig}.

\subsection{Quantum cluster algebras}
\label{q-cl}
{We work in a multiparameter setting which extends the uniparameter case originally developed by Berenstein and Zelevinsky \cite{BZ}.}
Fix a positive integer $N$, and denote by  $\{ e_1, \dots, e_N \}$ the standard basis of $\Zset^N$. Write elements of $\Zset^N$ as column vectors.

A matrix $\rb \in M_N(\kx)$ is called {\em{multiplicatively skew-symmetric}} if $r_{kk} = 1$ and $r_{jk} = r_{kj}^{-1}$, for all $j,k \in [1,N]$. Corresponding to any such matrix is a {\em{skew-symmetric bicharacter}}
\begin{equation}
\label{Om}
\Om_\rb : \Zset^N \times \Zset^N \rightarrow \kx, \quad \text{given by} \quad \Om_\rb(e_k, e_j) = r_{kj}, \; \; \forall\,  j,k \in [1,N].
\end{equation}

Let $\FF$ be a division algebra over $\KK$. A {\em{toric frame}} for $\FF$ is a mapping 
$$
M : \Zset^N \to \FF
$$
such that 
\begin{equation}
\label{MOm}
M(f) M(g) = \Om_\rb(f,g) M(f+g), \quad
\forall\,  f,g \in \Zset^N,
\end{equation}
for some multiplicatively  skew-symmetric bicharacter $\Om_\rb : \Zset^N \times \Zset^N \to \kx$, 
and such that the elements in the image of $M$ are linearly independent over $\KK$ and $\Fract \KK\langle M(\Zset^N) \rangle = \FF$,
\cite{BZ,GYbig}. The matrix $\rb$ is uniquely reconstructed from the toric frame $M$, 
and will be denoted by $\rb(M)$. The elements $M(e_1), \ldots M(e_N)$ are called \emph{cluster variables}. 
Fix a subset $\ex \subset [1,N]$, to be called the set of {\em{exchangeable  
indexes}}; the indices in $[1,N] \backslash \ex$ will be called {\em{frozen}}. An integral 
$N \times \ex$ matrix $\wt{B}$ will be called an {\em{exchange matrix}} if its principal part (the $\ex\times\ex$ submatrix) is 
skew-symmetrizable.
A \emph{quantum seed} for $\FF$ is a pair $(M, \wt{B})$ consisting of a toric frame for $\FF$ 
and an exchange matrix such that 
\begin{align*}
&
\Om_\rb(b^k, e_j) = 1, \; \; \forall\,  k \in \ex, \; j \in [1,N], \; k \neq j
\quad \mbox{and}
\\
&\Om_\rb(b^k, e_k) \; \; 
\mbox{are not roots of unity}, \; \; \forall\,  k \in \ex
\end{align*}
where $b^k$ denotes the $k$-th column of $\wt{B}$.  

{
The \emph{mutation in direction $k \in \ex$} of a quantum seed $(M,\wt{B})$ is the quantum seed $(\mu_k(M), \mu_k(\wt{B}))$ where $\mu_k(M)$ is described below and $\mu_k(\wt{B})$ is the $N\times\ex$ matrix $(b'_{ij})$ with entries
$$
b'_{ij} :=
\begin{cases}
- b_{ij}, &\mbox{if} \; \;
i=k \; \; \mbox{or} \; \; j=k
\\
b_{ij} + \frac{ |b_{ik}| b_{kj} + b_{ik} | b_{kj}|}{2},
&\mbox{otherwise},
\end{cases}
$$
\cite{FZ1}. Corresponding to the column $b^k$ of $\wt{B}$ are automorphisms $\rho_{b^k,\pm}$ of $\FF$ such that
$$
\rho_{b^k,\ep}(ME_\ep(e_j)) = 
\begin{cases}
ME_\ep(e_k) + ME_\ep(e_k + \ep b^k), \quad &\text{if} \ j = k  \\
ME_\ep(e_j), \quad &\text{if} \ j \ne k,
\end{cases}
$$
\cite[Proposition 4.2]{BZ} and 
\cite[Lemma 2.8]{GYbig}, where $E_\ep = E^{\wt{B}}_\ep$ is the $N\times N$ matrix with entries
$$
(E_\ep)_{ij} =
\begin{cases}
\delta_{ij}, & \mbox{if} \; j \neq k \\
-1, & \mbox{if} \; i=j=k \\
\max(0, - \ep b_{ik}), & \mbox{if} \; i \neq j = k.
\end{cases}
$$
The toric frame $\mu_k(M)$ is defined as
$$
\mu_k(M) := \rho_{b^k,\ep} M E_\ep : \Zset^N \rightarrow \FF.
$$
It is independent of the choice of $\ep$, and, paired with $\mu_k(\wt{B})$, forms a quantum seed \cite[Proposition 2.9]{GYbig}. (See also \cite[Corollary 2.11]{GYbig}, and compare with \cite[Proposition 4.9]{BZ} for the uniparameter case.)
}

Fix a subset $\inv$ of the set $[1,N] \backslash \ex$ of frozen indices -- the corresponding 
cluster variables will be inverted. The {\em{quantum cluster algebra}} 
$\AA(M, \wt{B}, \inv)$ is the unital subalgebra of $\FF$ generated by the cluster
variables of all seeds obtained from $(M, \wt{B})$ by mutation and by 
$\{M(e_k)^{-1} \mid k \in \inv \}$. To each quantum seed $(M, \wt{B})$ one associates the mixed quantum 
torus/quantum affine space algebra
$$
\TT_{(M, \wt{B})}:= \KK \lcor M(e_k)^{\pm 1}, \;  M(e_j) \mid k \in \ex \cup \inv, \; j \in [1,N] \backslash (\ex \cup \inv) \rcor
\subset \FF.
$$
The intersection of all such subalgebras of $\FF$ associated to all seeds that are obtained by mutation 
from the seed $(M,\wt{B})$ is called the {\em{upper quantum cluster algebra}} of 
$(M,\wt{B})$ and is denoted by $\UU(M, \wt{B}, \inv)$.

\subsection{Equivariant noncommutative unique factorization domains}
\label{ncufd}

Recall that a {\em{prime element}} of a domain $R$ is any nonzero normal element $p\in R$ (\emph{normality} meaning that $Rp = pR$)
such that $Rp$ is a completely prime ideal, i.e., $R/Rp$ is a domain.
Assume that in addition $R$ is a $\KK$-algebra and $\HH$ a group acting on $R$ 
by $\KK$-algebra automorphisms. An {\em{$\HH$-prime ideal}} of $R$ is any proper 
$\HH$-stable ideal $P$ of $R$ such that $(IJ\subseteq P \implies I\subseteq P$ 
or $J\subseteq P)$ for all $\HH$-stable ideals $I$ and $J$ of $R$. 
One says that $R$ is an {\em{$\HH$-UFD}} if each nonzero $\HH$-prime ideal 
of $R$ contains a prime $\HH$-eigenvector. This is an equivariant 
version of Chatters' notion \cite{Cha} of noncommutative unique 
factorization domain given in \cite[Definition 2.7]{LLR}.

The following fact is an equivariant version of results of 
Chatters and Jordan \cite[Proposition 2.1]{Cha},
\cite[p. 24]{ChJo}, see \cite[Proposition 2.2]{GY1} and
\cite[Proposition 6.18 (ii)]{Y-sqg}. 

\bpr{factorHUFD}
Let $R$ be a noetherian $\HH$-UFD. Every normal $\HH$-eigenvector in $R$ 
is either a unit or a product of prime $\HH$-eigenvectors. The factors are unique 
up to reordering and taking associates.
\epr

\subsection{CGL extensions}
\label{CGLext}
Let $S$ be a unital $\KK$-algebra. We use the standard notation $S[x;\theta,\de]$ for a \emph{skew polynomial algebra}, or \emph{skew polynomial extension}; it denotes a $\KK$-algebra generated by $S$ and an element $x$ satisfying $xs = \theta(s) x + \de(s)$ for all $s \in S$, where $\theta$ is an algebra endomorphism of $S$ and $\de$ is a {\rm(}left\/{\rm)} $\theta$-derivation of $S$. The algebra $S[x;\theta,\de]$ is a free left $S$-module, with the nonnegative powers of $x$ forming a basis. 

We focus on iterated skew polynomial extensions
\begin{equation} 
\label{itOre}
R := \KK[x_1][x_2; \theta_2, \delta_2] \cdots [x_N; \theta_N, \delta_N],
\end{equation}
where it is taken as implied that $\KK[x_1] = \KK[x_1; \theta_1, \delta_1]$ with $\theta_1 = \id_\KK$ and $\de_1 = 0$. For $k \in [0,N]$, set
$$
R_k := \KK \langle x_1, \dots, x_k \rangle = \KK[x_1][x_2; \theta_2, \delta_2] \cdots [x_k; \theta_k, \delta_k].
$$
In particular, $R_0 = \KK$.

\bde{CGL} An iterated skew polynomial extension \eqref{itOre}
is called a \emph{CGL extension} 
\cite[Definition 3.1]{LLR} or a \emph{quantum nilpotent algebra} if it is equipped with a rational action of a $\KK$-torus $\HH$ 
by $\KK$-algebra automorphisms satisfying the following conditions:
\begin{enumerate}
\item[(i)] The elements $x_1, \ldots, x_N$ are $\HH$-eigenvectors.
\item[(ii)] For every $k \in [2,N]$, $\de_k$ is a locally nilpotent 
$\theta_k$-derivation of the algebra $R_{k-1}$. 
\item[(iii)] For every $k \in [1,N]$, there exists $h_k \in \HH$ such that 
$\theta_k = (h_k \cdot)|_{R_{k-1}}$ and the $h_k$-eigenvalue of $x_k$, to be denoted by $\la_k$, is not a root of unity.
\end{enumerate}

Conditions (i) and (iii) imply that 
$$
\theta_k(x_j) = \la_{kj} x_j \; \; \mbox{for some} \; \la_{kj} \in \kx, \; \; \forall\,  1 \le j < k \le N.
$$
We then set $\la_{kk} :=1$ and $\la_{jk} := \la_{kj}^{-1}$ for $j< k$.
This gives rise to a multiplicatively skew-symmetric 
matrix $\lab := (\la_{kj}) \in M_N(\kx)$ and the corresponding skew-symmetric bicharacter $\Om_\lab$ from \eqref{Om}.

Define the \emph{length} of $R$ to be $N$ and the  \emph{rank} of $R$ by 
\begin{equation}
\label{rk}
\rk(R) := \{ k \in [1,N] \mid \de_k = 0 \} \in \Zset_{ > 0}
\end{equation}
(cf.~\cite[Eq. (4.3)]{GY1}).
Denote the character group of the torus $\HH$ by $\xh$.  The action of 
$\HH$ on $R$ gives rise to an $\xh$-grading of $R$, and the $\HH$-eigenvectors 
are precisely the nonzero homogeneous elements with respect to this grading. 
The $\HH$-eigenvalue of a nonzero homogeneous element $u \in R$ will be denoted by $\chi_u$. In other words, $\chi_u = \xh\mbox{-deg}(u)$ in terms of the $\xh$-grading.
\ede 

By \cite[Proposition 3.2, Theorem 3.7]{LLR}, every CGL extension is an 
$\HH$-UFD. A recursive description of the sets of homogeneous prime elements 
of the intermediate algebras $R_k$ of a CGL extension $R$ was obtained in  
\cite{GY1}. The statement of the result involves the standard predecessor and successor functions, $p = p_\eta$ and $s = s_\eta$, of a function $\eta : [1,N] \to \Zset$, defined as follows:
\begin{equation}
\label{pred.succ}
\begin{aligned}
p(k) &= \max \{ j <k \mid \eta(j) = \eta(k) \}, \\
s(k) &= \min \{ j > k \mid \eta(j) = \eta(k) \}, 
\end{aligned}
\end{equation}
where $\max \varnothing = -\infty$ and $\min \varnothing = +\infty$. Define corresponding order functions $O_\pm : [1,N] \rightarrow \Znn$ by
\begin{equation}
\label{O-+}
\begin{aligned}
O_-(k) &:= \max \{ m \in \Znn \mid p^m(k) \ne -\infty \},  \\
O_+(k) &:= \max \{ m \in \Znn \mid s^m(k) \ne +\infty \}.
\end{aligned}
\end{equation}

\bth{1} \cite[Theorem 4.3]{GY1} Let $R$ be a CGL extension of length $N$ and rank $\rk(R)$ as 
in \eqref{itOre}. There exist a function $\eta : [1,N] \to \Zset$ 
whose range has cardinality $\rk(R)$ and elements
$$
c_k \in R_{k-1} \; \; \mbox{for all} \; \; k \in [2,N] \; \; 
\mbox{with} \; \; p(k) \neq - \infty
$$
such that the elements $y_1, \ldots, y_N \in R$, recursively defined by 
\begin{equation}
\label{y}
y_k := 
\begin{cases}
y_{p(k)} x_k - c_k, &\mbox{if} \; \;  p(k) \neq - \infty \\
x_k, & \mbox{if} \; \; p(k) = - \infty,  
\end{cases}
\end{equation}
are homogeneous and have the property that for every $k \in [1,N]$,
\begin{equation}
\label{prime-elem}
\{y_j \mid j \in [1,k] , \, s(j) > k \}
\end{equation}
is a list of the homogeneous prime elements of $R_k$ up to scalar multiples.

The elements $y_1, \ldots, y_N \in R$ with these properties are unique.
The function $\eta$ satisfying the above 
conditions is not unique, but the partition of $[1,N]$ into a disjoint 
union of the level sets of $\eta$ is uniquely determined by $R$, as are the predecessor and successor functions $p$ and $s$.
The function $p$ has the property that $p(k) = - \infty$
if and only if $\de_k =0$.
\eth

The uniqueness of the level sets of $\eta$ was not stated in \cite[Theorem 4.3]{GY1}, 
but it follows at once from \cite[Theorem 4.2]{GY1}. This uniqueness immediately implies the uniqueness of $p$ and $s$.
In the setting of the theorem, the rank of $R$ is also given by
\begin{equation}
\label{rankR}
\rk(R) = |\{ j \in [1,N] \mid s(j) > N \}|
\end{equation}
\cite[Eq. (4.3)]{GY1}.

\bex{sl3a} {Consider the unital $\KK$-algebra $R$ with four generators $X_1^\pm$, $X_2^\pm$ and the relations
\begin{gather}
(X_i^\pm)^2 X_{3-i}^{\pm} - (q + q^{-1}) X_i^\pm X_{3-i}^\pm X_i^\pm + X_{3-i}^\pm (X_i^\pm)^2=0,  \label{ex-rels1}  \\
\begin{gathered}
X_i^+ X_{3-i}^- = q^{-1} X_{3-i}^- X_i^+,  \\ 
X_i^+ X_i^- = q^2 X_i^- X_i^+ + (q^{-1} - q)^{-1} 
\end{gathered}
\label{ex-rels23}
\end{gather}
for $i =1,2$, where $q \in \kx$ is not a root of unity.
It has an action of $(\kx)^2$ by algebra automorphisms given by
\[
(s_1, s_2) \cdot X_i^\pm := s_i^{\pm 1} X_i^\pm, \quad i =1,2 
\]
for $s_1, s_2 \in \kx$ 
and an involutive $\KK$-algebra antiautomorphism $\gamma$ given by 
\[
\gamma(X_i^\pm) := X_i^\mp, \quad i =1,2.
\]
Denote the elements 
\begin{equation}
\label{X12}
X_{12}^+:= X_1^+ X_2^+ - q^{-1} X_2^+ X_1^+ \quad \mbox{and} \quad X_{12}^-:= X_1^- X_2^- - q X_2^- X_1^-.
\end{equation}
The antiautomorphism $\gamma$ interchanges the two elements after rescaling: $\gamma(X_{12}^+) = - q^{-1} X_{12}^-$. The 
six elements 
\[
X_1^\pm,\ X_2^\pm,\ X_{12}^\pm
\]
satisfy the relations
\begin{equation}
\begin{gathered}
\begin{aligned}
X_2^+ X_1^+ &= q X_1^+ X_2^+ - q X_{12}^+, \qquad  &X_1^- X_2^- &= q X_2^- X_1^- + X_{12}^-, \\
X_{12}^+ X_1^+ &= q^{-1} X_1^+ X_{12}^+,  &X_2^+ X_{12}^+ &=  q^{-1} X_{12}^+ X_2^+, \\
X_{12}^- X_2^- &= q^{-1} X_2^- X_{12}^-,  &X_1^- X_{12}^- &=  q^{-1} X_{12}^- X_1^-, \\
X_{12}^+ X_1^- &= q X_1^- X_{12}^+,  &X_{12}^+ X_2^- &= q X_2^- X_{12}^+ - q^{-1} X_1^+, \\
X_1^+ X_{12}^-  &= q X_{12}^- X_1^+,  &X_2^+ X_{12}^-  &= q X_{12}^- X_2^+  + X_1^-,
\end{aligned}  \\
X_{12}^+ X_{12}^- = q^2 X_{12}^- X_{12}^+ + (q^2 -1) X_1^- X_1^+ + (q^{-1} - q)^{-1} .
\end{gathered}
\label{ex-Orerels}
\end{equation}
Conversely, the $\KK$-algebra presented by the six generators
\begin{equation}
\label{6gen}
X_2^-,\ X_{12}^-,\ X_1^-,\ X_1^+,\ X_{12}^+,\ X_2^+
\end{equation}
and the relations \eqref{ex-Orerels}, \eqref{ex-rels23} also satisfies \eqref{ex-rels1}. Thus, $R$ can be presented by the generators \eqref{6gen} and the $15$ relations \eqref{ex-Orerels}, \eqref{ex-rels23}.
One easily deduces from these relations that $R$ is a CGL extension with the 
above action of $(\kx)^2$ when the six generators are adjoined in the order \eqref{6gen}. There are other possible ways to adjoin these generators that satisfy \deref{CGL}, but we will see in \exref{sl3b}
that this presentation has an additional symmetry property.}

{For $j \in [3,6]$ we have
$\la_j = q^{-2}$. The choice of the elements $h_1, h_2 \in (\kx)^2$ for this algebra that satisfy the conditions in \deref{CGL} is non-unique. 
(The choice of the other elements $h_j  \in (\kx)^2$ is unique.) If one chooses
$h_1 := (1,q^2)$ and $h_2 := (q,q)$, then we also have $\la_j = q^{-2}$ for $j =1,2$.}

{The corresponding sequence of homogeneous prime elements from \thref{1} is 
\begin{equation}  \label{ex-yelts}
\begin{aligned}
y_1 &= X_1^-, \qquad y_2 = X_{12}^-, \qquad y_3= X_2^- X_1^- - t X_{12}^-, 
\\
y_4 &= X_2^- X_1^- X_1^+ - t X_{12}^- X_1^+ - q^{-1} t^2 X_2^-, \qquad 
y_5= X_{12}^- X_{12}^+ + X_1^- X_1^+ - q^{-1} t^2, \\
y_6 &= X_2^- X_1^- X_1^+ X_2^+ - t X_{12}^-X_1^+ X_2^+ - q^{-1} t^2 X_2^- X_2^+  \\
&\hphantom{= X_2^- X_1^- X_1^+ X_2^+} \ + q t X_2^- X_1^- X_{12}^+ - q t^2 X_{12}^- X_{12}^+ + q^{-2} t^4,
\end{aligned}
\end{equation}
where 
\begin{equation}
\label{t-notation}
t: = (q^{-1} - q)^{-1}.
\end{equation}
The sequence is recursively computed by determining the elements $c_k$ from \thref{1} by using \cite[Proposition 4.7(b)]{GY1}.
From this, one sees that the function $\eta : [1,6] \to \Zset$ can be chosen to be
\[
\eta(1) = \eta(3) = \eta(4) = \eta(6) :=1, \; \; \eta(2) = \eta(5) :=2. 
\]
}
\eex

\bde{lterm}
Denote by $\prec$ the reverse lexicographic order on $\Znn^N$:
\begin{multline}
\label{prec}
(m'_1, \ldots, m'_N) \prec (m_1, \ldots, m_N)
\quad \mbox{if there exists} \\
\mbox{$n \in [1,N]$ with $m'_n < m_n, \, m'_{n+1} = m_{n+1}, \, \ldots, \, m'_N = m_N$.}
\end{multline}

A CGL extension $R$ as in \eqref{itOre} has the $\KK$-basis 
$$
\{ x^f := x_1^{m_1} \cdots x_N^{m_N} \mid 
f = (m_1, \ldots, m_N)^T \in \Znn^N \}.
$$
We say that a nonzero element
$b \in R$ 
has \emph{leading term} $t x^f$ where $t \in \kx$ and $f \in \Znn^N$ 
if 
$$
b = t x^f + \sum_{g \in \Znn^N,\; g \prec f} t_g x^g
$$
for some $t_g \in \KK$, and we set $\lt(b) := t x^f$.

The leading terms of the prime elements $y_k$ in \thref{1} are given by 
\begin{equation}
\label{lt-y}
\lt(y_k) = x_{p^{O_-(k)}(k)} \ldots x_{p(k)} x_k, \; \; \forall\,  k \in [1,N].
\end{equation}
\ede

The leading terms of reverse-order monomials $x_N^{m_N} \ldots x_1^{m_1}$ involve symmetrization scalars in $\kx$ defined by
\begin{equation}
\label{Scrlab}
\Scr_\lab(f) := \prod_{1 \le j< k \le N} \la_{jk}^{- m_j m_k}, \; \; \forall\,  f = (m_1,\dots,m_N)^T \in \Zset^N.
\end{equation}
Namely,
\begin{equation}
\label{x-comm}
\lt( x_N^{m_N} \ldots x_1^{m_1} ) = \Scr_\lab( (m_1,\dots,m_N)^T ) x_1^{m_1} \ldots x_N^{m_N}, \; \; \forall\,  (m_1,\dots,m_N)^T \in \Zset^N.
\end{equation}

\subsection{Symmetric CGL extensions}
\label{symmCGLext}

Given an iterated skew polynomial extension $R$ as in \eqref{itOre}, denote
$$
R_{[j,k]} := \KK \langle x_i \mid j \le  i \le k \rangle, \; \; \forall\,  j,k \in [1,N].
$$
So, $R_{[j,k]} = \KK$ if $j \nleq k$.

\bde{symCGL} We call a CGL extension $R$ of length $N$ as in 
\deref{CGL} {\em{symmetric}} if the following two conditions hold:
\begin{enumerate}
\item[(i)] For all $1 \leq j < k \leq N$,
$$
\de_k(x_j) \in R_{[j+1, k-1]}.
$$
\item[(ii)] For all $j \in [1,N]$, there exists $h^\sy_j \in \HH$ 
such that 
$$
h^\sy_j \cdot x_k = \la_{kj}^{-1} x_k = \la_{jk} x_k, \; \; \forall\,  
k \in [j+1, N]
$$
and $h^\sy_j \cdot x_j = \la^\sy_j x_j$ for some $\la^\sy_j \in \kx$ which is not a root of unity.
\end{enumerate}

The above conditions imply that $R$ has a CGL extension presentation with the variables $x_k$ in descending order:
$$R = \KK[x_N] [x_{N-1}; \theta^*_{N-1}, \de^*_{N-1}] \cdots [x_1; \theta^*_1, \de^*_1],$$
see \cite[Corollary 6.4]{GY1}.
\ede
\bex{sl3b} {The algebra from \exref{sl3a} is a symmetric CGL extension when its generators $X_1^\pm, X_2^\pm, X_{12}^\pm$ are 
adjoined in the order \eqref{6gen}. This directly follows from the straightening relations between those generators listed in the example.}
\eex

\bpr{la} {\rm\cite[Proposition 5.8]{GYbig}}
Let $R$ be a symmetric CGL extension of length $N$.
If $l \in [1,N]$ with $O_+(l) = m > 0$, then
\begin{equation}
\label{la-eq}
\la^*_l = \la^*_{s(l)} = \ldots = \la^*_{s^{m-1}(l)} = 
\la^{-1}_{s(l)} = \la^{-1}_{s^2(l)}= \ldots = \la^{-1}_{s^m(l)}.
\end{equation}
\epr

\bde{XiN}
Denote the following subset of the symmetric group $S_N$:
\begin{equation}
\label{tau}
\begin{aligned}
\Xi_N := \{ \sig \in S_N \mid 
\sig(k) &= \max \, \sig( [1,k-1]) +1 \; \;
\mbox{or} 
\\
\sig(k) &= \min \, \sig( [1,k-1]) - 1, 
\; \; \forall\,  k \in [2,N] \}.
\end{aligned}
\end{equation}
In other words, $\Xi_N$ consists of those $\sig \in S_N$ 
such that $\sig([1,k])$ is an interval for all $k \in [2,N]$. 

If $R$ is a symmetric CGL extension of length $N$, then for each $\sig \in \Xi_N$ we have a CGL extension presentation
\begin{equation}
\label{tauOre}
R = \KK [x_{\sig(1)}] [x_{\sig(2)}; \theta''_{\sig(2)}, \de''_{\sig(2)}] 
\cdots [x_{\sig(N)}; \theta''_{\sig(N)}, \de''_{\sig(N)}],
\end{equation}
see \cite[Remark 6.5]{GY1}, \cite[Proposition 3.9]{GYbig}.

Denote by $\Ga_N$ the following subset of $\Xi_N$:
\begin{equation}
\label{GaN}  
\begin{aligned}
\Ga_N &:= \{ \sig_{i,j} \mid 1 \leq i \leq j \leq N \}, \; \; \text{where}  \\
\sig_{i,j} &:= [i+1, \ldots, j, i, j+1, \ldots, N, i-1, i-2, \ldots, 1] . 
\end{aligned}
\end{equation}
\ede

If $R$ is a symmetric CGL extension of length $N$ and $1 \le i \le k \le N$, then the subalgebra $R_{[i,k]}$ of $R$ is a symmetric CGL extension, and \thref{1} may be applied to it.  Moreover, for the case $k = s^m(i)$ we have

\bpr{yjk} {\rm\cite[Theorem 5.1]{GYbig}} Assume that $R$ is a symmetric CGL extension of length $N$, and $i \in [1,N]$ and $m \in \Znn$ are such that $s^m(i) \in [1,N]$. Then there is a unique homogeneous prime element $y_{[i,s^m(i)]} \in R_{[i,s^m(i)]}$ such that

{\rm{(i)}} $y_{[i,s^m(i)]} \notin R_{[i,s^m(i) -1]}$ and $y_{[i,s^m(i)]} \notin R_{[i+1, s^m(i) ]}$.

{\rm{(ii)}} $\lt( y_{[i,s^m(i)]} ) = x_i x_{s(i)} \cdots x_{s^m(i)}$.
\epr

Certain combinations of the homogeneous prime elements from the above proposition play an important role in the mutation formulas for quantum cluster variables of symmetric CGL extensions. They are given  in the following theorem, where we denote
\begin{equation}
\label{e-interval}
\begin{aligned}
e_{[j, s^l(j)]} := \; &e_j + e_{s(j)} + \cdots + e_{s^l(j)},  \\
 &\forall\,  j \in [1,N], \; l \in \Znn \; \; \text{such that} \; \; s^l(j) \in [1,N].
 \end{aligned}
\end{equation}

\bth{uismi}
{\rm\cite[Corollary 5.11]{GYbig}}
Assume that $R$ is a symmetric CGL extension of length $N$,
and $i \in [1,N]$ and $m \in \Zset_{>0}$ are such that $s^m(i) \in [1,N]$.
Then
\begin{equation}
\label{uuu}
u_{[i,s^m(i)]} := 
y_{[i, s^{m-1}(i)]} y_{[s(i), s^m(i)]} - \Om_\lab(e_i, e_{[s(i), s^{m-1}(i)]}) 
y_{[s(i), s^{m-1}(i)]} y_{[i,s^m(i)]} 
\end{equation}
is a nonzero homogeneous normal element of $R_{[i+1,s^m(i)-1]}$ 
which is not a multiple of $y_{[s(i), s^{m-1}(i)]}$ if $m \geq 2$. It normalizes the elements of $R_{[i+1,s^m(i)-1]}$
in exactly the same way as $y_{[s(i), s^{m-1}(i)]} y_{[i,s^m(i)]}$ does.
\eth
\bex{sl3c} {Continuing \exref{sl3a}, we illustrate \thref{uismi} for the algebra in the example and $i=1$, $m=3$. From the $\eta$-function in the 
example, we obtain $s(1)=3$, $s^2(1) = 4$ and $s^3(1) = 6$. By a direct calculation one gets
\begin{align*}
&y_{[1,4]} = y_4, \; \; y_{[1,6]} = y_6, \; \; y_{[3,4]} = X_1^- X_1^+ - q^{-1} t^2,
\\
&y_{[3,6]} = X_1^- X_1^+ X_2^+ - q^{-1} t^2 X_2^+ + q t X_1^- X_{12}^+ 
\end{align*}
using the notation \eqref{t-notation}. This gives
\[
u_{[1,6]} = y_{[1,4]} y_{[3,6]} - y_{[3,4]} y_{[1,6]} = - q^{-2} t^4 y_5 = - q^{-2} t^4 y_{[2,5]}. 
\]
}
\eex
\subsection{Rescaling of generators}
\label{rescalegens}

Suppose $R$ is a CGL extension of length $N$ as in \eqref{itOre}. Given any nonzero scalars $t_1,\dots, t_N \in \kx$, one can rescale the generators $x_j$ of $R$ in the fashion
\begin{equation}
\label{rescalexj}
x_j \longmapsto t_j x_j, \; \; \forall\,  j \in [1,N],
\end{equation}
meaning that $R$ can be viewed as an iterated skew polynomial extension with generators $t_j x_j$:
\begin{equation}
\label{rescaleitOre}
R := \KK[t_1 x_1][t_2 x_2; \theta_2, t_2 \delta_2] \cdots [t_N x_N; \theta_N, t_N \delta_N].
\end{equation}
This is also a CGL extension presentation of $R$, and if \eqref{itOre} is a symmetric CGL extension, then so is \eqref{rescaleitOre}. 

A rescaling as in \eqref{rescaleitOre} does not affect the $\HH$-action or the matrix $\lab$, but various elements computed in terms of the new generators are correspondingly rescaled, such as the homogeneous prime elements from \thref{1} and \prref{yjk}, which transform according to the rules
\begin{equation}
\label{rescaleyelems}
y_k \longmapsto \biggl( \, \prod_{l = 0}^{O_-(k)} t_{p^l(k)} \biggr) y_k \quad \text{and} \quad y_{[i,s^m(i)]} \longmapsto \biggl( \, \prod_{l=0}^m t_{s^l(i)} \biggr) y_{[i,s^m(i)]} .
\end{equation}
Consequently, the homogeneous normal elements \eqref{uuu} transform as follows:
\begin{equation}
\label{rescaleuelems}
u_{[i,s^m(i)]} \longmapsto \bigl( t_i t_{s(i)}^2 \cdots t_{s^{m-1}(i)}^2 t_{s^m(i)} \bigr) u_{[i,s^m(i)]} .
\end{equation}

\subsection{Normalization conditions}
\label{normalizations}

We next describe some normalizations which are required in order for the homogeneous prime elements $y_k$ from \thref{1} to function as quantum cluster variables. Throughout this subsection, assume that $R$ is a symmetric CGL extension of length $N$ as in Definitions \ref{dCGL} and \ref{dsymCGL}. Assume also that the following mild conditions on scalars are satisfied:

{\bf{Condition (A).}} The base field $\KK$ contains square roots $\nu_{kl} = \sqrt{\la_{kl}}$ of the scalars
$\la_{kl}$ for $1 \leq l < k \leq N$ such that the subgroup of $\kx$ generated 
by all of them contains no elements of order $2$. Then set $\nu_{kk} := 1$ and $\nu_{kl} := \nu_{lk}^{-1}$ for $k < l$, so that $\nub := (\nu_{kl})$ is a multiplicatively skew-symmetric matrix.

{\bf{Condition (B).}} There exist positive integers $d_n$, $n \in \eta([1,N])$,
such that 
$$
\la_k^{d_{\eta(l)}} = \la_l^{d_{\eta(k)}} 
$$
for all $k, l \in [1,N]$ with $p(k), p(l) \neq - \infty$. In view of \prref{la}, this is equivalent to the condition
$$
(\la^*_k)^{d_{\eta(l)}} = (\la^*_l)^{d_{\eta(k)}}, \; \; \forall\,  k,l \in [1,N] \; \text{with} \; s(k), s(l) \ne +\infty.
$$

In parallel with \eqref{Scrlab}, define
\begin{equation}
\label{Scrnub}
\Scr_\nub(f) := \prod_{1\le j < k \le N} \nu_{jk}^{-m_j m_k}, \; \; \forall\,  f = (m_1,\dots, m_N)^T \in \Zset^N.
\end{equation}
Then set
\begin{equation}
\label{eybarj}
\ol{e}_j := e_j + e_{p(j)} + \cdots + e_{p^{O_-(j)}(j)} \quad \text{and} \quad  \ol{y}_j := \Scr_\nub(\ol{e}_j) y_j, \; \; \forall\,  j \in [1,N].
\end{equation}

The homogeneous prime elements appearing in \prref{yjk} are normalized analogously to \eqref{eybarj}:
\begin{equation}
\label{ybarismi}
\begin{aligned}
\ol{y}_{[i,s^m(i)]} := \; &\Scr_\nub( e_{[i,s^m(i)]} ) y_{[i,s^m(i)]}, \\ &\forall\,  i \in [1,N], \; m \in \Znn \; \; \text{such that} \; \; s^m(i) \in [1,N].
\end{aligned}
\end{equation}

One additional normalization, for the leading coefficients of the homogeneous normal elements $u_{[i,s^m(i)]}$, is needed in order to establish mutation formulas for the quantum cluster variables $\ol{y}_k$. For $i \in [1,N]$ and $m \in \Znn$ such that $s^m(i) \in [1,N]$, write
\begin{equation}
\label{ltu}
\begin{aligned}
\lt( u_{[i,s^m(i)]} ) =\;  &\pi_{[i,s^m(i)]} x^{f_{[i,s^m(i)]}},  \\
  &\pi_{[i,s^m(i)]} \in \kx, \; \; f_{[i,s^m(i)]} \in \sum_{j=1+1}^{s^m(i)-1} \Znn \, e_j \subset \Znn^N .
  \end{aligned}
\end{equation}
We will require the condition
\begin{equation}
\label{cond}
\pi_{[i,s(i)]} = \Scr_\nub( -e_i + f_{[i,s(i)]} ), \; \; \forall\,  i \in [1,N] \; \; \text{such that} \; \; s(i) \ne +\infty.
\end{equation} 
This can always be satisfied after a suitable rescaling of the generators of $R$, as follows.

\bpr{rescalecond}
{\rm\cite[Propositions 6.3, 6.1]{GYbig}}
Let $R$ be a symmetric CGL extension of length $N$, satisfying condition {\rm(A)}. 

{\rm (i)} There exist $N$-tuples $(t_1,\dots, t_N) \in (\kx)^N$ such that after the rescaling \eqref{rescalexj}, condition \eqref{cond} holds.

{\rm (ii)} The set of $N$-tuples occurring in {\rm (a)} is parametrized by $(\kx)^{\rk(R)}$.

{\rm (iii)} Once the $x_j$ have been rescaled so that \eqref{cond} holds, then
$$
\pi_{[i,s^m(i)]} = \;  \Scr_\nub( e_{[s(i),s^m(i)]} )^{-2} \Scr_\nub( -e_i + f_{[i,s^m(i)]} )
$$
for all $i \in [1,N]$, $m \in \Znn$ with $s^m(i) \in [1,N]$.
\epr

\bex{sl3d} {In the setting of \exref{sl3a} assume that $\sqrt{q} \in \KK$. 
The rescaling from \prref{rescalecond} of the generators of the symmetric CGL extension in the example, to obtain \eqref{cond}, is 
\[
X_k^\pm \mt q^{1/2} t^{-1} X_k^\pm \quad \mbox{for} \quad k =1,2, \quad X_{12}^- \mt  q^{1/2} t^{-1} X_{12}^-, \quad X_{12}^+ \mt - q^{3/2} t^{-1} X_{12}^+,
\]
using the notation \eqref{t-notation}. (The difference in the rescalings of $X_{12}^\pm$ symmetrizes the action of the automorphism $\ga$ described in \exref{sl3a}, which interchanges the rescaled elements $X_{12}^\pm$.) The corresponding normalized sequence of prime elements, expressed in terms of the original $y$-elements \eqref{ex-yelts}, is
\begin{equation}
\label{ex-yseq}
\begin{aligned}
\ol{y}_k = r_k y_k, \; \;  k \in [1,6], \quad \mbox{for} \quad
r_1 &= r_2 := q^{1/2} t^{-1}, \quad r_3 := q^{3/2} t^{-2},  \\
 r_4 &:= q^{5/2} t^{-3}, \quad r_5 := - q^3 t^{-2}, \quad r_6 := q^4 t^{-4}.
\end{aligned}
\end{equation}
}  
\eex
\subsection{Main theorem}
\label{thm8.2GYbig}

In this subsection, we present the main theorem from \cite{GYbig}, which gives quantum cluster algebra structures for symmetric CGL extensions. 

Recall the notation on quantum cluster algebras from \S \ref{q-cl}. There is a right action 
of $S_N$ on the set of toric frames for a division algebra $\FF$, given by re-enumeration,
\begin{equation}
\label{r-act}
(M \cdot \tau) (e_k) := M(e_{\tau(j)}), \; \;  \rb(M \cdot \tau)_{jk} = \rb(M)_{\tau(j), \tau(k)}, \quad \tau \in S_N, \;\; k \in [1,N].
\end{equation}

Fix a symmetric CGL extension $R$ of length $N$ such that Conditions (A) and (B) hold. Define the multiplicatively skew-symmetric matrix $\nub$ as in Condition (A),  with associated bicharacter $\Om_\nub$, and define a second multiplicatively skew-symmetric matrix $\rb = (r_{kj})$ by
\begin{equation}
\label{defrb}
 r_{kj} := \Om_\nub( \ol{e}_k, \ol{e}_j ), \; \; \forall\,  k,j \in [1,N].
\end{equation}
Let $\ol{y}_1, \ldots, \ol{y}_N$ be the sequence of normalized homogeneous prime elements 
given in \eqref{eybarj}. (We recall that generally 
each of those is a prime element of some of 
the subalgebras $R_l$, not of the full algebra $R=R_N$.) There is  a unique toric frame $M : \Zset^N \to \Fract(R)$ whose matrix is $\rb(M) := \rb$ 
and such that $M(e_k) := \ol{y}_k$, for all $k \in [1,N]$, \cite[Proposition 4.6]{GYbig}.

Next, consider an arbitrary element $\sig \in \Xi_N \subset S_N$, recall \eqref{tau}. By the definition of $\Xi_N$, 
for all $k \in [1,N]$,
\begin{equation}
\label{tau-int}
\eta^{-1} \eta \sig(k) \cap \sig([1,k])
= \begin{cases}
\{p^n(\sig(k)), \ldots, p(\sig(k)), \sig(k)\}, & \mbox{if} \; \; \sig(1) \le \sig(k)
\\ 
\{\sig(k), s(\sig(k)), \ldots, s^n(\sig(k))\}, & \mbox{if} \; \; \sig(1) \ge \sig(k)
\end{cases}
\end{equation}
for some $n \in \Znn$. Corresponding to $\sig$, we have the CGL extension presentation \eqref{tauOre}, whose $\lab$-matrix is  the matrix $\lab_\sig$ with entries $(\lab_\sig)_{ij} := \lab_{\sig(i) \sig(j)}$. Analogously we define the 
matrix $\nub_\sig$, and denote by $\rb_\sig$ the corresponding multiplicatively skew-symmetric 
matrix derived from $\nub_\sig$ by applying \eqref{defrb} to the presentation  \eqref{tauOre}. It is explicitly given by
\begin{equation}
\label{tau-frame1}
(\rb_\sig)_{kj} = \prod \{ \nu_{il} \mid i \in \sig([1,k]),\; \eta(i) = \eta \sig(k), \;
l \in \sig([1,j]), \; \eta(l) = \eta \sig (j) \},
\end{equation}
cf. \eqref{tau-int}.
Let $\ol{y}_{\sig,1}, \ldots, \ol{y}_{\sig,N}$ 
be the sequence of normalized prime elements given by \eqref{eybarj}
applied to the presentation \eqref{tauOre}. By \cite[Proposition 4.6]{GYbig}, 
there is a unique toric frame $M_\sig : \Zset^N \to \Fract(R)$ whose matrix 
is $\rb(M_\sig) := \rb_\sig$ and such that for all $k \in [1,N]$
\begin{equation}
\label{tau-frame2}
M_\sig(e_k) := \ol{y}_{\sig, k} =
\begin{cases} 
\ol{y}_{[p^n(\sig(k)), \sig(k)]}, 
\\ 
\ol{y}_{[\sig(k), s^n(\sig(k))]}
\end{cases}
\end{equation}
in the two cases of \eqref{tau-int}, respectively. The last equality 
is proved in \cite[Theorem 5.2]{GYbig}.

Recall that $P(N) := \{ k \in [1,N] \mid s(k) = + \infty \}$
parametrizes the set of homogeneous prime elements of $R$,
i.e.,  
$$
\{y_k \mid k \in P(N) \} \; \; \mbox{is a list of the homogeneous prime
elements of $R$}
$$ 
up to scalar multiples (\thref{1}). Define
$$
\ex := [1,N] \setminus P(N) = \{ l \in [1,N] \mid s(l) \neq + \infty \}.
$$
Since $|P(N)| = \rk(R)$,
the cardinality of this set is $|\ex| = N - \rk(R)$. 
More generally, for $\sig \in \Xi_N$, define the set
$$
\ex_\sig = \{ l \in [1, N] \mid \exists\, k > l \;\; \text{with} \;\; \eta \sig(k) = \eta \sig (l) \}
$$
of the some cardinality. Finally, recall 
that for a homogeneous element $u \in R$, $\chi_u \in \xh$ denotes 
its $\HH$-eigenvalue.

In \cite[Theorem 8.2]{GYbig} we re-indexed all toric frames $M_\sig$ in such a way that 
the right action in \thref{maint} (c) was trivialized and the exchangeable variables in all 
such seeds were parametrized just by $\ex$, rather than by $\ex_\sig$. We will not do the re-indexing here 
to simply the exposition.

\bth{maint} {\rm\cite[Theorem 8.2]{GYbig}}
Let $R$ be a symmetric CGL extension of length $N$ and rank $\rk(R)$
as in Definitions {\rm\ref{dCGL}, \ref{dsymCGL}}. Assume that Conditions {\rm (A), (B)} hold, and that the sequence of 
generators $x_1, \ldots, x_N$ of $R$ is normalized {\rm(}rescaled\/{\rm)}
so that condition \eqref{cond} is satisfied. 
Then the following hold:

{\rm(a)} For all $\sig \in \Xi_N$ {\rm(}see \eqref{tau}{\rm)} and $l \in \ex_\sig$, 
there exists a unique vector $b_\sig^l \in \Zset^N$ such that $\chi_{M_\sig(b_\sig^l)}=1$ and 
\begin{equation}
\label{linear-eq}
\Om_{\rb_\sig} ( b^l_\sig, e_j) = 1, \; \; \forall\,  j \in [1,N], \; j \neq l
\quad \mbox{and} \quad
\Om_{\rb_\sig} (b^l_\sig, e_l)^2 = \la^*_{ \min \eta^{-1} \eta (\sig(l))}.
\end{equation}
Denote by $\wt{B}_\sig \in M_{N \times |\ex| }(\Zset)$ the matrix with columns 
$b^l_\sig$, $l \in \ex_\sig$. Let $\wt{B}:= \wt{B}_\id$.
 
{\rm(b)} For all $\sig \in \Xi_N$, the pair $(M_\sig, \wt{B}_\sig)$ is a 
quantum seed for $\Fract(R)$. The principal part of $\wt{B}_\sig$ 
is skew-symmetrizable via the integers $d_{\eta(k)}$, $k \in \ex_\sig$ from Condition {\rm (B)}.

{\rm(c)} All such quantum seeds are mutation-equivalent to each other up to the $S_N$ action.
They are linked by the following one-step mutations.
Let $\sig, \sig' \in \Xi_N$ be such that
$$
\sig' = ( \sig(k), \sig(k+1)) \sig = \sig (k, k+1)
$$
for some $k \in [1,N-1]$.
If $\eta(\sig(k)) \neq \eta (\sig(k+1))$, then $M_{\sig'} = M_\sig \cdot (k,k+1)$
in terms of the action \eqref{r-act}.
If $\eta(\sig(k)) = \eta (\sig(k+1))$, 
then $M_{\sig'} = \mu_k (M_\sig)  \cdot (k,k+1)$.

{\rm(d)} We have the following equality between the CGL extension $R$ 
and the quantum cluster and upper cluster algebras associated to $M$, $\wt{B}$, $\varnothing$:
$$
R = \AA(M, \wt{B}, \varnothing) = \UU(M, \wt{B}, \varnothing).
$$
In particular, $\AA(M, \wt{B}, \varnothing)$ and $\UU(M, \wt{B}, \varnothing)$ 
are affine and noetherian, and more precisely $\AA(M, \wt{B}, \varnothing)$ is 
generated by the cluster variables in the seeds parametrized by the finite subset 
$\Ga_N$ of $\Xi_N$, recall \eqref{GaN}.

{\rm(e)} Let $\inv$ be any subset of the set $P(N)$ of frozen variables.
Then 
$$
R[y_k^{-1} \mid k \in \inv] = \AA(M, \wt{B}, \inv) = \UU(M, \wt{B}, \inv).
$$
\eth 
\bex{sl3e} 
{
Consider the symmetric CGL extension $R$ from \exref{sl3a}. 
}

{
(1) The initial seed for the quantum cluster algebra structure constructed on $R$ from \thref{maint} corresponding to $\sigma = \id$
has cluster variables $(\ol{y}_1, \ldots, \ol{y}_6)$
given by \eqref{ex-yseq}. Its exchange matrix is the one associated to the following mutation quiver with $5$ and $6$ being 
frozen vertices:
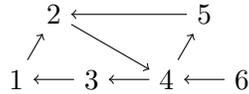
\[
\begin{tikzpicture}
\node (1) at (0,0) {1}; 
\node (3) at (1,0) {3}; 
\node (4) at (2,0) {4}; 
\node (6) at (3,0) {6};
\node (2) at (0.5,1.732/2) {2};
\node (5) at (2.5,1.732/2) {5};
\draw[<-] (1) -- (3);
\draw[<-] (3) -- (4);
\draw[<-] (4) -- (6);
\draw[<-] (2) -- (5);
\draw[<-] (2) -- (1);
\draw[<-] (4) -- (2);
\draw[<-] (5) -- (4);
\end{tikzpicture}
\]
}

{
(2) On the other hand, the quantum seed of $R$ from \thref{maint} corresponding to $\sigma = (123)$ will be linked to the 
seeds of Berenstein, Fomin and Zelevinsky \cite{BFZ,BZ}, see \exref{sl3g} below. This quantum seed has cluster variables
$(\ol{y}_2, q^{1/2} t^{-1} X_1^-, \ol{y}_3, \ldots, \ol{y}_6)$ and mutation quiver
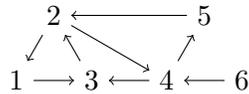
\[
\begin{tikzpicture}
\node (1) at (0,0) {1}; 
\node (3) at (1,0) {3}; 
\node (4) at (2,0) {4}; 
\node (6) at (3,0) {6};
\node (2) at (0.5,1.732/2) {2};
\node (5) at (2.5,1.732/2) {5};
\draw[->] (1) -- (3);
\draw[->] (4) -- (3);
\draw[->] (6) -- (4);
\draw[->] (5) -- (2);
\draw[->] (2) -- (1);
\draw[->] (3) -- (2);
\draw[->] (2) -- (4);
\draw[->] (4) -- (5);
\end{tikzpicture}
\]
where again 5 and 6 are frozen vertices.
This second seed is the mutation of the first seed at the vertex 1.
}
\eex

The next theorem gives an explicit formula 
for each of the exchange matrices $\wt{B}_\sig$ in terms of the original exchange matrix 
$\wt{B}:=\wt{B}_\id$. This implies that the only input needed to describe the quantum seeds $(M_\sig, \wt{B}_\sig)$ 
for a given algebra in \thref{maint} are the normalized interval prime elements (used in \eqref{tau-frame1}--\eqref{tau-frame2}) 
and the explicit form of the exchange matrix $\wt{B}$. Define the matrices
$$
Z= (z_{jk}), \; \; Z_\sig = (z_{jk}^\sig) \in M_N(\Zset), \quad \sig\in \Xi_N,
$$
by $z_{jk} = 1$ if $j = p^m(k)$ for some $m \in \Znn$ and $z_{jk} = 0$ otherwise, 
and $z_{jk}^\sig = 1$ if  $j  \in \sig([1,k]) \cap \eta^{-1}\eta (\sig(k))$  and
$z_{jk} = 0$ otherwise.

\bth{Btau} For each $\sig \in \Xi_N$, the columns of the exchange matrix $\wt{B}_\sig$ are 
explicitly given in terms of the columns of the original exchange matrix $\wt{B}$
as follows. For each $l \in \ex_\sig$, there exists $k >l$ such that $\eta \sig(k) = \eta \sig(l)$; 
set $k := \min \{ j > l \mid \eta \sig(j) = \eta \sig(l) \}$. Then 
either $\sig(k) = s^m \sig(l)$ or $\sig(k) = p^m \sig(l)$ for some 
$m \in \Zset_{>0}$, and 
\begin{equation}
\label{b-sig}
b_\sig^l =
\begin{cases} 
Z_\sig^{-1} Z (b^{\sig(l)} + \cdots + b^{s^{m-1} \sig(l)}),
\\
- Z_\sig^{-1} Z (b^{p^m \sig(l)} + \cdots + b^{p \sig(l)}), 
\end{cases}
\end{equation}
in the two cases respectively.
\eth

\begin{proof} For $j \in [1,N]$, set $\mu(j):= \min \eta^{-1} \eta (j)$ and
$$
\ol{e}_{\sig,j} := \sum \{ e_i \mid i \in \sig([1,j]) \cap \eta^{-1}\eta (\sig(j)) \} \in \Zset^N.
$$
Thus, in the two cases of \eqref{tau-int}, $\ol{e}_{\sig,j} =e_{[p^n(\sig(j)), \sig(j)]}$ and 
$\ol{e}_{\sig,j} = e_{[\sig(j), s^n(\sig(j))]}$, respectively (using the notation \eqref{e-interval}). In particular, $\ol{e}_j = \ol{e}_{\id,j}$. Define the homomorphism $\chi : \Zset^N \to \xh$ by $\chi( \sum_{j} m_j e_j ) : = \prod_j  \chi_{x_j}^{m_j}$. Thus, $\chi(\ol{e}_{\sig,j}) = \chi_{M_\sig(e_j)}$.
For $\sig \in \Xi_N$ and $l \in \ex_\sig$, define the vectors
$$
v_\sig^l := \sum_{j=1}^N (\wt{B}_\sig)_{jl} \ol{e}_{\sig, j} \in \Zset^N, \quad v^l:= v^l_\id.
$$
The sets $\{e_j\}$, $\{ \ol{e}_j \}$ and $\{\ol{e}_{\sig,j}\}$ are bases of $\Zset^N$ and $b_\sig^l = [v_\sig^l]_{\{\ol{e}_{\sig,j} \}}$
where $[v]_{ \{ \bullet\} }$ denotes the vector of coordinates of $v \in \Zset^N$ with respect to a given basis. 

By \cite[Lemma 8.14]{GYbig}, the statement of \thref{maint} (a) 
is valid if the first condition in \eqref{linear-eq} is replaced by 
$$
\Om_{\rb_\sig} ( b^l_\sig, e_j)^2 = (\la^*_{\mu(\sig(l))})^{\de_{jl}}, \quad \forall\, l \in \ex_\sig, \;\; j \in [1,N].
$$
Therefore, for $\sig \in \Xi_N$ and $l \in \ex_\sig$, the vector $v_\sig^l \in \Zset^N$ is the unique solution of 
\begin{equation}
\label{vsig}
\chi(v^l_\sig) = 1 \quad \mbox{and} \quad
\Om_{\nub} ( v^l_\sig, \ol{e}_{\sig,j})^2 = (\la^*_{\mu(\sig(l))})^{\de_{jl}}, \quad \forall\, j \in [1,N]. 
\end{equation}
In the case $\sig = \id$, the second condition is equivalent to 
$\Om_{\nub} ( v^l, e_j)^2 = (\la^*_{\mu(l)})^{a_{jl}}$ where $a_{jl} =1$ if $j=l$, $a_{jl} = -1$ if 
$j = s(l)$ and $a_{jl} = 0$ for $j \in [1,N] \backslash \{j, s(j) \}$. From the uniqueness statement in \eqref{vsig} it follows that
$$
v_\sig^l =
\begin{cases} 
v^{\sig(l)} + \cdots + v^{s^{m-1} \sig(l)},
\\
- v^{p^m \sig(l)} - \cdots - v^{p \sig(l)}, 
\end{cases}
$$
in the respective cases $\sig(k) = s^m\sig(l)$ and $\sig(k) = p^m\sig(l)$.
The theorem follows from this by taking into account that $Z$ and $Z_\sig$
are change of bases matrices 
$$
[v]_{ \{ e_j \} } = Z [v]_{ \{ \ol{e}_j \} }, \quad [v]_{ \{ e_j \} } = Z_\sig [v]_{ \{ \ol{e}_{\sig, j} \} }, \quad \forall\, v \in \Zset^N, \;\; j \in [1,N]
$$
and the fact that $b_\sig^l = [v_\sig^l]_{\{\ol{e}_{\sig,j} \}}$.
\end{proof}

\sectionnew{Quantum groups}
\label{qg}
In this section we gather some facts on quantum groups, quantized coordinate rings of 
double Bruhat cells and quantum Schubert cell algebras. We describe the connections between 
these classes of algebras, as well as their relation to the axiomatic 
class of quantum nilpotent algebras from the previous section.
\subsection{Quantized universal enveloping algebras}
\label{3.1}
Throughout $\g$ will denote an arbitrary finite dimensional complex simple Lie algebra with 
set of simple roots $\Pi= \{ \al_1, \ldots, \al_r \}$
and Weyl group $W$. 
Let $\{s_i\}$, $\{\al_i\spcheck \}$ and $\{\vpi_i\}$ be 
the sets of simple reflections, simple coroots and fundamental weights of $\g$.
The weight and root lattices of $\g$ will be denoted by $\PP$ and 
$\QQ$. Let $\PP^+$ be the set of dominant integral weights of $\g$
and $\QQ^+:= \sum_i \Zset_{\geq 0} \al_i$.
Denote by $\lcor.,. \rcor$ the invariant bilinear form on $\Rset \Pi$
normalized by $\lcor \al_i, \al_i \rcor = 2$ for short roots 
$\al_i$. The Cartan matrix of $\g$ is given by 
\begin{equation}
\label{Cartan}
(c_{ij}) := ( \lcor \al_i\spcheck, \al_j \rcor ) \in M_r(\Zset).
\end{equation}
For $\ga \in \Rset \Pi$, set
$$
\| \ga \|^2 = \lcor \ga, \ga \rcor.
$$
For a set of Weyl group elements $w_1, \ldots, w_m \in W$ denote the union of their supports
$$
\SS(w_1, \ldots, w_m) := \{ i \in [1,r] \mid s_i \; \; 
\mbox{appears in a reduced expression of some} \; \; w_j \}.
$$

For an arbitrary (infinite) field $\KK$ and a non-root of unity $q \in \KK^*$,
denote by $\UU_q(\g)$ the quantized universal enveloping algebra of $\g$ 
over $\KK$ with deformation parameter $q$. We will follow the conventions of 
Jantzen's book \cite{Ja} with the exception that we will denote
the standard generators of $\UU_q(\g)$ by $K_i^{\pm 1}$, $X_i^+$ and $X_i^-$ instead of 
$K_{\pm \al_i}$, $E_{\al_i}$ and $F_{\al_i}$ to avoid repeating formulas for the positive and negative parts of $\UU_q(\g)$, {and we will write $q_i$ for $q^{ \langle \al_i, \al_i \rangle/2 }$.}  
In particular we will use the form of the Hopf algebra $\UU_q(\g)$ with relations, 
coproduct, etc. as in \cite[\S 4.3, 4.8]{Ja}. 

The subalgebras of $\UU_q(\g)$ generated by $\{ X^\pm_i \}$ and $\{K_i\}$ 
will be denoted by $\UU^\pm$ and $\UU^0$, respectively.
The algebra $\UU_q(\g)$ is $\QQ$-graded by 
$$
\deg K_i = 0, \;\; \deg X_i^\pm = \pm \al_i,  \quad i \in [1,r].
$$
Its graded components will be denoted by $\UU_q(\g)_\ga$, $\ga \in \QQ$.

We will use the standard notation for $q$-factorials and binomial coefficients 
based on setting $[n]_q := (q^n - q^{-n})/(q - q^{-1})$ for $q$-integers. 
Let $q_i := q^{\| \al_i \|^2/2}$.

\subsection{Quantum groups}
\label{3.2}
A (left) $\UU_q(\g)$-module $V$ is called a \emph{type one} module if it is a direct sum 
of its weight spaces defined by 
$$
V_\mu := \{ v \in V \mid K_i v = q^{ \lcor \mu, \al_i \rcor} v, \; \; 
\forall\,  i \in [1,r] \}, \quad \mu \in \PP.
$$
The category of finite dimensional type one 
$\UU_q(\g)$-modules is semisimple 
(see  \cite[Theorem 5.17]{Ja} and the remark on p. 85 of 
\cite{Ja} for the case of general fields $\KK$) 
and is closed under taking tensor products and duals.
The irreducible modules in this category are parametrized by 
the dominant integral weights $\PP^+$, \cite[Theorem 5.10]{Ja}. 
For $\mu \in \PP^+$ we will denote by $V(\mu)$ the corresponding irreducible module
and we will fix a highest weight vector $v_\mu$ of $V(\mu)$.

Denote by $G$ the connected, simply connected, complex simple {algebraic} group with Lie algebra $\g$.
The quantum function algebra $R_q[G]$ is defined as the Hopf subalgebra 
of the restricted dual $\UU_q(\g)^\circ$ consisting of the matrix coefficients 
of all finite dimensional type one $\UU_q(\g)$-modules; that is, $R_q[G]$ is spanned by  
the matrix coefficients $c_{\xi, v}$ of the modules $V(\mu)$, 
$\mu \in \PP^+$, given by
\begin{equation} 
\label{c-notation}
c_{\xi, v}(x) := \xi ( x v ) \quad \mbox{for all} \; \; x \in \UU_q(\g)
\end{equation}
where $v \in V(\mu)$, $\xi \in V(\mu)^*$. The algebra $R_q[G]$ is $\PP \times \PP$-graded by 
\begin{equation}
\label{Rgrade}
R_q[G]_{\nu, \nu'} = \{ c_{\xi, v} \mid
\mu \in \PP^+, \; \xi \in (V(\mu)^*)_{\nu}, \; v \in V(\mu)_{\nu'} \}, 
\quad \nu, \nu' \in \PP.
\end{equation}

Throughout the paper, we will make use of the compatible Lusztig actions of the braid group $\BB_\g$ 
on the finite dimensional type one $\UU_q(\g)$-modules \cite[\S 8.6]{Ja} and 
on $\UU_q(\g)$ \cite[\S 8.14]{Ja}. For $\mu \in \PP^+$ and $w \in W$, set
\begin{equation}
\label{vwmu}
v_{w \mu} := T^{-1}_{w^{-1}} v_\mu \in V(\mu)_{w\mu}
\end{equation}
which by \cite[Proposition 39.3.7]{L} only depends on $w \mu$. Let 
$\xi_{w \mu} \in (V(\mu)^*)_{- w\mu}$ be the unique vector such that
$\lcor \xi_{w \mu}, T^{-1}_{w^{-1}} v_\mu \rcor =1$. 
For $w, u \in W$ and $\mu \in \PP^+$ define the (generalized) quantum minors
$$
\De_{w \mu, u \mu} := c_{\xi_{w\mu}, v_{u\mu}} \in R_q[G]_{-w\mu, u \mu}.
$$
They satisfy
\begin{equation}
\label{q-minor}
\De_{w \mu, u \mu} (x) = \De_{\mu, \mu} ( T_{w^{-1}} \cdot x \cdot T_{u^{-1}}^{-1} ), \quad \forall\,  x \in \UU_q(\g)
\end{equation}
and are equal \cite[\S 9.3]{GYbig} to the Berenstein--Zelevinsky quantum minors \cite[Eq. (9.10)]{BZ}.

\bre{T11} The use of $T_{w^{-1}}^{-1}$  (vs. $T_w$) 
is motivated by the fact that in the specialization 
$q=1$, the action of $T_{w^{-1}}^{-1}$ on $\UU_q(\g)$-modules 
specializes to the action of the Weyl group representative $\ol{w}$ of Fomin and Zelevinsky \cite{FZ}, 
while $T_w$ spezializes to $\ol{\ol{w}}$. 
\ere

The quantum minors satisfy
\begin{equation}
\label{prod}
\De_{w \mu, u \mu} \De_{w \mu', u \mu'} = \De_{w(\mu+\mu'), u(\mu + \mu')}
\quad \mbox{for} \; \; 
\mu, \mu' \in \PP^+. 
\end{equation}

\subsection{Quantized coordinate rings of double Bruhat cells}
\label{3.3}
Define the subalgebras of $R_q[G]$
\begin{align*}
R^+ &= 
\Span \{c_{\xi, v} \mid \mu \in \PP^+, \;
v \in V(\mu)_\mu, \; \xi \in V(\mu)^* \},
\\
R^- & = \Span \{ c_{\xi, v} \mid \mu \in \PP^+, \;
v \in V(\mu)_{w_0 \mu}, \; \xi \in V(\mu)^* \}
\end{align*}
where $w_0$ denotes the longest element of $W$. Joseph proved \cite[Proposition 9.2.2]{J} that 
$R_q[G]= R^+ R^- = R^- R^+$. A special case of the notation \eqref{vwmu} gives the lowest weight vectors 
\begin{equation}
\label{lowest}
v_{-\mu} = T^{-1}_{w_0} v_{- w_0 \mu} \in V(-w_0 \mu)_{-\mu}, \quad \mu \in \PP^+
\end{equation}
in  terms of which we have
$$
v_{-w \mu} = v_{w w_\ci (- w_\ci \mu)}= T_w v_{- \mu}
$$
(the point here is that $- w_\ci \mu \in \PP^+$). The latter enter in the definition of the quantum minors $\De_{-w \mu, -u \mu} \in R_q[G]$, $\mu \in \PP$ 
which are another way of writing those in \eqref{q-minor}.
For $\xi \in V(\mu)^*$ and $\xi' \in V(- w_0 \mu)^*$,
denote for brevity
\begin{equation}
\label{special-c}
c^+_\xi := c_{\xi, v_\mu } \; \; 
\mbox{and} \; \; 
c^-_{\xi'} := c_{\xi', v_{-\mu} }.
\end{equation}
Consider the Demazure modules
\begin{equation}
\label{Demazure}
V_w^+(\mu) = \UU^+ V(\mu)_{w \mu} \subseteq V(\mu) \quad \mbox{and} \quad 
V_u^-(\mu) = \UU^- V(- w_0 \mu)_{- u \mu} \subseteq V(- w_0 \mu).
\end{equation}
The ideals 
\begin{align*}
&I_w^\pm := \Span \{ c^\pm_\xi \mid \xi \in V_w^\pm(\mu)^\perp, \; \mu \in \PP^+ \} 
\subset R^\pm, \\
&I_{w, u} := I_w^+ R^- + R^+ I_u^- \subset R_q[G]
\end{align*}
are completely prime and homogeneous with respect to the $\PP \times \PP$-grading, \cite[Propositions 10.1.8, 10.3.5]{J}.
For all $\mu \in \PP^+$, $w \in W$,
$\nu, \nu' \in \PP$ and $a \in R_q[G]_{-\nu, \nu'}$
\begin{align}
\label{n1}
\De_{w \mu, \mu} a &= q^{\lcor w \mu, \nu \rcor - \lcor \mu, \nu' \rcor} 
a \De_{w \mu, \mu} \mod I_{w}^+ R^-,  
\\
\label{n2}
\De_{-w \mu, - \mu} a &= q^{ \lcor w \mu, \nu \rcor - \lcor \mu, \nu' \rcor} 
a \De_{-w \mu, - \mu} \mod R^+ I_w^-, 
\end{align}
see \cite[(2.22), (2.23)]{Y-sqg}. 
{Because of these commutation relations, these quantum minors will play the role of frozen variables later.}
Denote $E^\pm_w := \{ \De_{\pm w \mu, \pm \mu} \mid \mu \in \PP^+ \}$ and
$E_{w,u} := q^\Zset E_w^+ E_u^- \subset R_q[G]$.
By abuse of notation we will denote by the same symbols the images 
of elements and subsets in factor rings. The set $E_{w,u}$ is a multiplicative 
subset of normal elements of $R_q[G]/I_{w,u}$ and 
$E^\pm_w$ are multiplicative subsets of normal elements 
of $R^\pm/I_w^\pm$ (by \eqref{prod}, \eqref{n1}-\eqref{n2}). In particular,
\begin{equation}
\label{empty}
E_{w,u} \cap I_{w,u} = \varnothing \quad 
\mbox{and} \quad 
E_w^\pm \cap I_w^\pm = \varnothing.
\end{equation}
The localization 
$$
R_q[G^{w,u}] := (R_q[G]/I_{w,u})[E_{w,u}^{-1}]
$$
is called the \emph{quantized coordinate ring of the double Bruhat cell $G^{w,u}$}.

\subsection{The algebras $S_w^\pm$ and a presentation of the quantum double Bruhat cell algebras}
\label{3.4}
The localizations $R^\pm_w := (R^\pm/I^\pm_w)[(E^\pm_w)^{-1}]$ are $\PP \times \PP$-graded.
Their subalgebras consisting of elements in degrees $\PP \times \{ 0 \}$ will be denoted 
by $S_w^\pm$. A simple computation shows that they are effectively $\QQ \cong \QQ \times \{ 0 \}$--graded.
The graded components of $S_w^\pm$ will be denoted by $(S_w^\pm)_\ga= (S_w^\pm)_{\ga, 0}$, 
$\ga \in \QQ$. Each element of $S^+_w$ (resp. $S^-_u$) can be represented 
in the form $c^+_\xi \De_{w \mu, \mu}^{-1}$ 
(resp. $\De_{- w \mu, - \mu}^{-1} c^-_{\xi'}$)
for some $\mu \in \PP^+$ and $\xi \in V(\mu)^*$ 
(resp. $\xi' \in V(- w_0 \mu)^*$), see \cite[\S 10.3.1]{J}.
Here and below the coset notations for factor rings are omitted. From the form of the elements of $S_w^\pm$ we see that $(S_w^\pm)_0 = \KK$.

For $\ga \in \QQ^+$ set
$m_w(\ga) := \dim \UU^+_\ga= \dim \UU^-_{-\ga}$. Let
$\{x_{\ga, n} \}_{n=1}^{m_w(\ga)}$ and 
$\{x_{-\ga, n} \}_{n=1}^{m_w(\ga)}$ be 
dual bases of $\UU^+_\ga$ and 
$\UU^-_{-\ga}$ with respect to the Rosso--Tanisaki form,
cf. \cite[Ch. 6]{Ja}. 

Joseph proved \cite[\S 10.3.2]{J} that there exists a unique $\QQ \cong \QQ \times \{ 0 \}$-graded 
$\KK$-algebra $S^+_w \bowtie S^-_u$ 
which satisfies the following properties:
\begin{enumerate}
\item[(i)] $S^+_w \bowtie S^-_u \cong S^+_w \otimes_\KK S^-_u$ as $\KK$-vector spaces and
the canonical embeddings of $S^+_w$ and $S^-_u$ in it are graded algebra embeddings,
\item[(ii)] the elements of $S^+_w$ and $S^-_u$ satisfy the commutation relations
\begin{align}
\label{commRR}
&\big[ \De_{- u \mu', -\mu'}^{-1} c^-_{\xi'}  \big] 
\big[ c^+_\xi \De_{w \mu, \mu}^{-1} \big]
= q^{ - \lcor \nu' + u \mu' , \nu - w \mu  \rcor} 
\big[ c^+_\xi \De_{w \mu, \mu}^{-1} \big]
\big[ \De_{ - u \mu', - \mu'}^{-1} c^-_{\xi'}  \big]
\\
\nn
+ & \sum_{\ga \in \QQ^+, \ga \neq 0}
\sum_{i=1}^{m(\ga)}
q^{ - \lcor \nu' + \ga + u \mu', \nu - \ga - w \mu \rcor } 
\big[ c^+_{S^{-1}(x_{\ga, i}) \xi} \De_{w \mu, \mu}^{-1} \big] .
\big[ \De_{-u \mu', - \mu'}^{-1} c^-_{S^{-1}(x_{-\ga, i}) \xi'} \big],
\end{align}
for all $\mu, \mu' \in \PP^+$, $\xi \in( V(\mu)^*)_{-\nu}$, $\xi' \in (V(- w_0 \mu')^*)_{- \nu'}$.
\end{enumerate}
The commutation relation 
\eqref{commRR} is similar to a bicrossed product of Hopf algebras.
The multiplication in the algebra induces the linear isomorphism
\begin{equation}
\label{SwSu2}
S^+_w \bowtie S^-_u \cong S^-_u \otimes_\KK S^+_w.
\end{equation}

Furthermore, Joseph proved \cite[\S 10.3.2]{J} that there exists a graded algebra embedding 
\begin{equation}
\label{SRembed}
S^+_w \bowtie S^-_u \hra R_q[G^{w,u}]
\end{equation}
whose restrictions to the subalgebras $S^+_w$ and $S^-_u$ are the compositions 
\begin{align*}
&S^+_w \hra (R^+/I_w^+)[(E^+_w)^{-1}] \hra (R_q[G]/I_{w,u})[E_{w,u}^{-1}]=R_q[G^{w,u}] 
\quad \mbox{and} \quad \\
&S^-_u \hra (R^-/I_u^-)[(E^-_u)^{-1}] \hra (R_q[G]/I_{w,u})[E_{w,u}^{-1}]=R_q[G^{w,u}].
\end{align*}
We will identify $S^+_w \bowtie S^-_u$ with its image.
Consider the identification of $\UU_q(\g)$-modules
\begin{equation}
\label{VVidentify}
V(\vpi_i)^* \cong V(- w_0 \vpi_i)
\quad {\mbox{normalized by}} \quad 
\xi_{\vpi_i} \mt v_{- \vpi_i}.
\end{equation}
Let $\{\xi_{i, j} \}$ and $\{\xi'_{i,j} \}$
be a pair of dual bases of $V(\vpi_i)^* \cong V(- w_0 \vpi_i)$ and $V(- w_0 \vpi_i)^*$ with respect to the corresponding dual pairing.
Define the elements 
\begin{equation}
\label{p_idef}
p_i:= \sum_j \big[ \De_{- u \vpi_i, - \vpi_i}^{-1} c^-_{\xi'_{i,j}}  \big]
\big[ c^+_{\xi_{i,j}} \De_{w \vpi_i, \vpi_i}^{-1} \big]  \in (S^+_w \bowtie S^-_u)_{(w-u)\vpi_i}.
\end{equation}
The element $p_i$ is not a scalar if and only if $i \in \SS(w,u)$, \cite[Proposition 7.6]{Y-sqg}. 
We have $\sum_j c^-_{\xi'_{i,j}} c^+_{\xi_{i,j}} =1$ because
\begin{align}
\label{cc1}
\sum_j c^-_{\xi'_{i,j}}(x) c^+_{\xi_{i,j}} (y) &= \sum_j \lcor \xi'_{i,j},  x  v_{- \vpi_i} \rcor \lcor \xi_{i,j}, y v_{\vpi_i} \rcor 
\\
&= \lcor x  \xi_{\vpi_i}, y v_{\vpi_i} \rcor = \lcor \xi_{\vpi_i}, S(x) y v_{\vpi_i} \rcor,  \quad 
\forall\,  x, y \in \UU_q(\g).
\nn
\end{align}
Therefore, in the embedding \eqref{SRembed},
\begin{equation}
\label{psi-pi}
p_i = \De_{- u \vpi_i, - \vpi_i}^{-1} \De_{w \vpi_i, \vpi_i}^{-1}.
\end{equation}
By \eqref{n1}--\eqref{n2}, all $p_i$ are normal elements of $S^+_w \bowtie S^-_u$. Extend \eqref{SRembed} to the embedding
\begin{equation}
\label{SRembed2}
(S^+_w \bowtie S^-_u)[p_i^{-1}, i \in \SS(w,u)]  \hra R_q[G^{w,u}].
\end{equation}

Denote by $\LL_w^+$ and $\LL_u^-$ the respective subalgebras of $R_q[G^{w,u}]$ generated by the elements of $E_w^+ \cup (E_w^+)^{-1}$ 
and the elements of $E_u^- \cup (E_u^-)^{-1}$.
It follows from \eqref{prod} and \eqref{empty} that $\LL_w^+$ and $\LL_u^-$ are Laurent polynomial rings
\begin{equation}
\label{Laur}
\LL_w^+, \; \LL_u^- \cong \KK [z_1^{\pm1}, \ldots, z_r^{\pm1}] \quad 
\mbox{by} \quad 
\De_{w \vpi_i, \vpi_i}^{ \pm 1}, \;
\De_{-u \vpi_i, - \vpi_i}^{ \pm 1} \mt z_i^{\pm 1}. 
\end{equation}
Finally, Joseph proved \cite[10.3.2(5)]{J} that 
\begin{align}
\label{dBisom}
R_q[G^{w,u}] &\cong
\big( (S^+_w \bowtie S^-_u)[p_i^{-1}, i \in \SS(w,u)] \big) \# \LL_w^+
\\
&\cong
\big( (S^+_w \bowtie S^-_u)[p_i^{-1}, i \in \SS(w,u)] \big) \# \LL_u^-
\nn
\end{align}
where the smash products are defined via the commutation relations \eqref{n1}-\eqref{n2}. 
Since $\De_{w \mu, \mu}$ is a nonzero element of $R_q[G^{w,u}]_{-w \mu, \mu}$ for $\mu \in \PP^+$, 
\eqref{SRembed2} induces the isomorphism
\begin{equation}
\label{redDB}
(S^+_w \bowtie S^-_u) [p_i^{-1}, i \in \SS(w,u)] \stackrel{\cong}{\longrightarrow} \bigoplus_{\nu \in \QQ} R_q[G^{w,u}]_{\nu, 0} = 
R_q[G^{w,u}/H].
\end{equation}
The algebra on the right is called the \emph{quantized coordinate ring of 
the reduced double Bruhat cell $G^{w,u}/H$} where $H = B^+ \cap B^-$.

Some of the quoted results of Joseph were stated in \cite{J} for base fields 
of characteristic 0 and $q \in \KK$ which is transcendental over $\Qset$. 
However all hold in full generality \cite{Y-multi,Y-sqg}, see 
in particular \cite[\S 3.4]{Y-sqg} for the above set of embeddings and isomorphisms.   

\subsection{Quantum Schubert cell algebras}
\label{3.5}
We will identify the rational character lattice $\xh$ of the torus 
$\HH:= ( \KK^*)^r$ with $\QQ$ by mapping $\ga \in \QQ$ to the character
\begin{equation}
\label{Hchar}
h=(t_1, \ldots, t_r) \in \HH \longmapsto h^\ga = \prod t_i^{ \lcor \ga, \vpi_i\spcheck \rcor}.
\end{equation}
A $\QQ$-graded algebra $R$ (e.g., $\UU_q(\g)$) has a canonical 
rational $\HH$-action by algebra automorphisms where $h \cdot x = h^\ga x$ for $x \in R_\ga$, 
$h \in \HH$. Define the embedding $\PP\spcheck \hra \HH$ by
taking $\nu \in \PP\spcheck$ to the unique element $h(\nu) \in \HH$ such that 
$h(\nu)^\ga = q^{\lcor \ga, \nu \rcor}$ for all $\ga \in \QQ$. 

Consider a Weyl group element $w$ and a reduced expression 
$w = s_{i_1} \ldots s_{i_N}$. Define the roots
\begin{equation}
\label{beta}
\beta_k := s_{i_1} \ldots s_{i_{k-1}} (\al_{i_k}), \quad \forall\,  k \in [1, N].
\end{equation}
The \emph{quantum Schubert cell algebras} $\UU^\pm[w]$, $w \in W$,
were defined by De Concini, Kac, and Procesi \cite{DKP}, and Lusztig \cite[\S 40.2]{L} 
as the subalgebras of $\UU^\pm$ generated by 
\begin{equation}
X_{\be_k}^\pm := T_{i_1} \ldots T_{i_{k-1}} (X_{i_k}^\pm), 
\quad k \in [1,N]
\label{rootv}
\end{equation}
see \cite[\S 39.3]{L}. These algebras do not depend on the choice of a 
reduced expression for $w$. They are $\QQ$-graded subalgebras 
of $\UU_q(\g)$ and are thus stable under the $\HH$-action. 

A reduced expression of $w$ gives rise to presentations of the algebras $\UU^\pm[w]$ as 
symmetric CGL extensions of the form 
\begin{equation}
\label{UwCGL1}
\UU^\pm[w] = \KK [X_{\be_1}^\pm] [X_{\be_2}^\pm; (h(\mp \be_2) \cdot), \de_2)] \ldots
[X_{\be_N}^\pm; (h(\mp \be_N) \cdot), \de_N]
\end{equation}
\cite[Lemma 9.1]{GYbig}. The reverse presentations have the forms
\begin{equation}
\label{UwCGL2}
\UU^\pm[w] = \KK [X_{\be_N}^\pm] [X_{\be_{N-1}}^\pm; (h(\pm \be_{N-1}) \cdot), \de^\sy_{N-1})] \ldots
[X_{\be_1}^\pm; (h(\pm \be_1) \cdot), \de^\sy_1].
\end{equation}
This is derived from the Levendorskii--Soibelman 
straightening law 
\begin{equation}
\label{LS}
X_{\be_k}^\pm X_{\be_j}^\pm - 
q^{ - \lcor \be_k, \be_j \rcor }
X_{\be_j}^\pm X_{\be_k}^\pm  \in \KK \lcor X_{\be_{j+1}}^\pm, \ldots, X_{\be_{k-1}}^\pm \rcor,
\quad \forall \; 1 \le j < k \le N.
\end{equation} 
The interval subalgebras for the algebras $\UU^\pm[w]$ are given by 
\begin{equation}
\label{intU}
\UU^\pm[w]_{[j,k]}= T_{i_1} \ldots T_{i_{j-1}}(\UU^\pm[s_{i_j} \ldots s_{i_k}])
\end{equation}
for all $1 \leq j \leq k \leq N$.

There is an algebra isomorphism $\om : \UU^+[w] \to \UU^-[w]$ 
obtained by restricting the algebra automorphism $\omega$ of $\UU_q(\g)$ given by 
$$
\om(K_i) = K_i^{-1}, \;\;
\om(X_i^\pm) = X_i^{\mp}, \quad
\forall\,  i \in [1,r].
$$
It is graded if the grading of one of the algebras $\UU^\pm[w]$ is reversed (i.e., $\ga \mt - \ga$).

\subsection{Quantum function algebras and isomorphisms}
\label{3.6}
We will make extensive use of the fact that 
the quantum Schubert cell algebras $\UU^\pm[w]$ are (anti)isomorphic to the algebras 
$S^\mp_w$. To describe these (anti)isomorphisms we will need the quantum $R$-matrix $\RR^w$ 
corresponding to $w$. It equals the sum of the tensor products of the elements in 
any two dual bases of $\UU^+[w]$ and $\UU^-[w]$ with respect to the 
Rosso-Tanisaki form and is explicitly given by 
\begin{multline}
\label{Rw}
\RR^w := \sum_{m_1, \ldots, m_N \, \in \, \Zset_{\geq 0}}
\left( \prod_{j=1}^N
\frac{ (q_{i_j}^{-1} - q_{i_j})^{m_j}}
{q_{i_j}^{m_j (m_j-1)/2} [m_j]_{q_{i_j}}! } \right) 
\\
(X^+_{\be_N})^{m_N} \ldots (X^+_{\be_1})^{m_1} \otimes  
(X^-_{\be_N})^{m_N} \ldots (X^-_{\be_1})^{m_1}.
\end{multline}
It belongs to the completion $\UU^+[w] \, \wh{\otimes} \, \UU^-[w]$ of 
$\UU^+[w] \otimes \UU^-[w]$ with respect to the descending 
filtration \cite[\S 4.1.1]{L}. We will denote the flip of the two components of $\RR^w$ by
\begin{equation}
\label{Rw2}
\RR^w_{\opp} \quad \mbox{and also set} \quad
\RR := \RR^{w_\ci}, \;\; \RR_{\opp} = \RR^{w_\ci}_{\opp}. 
\end{equation}
We will also need the 
unique graded algebra antiautomorphism 
$\tau$ of $\UU_q(\g)$ given by 
\begin{equation}
\label{Uq-tau}
\tau(X^\pm_i) = X^\pm_i,
\; \;
\tau(K_i) = K_i^{-1}, \quad
\forall\,  i \in [1,r],
\end{equation}
see \cite[Lemma 4.6(b)]{Ja}.
It is compatible with the braid group action:
\begin{equation}
\label{iota-ident}
\tau (T_w x) = T_{w^{-1}}^{-1} ( \tau (x)), \quad
\forall\,  x \in \UU_q(\g), \; w \in W,
\end{equation}
see \cite[Eq. 8.18(6)]{Ja}. For a linear map $\theta : \UU_q(\g) \to \UU_q(\g)$, denote
\begin{equation}
\label{dp}
{}^{\theta \times 1}\RR:= (\theta \otimes \id)\RR. 
\end{equation}
Analogously, we define ${}^{\theta \times 1}\RR^w$, ${}^{\theta \times 1}\RR_{\opp}$,
${}^{1 \times \theta}\RR$, etc.

\bth{isom} \cite[Theorem 2.6]{Y-sqg}
For all finite dimensional simple Lie algebras $\g$ and Weyl group elements $w \in W$,
the maps $\vp^\pm_w : S^\pm_w \to \UU^\mp[w]$ given by
\begin{align*}
& \vp^+_w \big( c^+_\xi \De_{w\mu, \mu}^{-1} \big) 
:= \big( c_{\xi, v_{w \mu}} \otimes \id \big)
{}^{\tau \times 1} \RR^w, \; \; 
\forall\,  \mu \in \PP^+, \; \xi \in V(\mu)^*
\\
& \vp^-_w \big( \De_{-w\mu, -\mu}^{-1} c^-_\xi  \big) 
:= \big( \id \otimes c_{\xi, v_{- w \mu}} \big)
{}^{1 \times \tau} \RR^w_{\opp}, \; \; 
\forall\,  \mu \in \PP^+, \; \xi \in V(- w_0\mu)^*
\end{align*}
are well defined $\QQ$-graded algebra antiisomorphisms and isomorphisms, 
respectively. In the definition of $\vp_w^\pm$ one can replace $\RR^w$ with $\RR$.
\eth

\sectionnew{Quantum double Bruhat cells and symmetric CGL extensions}
\label{qdoubB}
In this section we prove that all algebras $S^+_w \bowtie S^-_u$ are 
symmetric CGL extensions. Their interval subalgebras are identified with 
algebras of the same kind. We also describe an explicit model for the isomorphic 
bicrossed products $\UU^-[w]_{\opp} \bowtie \UU^+[u]$.
\subsection{Statements of the results}
\label{4.1}
Throughout the section we fix reduced expressions of the 
Weyl group elements $w, u \in W$,
\begin{equation}
\label{reduced}
w = s_{i_1} \ldots s_{i_N} \quad \mbox{and} \quad
u = s_{i'_1} \ldots s_{i'_M}.
\end{equation}
Consider the root vectors \eqref{beta}, \eqref{rootv} and 
\begin{equation}
\label{rootv2}
\be'_k := s_{i'_1} \ldots s_{i'_{k-1}}(\al_{i'_k}), \quad X^{\pm}_{\be'_k} := T_{i'_1} \cdots T_{i'_{k-1}} (X^{\pm}_{i'_k}), \quad
k \in [1,M].
\end{equation}
Define 
$$
w_{\leq k} := s_{i_1} \ldots s_{i_k} 
$$
and
$$
w_{[j,k]}
:= 
\begin{cases} 
s_{i_j} \ldots s_{i_k}, & \mbox{if} \; \; j \leq k
\\
1, & \mbox{if} \; \; j>k.
\end{cases}
$$
(For the simplicity of notation we will not explicitly show the dependence of this notation on the  
choice of reduced decomposition.) Define in a similar fashion $u_{\leq k}$ and $u_{[j,k]}$. 

Let 
\begin{equation}
\label{Celem}
C_{\be_k}^+:= \left( \vp^+_w \right)^{-1}(X^-_{\be_k}) \in (S_w^+)_{-\be_k} \quad \mbox{and} \quad
C_{\be'_l}^-:= \left( \vp^-_u \right)^{-1}(X^+_{\be'_l}) \in (S_u^-)_{\be'_l}
\end{equation}
for $k \in [1,N]$, $l \in [1,M]$.

\bth{SwSu} Let $\KK$ be an arbitrary base field and $q \in \kx$ a non-root of unity.
For all finite dimensional simple Lie algebras $\g$ and pairs of Weyl group elements $(w,u)$, the algebra 
$S^+_w \bowtie S^-_u$ is a symmetric CGL extension. It has the 
presentation
\begin{multline}
\label{SwSu1}
S^+_w \bowtie S^-_u = \KK [C_{\be_N}^+] [C_{\be_{N-1}}^+, (h(\be_{N-1}) \cdot), \de_{N-1}] \ldots 
[C_{\be_1}^+, (h(\be_1) \cdot), \de_1] 
\\
[C_{\be'_1}^-, (h(-\be'_1) \cdot), \del_1] 
\ldots [C_{\be'_M}^-, (h(-\be'_M) \cdot), \del_M]
\end{multline}
and the reverse presentation
\begin{multline}
\label{SwSu1b}
S^+_w \bowtie S^-_u = \KK [C_{\be'_M}^-] [C_{\be'_{M-1}}^-, (h(\be'_{M-1}) \cdot), \del^\sy_{M-1}] \ldots 
[C_{\be'_1}^-, (h(\be'_1) \cdot), \del^\sy_1] 
\\
[C_{\be_1}^+, (h(-\be_1) \cdot), \de^\sy_1] 
\ldots [C_{\be_N}^+, (h(-\be_N) \cdot), \de^\sy_N]
\end{multline}
where $\de_k, \de^\sy_k, \del_l , \del^\sy_l$ are locally nilpotent skew derivations.
\eth
Denote for brevity
\begin{equation}
\label{abs}
|k| := 
\begin{cases}
N-k+1, & \mbox{if} \; \; k \in [1,N]
\\
k-N, & \mbox{if} \; \; k \in [N+1, N+M].
\end{cases}
\end{equation}

\bco{la-k} The scalars $\la_k, \la_k^\sy$ {\rm{(}}from Definition {\rm\ref{dCGL} {(iii))}} for the symmetric CGL extension $S^+_w \bowtie S^-_u$ are 
given by 
\begin{equation}  \label{la*bowtie}
\la_k^\sy = \la_k^{-1} = 
\begin{cases}
q_{i_{ |k| }}^2, & \mbox{if} \; \; 1 \leq k \leq N
\\
q_{i'_{ |k| }}^2, & \mbox{if} \; \; N < k \leq M+N. 
\end{cases}
\end{equation}
The multiplicatively skew-symmetric matrix $\lab = (\la_{kj})$ for this CGL extension is given by
\begin{equation}  \label{lakjbowtie}
\la_{kj} = \begin{cases}
q^{ - \langle \be_{|k|}, \be_{|j|} \rangle }, &\text{if} \; \; 1 \le j < k \le N  \\
q^{ - \langle \be'_{|k|}, \be'_{|j|} \rangle }, &\text{if} \; \; N < j < k \le N+M  \\
q^{ \langle \be'_{|k|}, \be_{|j|} \rangle }, & \text{if} \; \; 1 \le j \le N < k \le N+M.
\end{cases}
\end{equation}
\eco

\begin{proof} By \thref{SwSu} 
$$
\la_k C^+_{\be_{|k|}}= h(\be_{|k|}) \cdot C^+_{\be_{|k|}} = 
q^{- \| \be_{ |k| } \|^2 } C^+_{\be_{|k|}} = 
q^{- \| \al_{i_{|k|}} \|^2 } C^+_{\be_{|k|}} = q_{i_{|k|}}^{-2} C^+_{\be_{|k|}}
$$ 
for $1 \leq k \leq N$. The other formulas and cases of \eqref{la*bowtie} are analogous. The first two cases of \eqref{lakjbowtie} follow from \eqref{LS} and \eqref{Celem}, while the third one follows from the upcoming \prref{int}.
\end{proof} 

\bth{SwSu2} In the setting of Theorem {\rm\ref{tSwSu}} the interval subalgebras of 
$S^+_w \bowtie S^-_u$ satisfy
\begin{align*}
&\big( S^+_w \bowtie S^-_u \big)_{[j,k]} \cong \UU^-[w_{[|k|, |j|]}]^{\opp} \cong S^+_{w_{[|k|,|j|]}}  \quad 
\mbox{for} \quad 1 \leq j \leq k \leq N,  
\\
&\big( S^+_w \bowtie S^-_u \big)_{[j,k]} \cong \UU^+[u_{[|j|,|k|]}] \cong S^-_{u_{[|j|,|k|]}}   \quad 
\mbox{for} \quad N < j \leq k \leq N + M,  
\\
&\big( S^+_w \bowtie S^-_u \big)_{[j,k]} \cong S^+_{w_{\leq |j| }} \bowtie S^-_{u_{\leq |k| }} \quad 
\mbox{for} \quad 1 \leq j \leq N < k \leq N + M.
\end{align*}
\eth

The first two sets of isomorphisms are obtained from the ones 
in \thref{isom}. For example, for $1 \leq j \leq k \leq N$, the isomorphism 
$$
\big( S^+_w \bowtie S^-_u \big)_{[j,k]} \cong \UU^-[w_{[N-k+1, N-j +1]}]^{\opp}
$$
is the restriction of $\vp^+_w$ to $\KK \lcor C^-_{N-k+1}, \ldots, C^-_{N-j+1} \rcor$. The isomorphism
$$
\UU^-[w_{[N-k+1, N-j +1]}]^{\opp} \cong S^+_{w_{[N-k+1,N-j+1]}}
$$
is given by $(\vp^-_{w_{[N-k+1, N-j +1]}})^{-1}$. The last isomorphism
is constructed in \S \ref{4.3}.

\subsection{Proof of \thref{SwSu}}
\label{4.2}
\thref{isom} and Eq. \eqref{LS} imply that 
\begin{align}
& C^+_{\be_k} C^+_{\be_j} - q^{- \lcor \be_k, \be_j \rcor} C^+_{\be_j} C^+_{\be_k} \in 
\KK \lcor C^+_{\be_{j+1}}, \ldots, C^+_{\be_{k-1}} \rcor \quad 
\mbox{for} \; \; 1 \leq j < k \leq N, 
\label{LS1}
\\
& C^-_{\be'_k} C^-_{\be'_j} - q^{- \lcor \be'_k, \be'_j \rcor} C^-_{\be_j} C^-_{\be_k} \in 
\KK \lcor C^-_{\be'_{j+1}}, \ldots, C^-_{\be'_{k-1}} \rcor \quad 
\mbox{for} \; \; 1 \leq j < k \leq M.
\label{LS2}
\end{align}
\thref{SwSu} easily follows by induction on $j=N, \ldots, 1$, $k=N+1, \ldots, N+M$
from the following proposition, these two identities, and the definition 
of the elements $h(\ga)$, $\ga \in \QQ\spcheck$. We leave the details to the reader.

\bpr{int} In the setting of Theorem {\rm\ref{tSwSu}}, for all $j \in [1,N]$ and $k \in [1,M]$,
$$
C^-_{\be'_k} C^+_{\be_j} - q^{ \lcor \be'_k, \be_j \rcor } C^+_{\be_j} C^-_{\be'_k} 
\in \KK \lcor C^+_{\be_{j-1}}, \ldots, C^+_{\be_1}\rcor . \KK \lcor C^-_{\be'_1}, \ldots, C^-_{\be'_{k-1}} \rcor.
$$
\epr

Before we proceed with the proof of the proposition we obtain some 
auxiliary results. Extend the reduced expression of $w$ 
to a reduced expression of the longest element 
of $W$
$$
w_0 = s_{i_1} \ldots s_{i_N} \ldots s_{i_l}.
$$
Define the roots $\be_k$ and the root vectors $X^\pm_{\be_k}$, $k \in [1,l]$ 
by extending \eqref{beta} and \eqref{rootv} to $k\in [1,l]$. 
For $\mb=(m_1, \ldots, m_l) \in \Znn^l$ set 
$$
X^{\mb, \pm} = (X^\pm_{\be_l})^{m_l} \ldots (X^\pm_{\be_1})^{m_1}.
$$
The Levendorskii--Soibelman straightening law \eqref{LS} and the fact that $\UU^\pm$ are 
$\QQ^\pm$-graded imply that 
\begin{equation}
\label{geq}
X^{\mb', \pm} X^{\mb, \pm} 
\in \Span \{ X^{\mb'', \pm} \mid 
\mb'' \in \Znn^l, \; \mb' \preceq \mb'' \},
\quad
\forall\,  \mb, \mb' \in \Znn^l
\end{equation}
with respect to the reverse lexicographic order \eqref{prec} on $\Znn^l$.

\ble{Xact} If $j \in [1,N]$, $\mu \in \PP^+$, and $\xi \in V(\mu)^*$ are such that 
$$
C^+_{\be_j} = c^+_\xi \De_{w \mu, \mu}^{-1}, 
$$
then
$$
\vp^+_w(c^+_{S^{-1}(\tau (X^{\mb,+})) \xi} \De_{w \mu, \mu}^{-1} )
\in \UU^-[w_{\leq j-1}], 
\quad
\forall\,  \mb \in \Znn^l, \; \mb \neq 0.
$$
\ele

\begin{proof} The last statement in \thref{isom} implies that
\begin{equation}
\label{assume}
\lcor \xi, \tau (X^{\mb',+}) v_{w \mu} \rcor = 0 \quad
\forall \mb' \in \Znn^l, \; \; \mb' \neq e_j
\end{equation} 
where $e_1, \ldots, e_l$ is the standard basis of $\Znn^l$.
Assume that $\mb, \mb'=(m'_1, \ldots, m'_l) \in \Znn^l$ are such that
$\mb \neq 0$ and $m'_j + \cdots + m'_N > 0$. 
By \eqref{geq},
$$
X^{\mb', +} X^{\mb, +} 
\in \Span \{ X^{\mb'', +} \mid 
\mb'' \in \Znn^l, \; \mb' \prec \mb'' \}.
$$
Combining this with \eqref{assume} leads to 
$$
\lcor S^{-1}(\tau (X^{\mb,+})) \xi, \tau (X^{\mb',+}) v_{w \mu} \rcor = 
\lcor \xi, \tau (X^{\mb', +} X^{\mb,+}) v_{ w \mu} \rcor =0.
$$
This implies the statement of the lemma in view of the definition of the map $\vp^+_w$.
\end{proof}

\begin{proof}[Proof of Proposition {\rm\ref{pint}}] 
The fact that $S^w_+ \bowtie S^u_- \cong S^w_+ \otimes S^u_-$ 
as $\KK$-vector spaces, \thref{isom} and the iterated skew polynomial extension 
presentations of $\UU^\pm[w]$ from \S \ref{3.5} imply 
\begin{multline}
\KK \lcor C^+_{\be_{j-1}}, \ldots, C^+_{\be_1} \rcor . S^-_u 
\cap S^+_w . \KK \lcor C^-_{\be'_1}, \ldots, C^-_{\be'_{k-1}} \rcor \\
= \KK \lcor C^+_{\be_{j-1}}, \ldots, C^+_{\be_1}\rcor . \KK \lcor C^-_{\be'_1}, \ldots, C^-_{\be'_{k-1}} \rcor.
\label{inter}
\end{multline}

We will prove that 
\begin{equation}
C^-_{\be'_k} C^+_{\be_j} - q^{ \lcor \be'_k, \be_j \rcor } C^+_{\be_j} C^-_{\be'_k} 
\in \KK \lcor C^+_{\be_{j-1}}, \ldots, C^+_{\be_1} \rcor . S^-_u
\label{C+-1}
\end{equation}
Analogously one proves 
\begin{equation}
C^-_{\be'_k} C^+_{\be_j} - q^{ \lcor \be'_k, \be_j \rcor } C^+_{\be_j} C^-_{\be'_k} 
\in S^+_w . \KK \lcor C^-_{\be'_1}, \ldots, C^-_{\be'_{k-1}} \rcor.
\label{C+-2}
\end{equation}
The proposition follows from the combination of Eqs. \eqref{inter}, \eqref{C+-1} and \eqref{C+-2}.

Let $\mu, \mu' \in \PP^+$, $\xi \in V(\mu)_{-\nu}$, $\xi' \in V(- w_0 \mu')_{- \nu'}$ 
be such that 
$$
C^+_{\be_j} = c^+_\xi \De_{w \mu, \mu}^{-1} \quad \mbox{and} \quad
C^-_{\be'_k} = \De_{- u \mu', - \mu'}^{-1} c^-_{\xi'}.
$$
Since $\vp_w^+$, $\vp_u^-$ are graded maps, 
\begin{equation}
\label{degg}
-\be_j = - \nu + w \mu \quad \text{and} \quad \be'_k = - \nu' - u \mu'.  
\end{equation}
The sets 
\begin{equation}
\label{bases0}
\{ X^{\mb, +} \mid \mb \in \Znn^l, \; \mb \neq 0 \} \quad 
\mbox{and} \quad 
\{ \tt_{\mb} X^{\mb, -} \mid \mb \in \Znn^l, \; \mb \neq 0 \}
\end{equation}
are dual bases of 
$\oplus_{\ga \in \QQ^+, \ga \neq 0 } \, \UU^+_\ga$ 
and $\oplus_{\ga \in \QQ^+, \ga \neq 0 } \, \UU^-_{-\ga}$ 
with respect to the Rosso--Tanisaki form for some scalars $\tt_\mb \in \kx$ (equal to those in \eqref{Rw}). 
By \cite[Lemma 6.16]{Ja}, 
the antiautomorphism $\tau$ of $\UU_q(\g)$ satisfies 
$$
\lcor \tau(x_+), \tau(x_-) \rcor = \lcor x_+, x_- \rcor, \quad 
$$
for all $x_\pm \in \UU^\pm$. So applying $\tau$ to \eqref{bases0} produces another set 
of dual bases. Invoking the commutation relation \eqref{commRR} with the use of the last pair of dual bases,
\leref{Xact} and Eq. \eqref{degg} gives
$$
C^-_{\be'_k} C^+_{\be_j} - q^{ \lcor \be'_k, \be_j \rcor } C^+_{\be_j} C^-_{\be'_k} 
\in \KK \lcor C^+_{\be_{j-1}}, \ldots, C^+_{\be_1} \rcor S^u_-.
$$
This proves \eqref{C+-1} which completes the proof of the proposition.
\end{proof}

\subsection{Proof of \thref{SwSu2}}
\label{4.3}
The first two sets of isomorphisms in \thref{SwSu2} were defined and proved 
in \S 4.1. We proceed with the last isomorphism. Fix two integers $j \in [1,N]$ and 
$k \in [N+1, N+M]$. Set for brevity 
$$
A^+:= \KK \lcor C^+_{\be_{N-j+1}}, \ldots, C^+_{\be_1} \rcor, 
\quad 
A^-:= \KK \lcor C^-_{\be'_1}, \ldots, C^-_{\be'_{k-N}} \rcor.
$$
Identify 
$$
(S_w^+ \bowtie S_u^-)_{[j,k]} = 
\KK \lcor  C^+_{\be_{N-j+1}}, \ldots, C^+_{\be_1}, 
C^-_{\be'_1}, \ldots, C^-_{\be'_{k-N}} \rcor \cong A^+ \otimes_\KK A^-
$$
as $\KK$-vector spaces. Define the $\KK$-linear isomorphism 
$$
\pi_{[j,k]} : 
(S_w^+ \bowtie S_u^-)_{[j,k]} \to 
S_{w_{\leq N-j+1}}^+ \bowtie S_{u_{\leq k-N}}^-
$$
by 
$$
\pi_{[j,k]}|_{A^+} := \big( \vp_{w_{\leq N-j+1}}^+ \big)^{-1} \vp_w^+ |_{A^+}
\quad \mbox{and} \quad 
\pi_{[j,k]}|_{A^-} := \big( \vp_{u_{\leq k-N}}^- \big)^{-1} \vp_u^-.
$$
It provides the last isomorphism in \thref{SwSu2} by the following result.

\bth{SwSu3} The map
$$
\pi_{[j,k]} : (S_w^+ \bowtie S_u^-)_{[j,k]} \to
S_{w_{\leq N-j+1}}^+ \bowtie S_{u_{\leq k-N}}^-
$$ 
is a $\QQ$-graded algebra isomorphism
for all $1 \leq j \leq N < k \leq N+M$.
\eth

\begin{proof} By \thref{isom}, the restrictions of $\pi_{[j,k]}$ to $A^\pm$
are $\QQ$-graded algebra isomorphisms. The statement of the theorem amounts 
to verifying that the commutation rule 
\eqref{commRR} is transformed appropriately under $\pi_{[j,k]}$.

It is sufficient to establish the theorem in the two cases
$j\in [1,N]$, $k=N+M$ and $j=1$, $k \in [N+1,N+M]$, and apply a composition of the maps $\pi$ in the general case.
We will restrict to the first case, the second being analogous. 
Thus, from now on $k = N+M$.

Fix $m \in [1,N-j +1]$.
Let $\mu \in \PP^+$, $\xi \in V(\mu)^*$, $\ol{\xi} \in V(\mu)^*$ 
be such that 
$C^+_{\be_m} = c^+_\xi \De_{w \mu, \mu}^{-1}$ 
and
$
\pi_{[j, N+M]}(C^+_{\be_m}) = 
\big( \vp_{w_{\leq N-j + 1}}^+ \big)^{-1} (X^-_{\be_m}) 
= c^+_{\ol{\xi}} \De_{w_{\leq N-j +1} \mu, \mu}^{-1}.
$
Consider the $\KK$-basis $\{X^{\mb, +} \mid \mb \in \Znn^l \}$ of $\UU^+$
from \S 4.2. As in \S 4.2, the definition of $\vp^+_w$ implies that 
$$
\lcor \xi, \tau(X^\mb) v_{w \mu} \rcor = 
\lcor \ol{\xi}, \tau(X^\mb) v_{w_{\leq N-k+1} (\mu)} \rcor = \de_{\mb, e_m}
$$
for all $\mb \in \Znn^l$. Therefore 
$$
\lcor \xi, x v_{w \mu} \rcor = 
\lcor \ol{\xi}, x V_{ w_{\leq N-k+1} (\mu)} \rcor, 
\quad \forall\,  x \in \UU^+
$$ 
and
$$
\vp_w^+ \left( c^+_{S^{-1}(x) \xi} \De_{w \mu, \mu}^{-1} \right)= 
\vp_{w_{\leq N-j+1}}^+ \left( c^+_{S^{-1}(x) \ol{\xi} } \, \De_{w_{\leq N-j +1} \mu, \mu}^{-1} \right),
\quad \forall\,  x \in \UU^+.
$$
This implies that $\pi_{[j,N+M]}$ transforms appropriately the commutation relation \eqref{commRR}
in the case when the first term is taken to be an arbitrary element of $S_u^-$ and the second term is 
taken to be $C^+_{\be_m}$, $m \in [1,N-j+1]$. Hence, $\pi_{[j,N+M]}$ is a $\QQ$-graded algebra isomorphism.
\end{proof}

\subsection{The algebras $\UU^-[w]_{\opp} \bowtie \UU^+[u]$} We use the linear isomorphism 
$\vp_w^+ \otimes \vp_u^-$ from $S^+_w \bowtie S^-_u \cong  S^+_w \otimes S^-_u$ to $\UU^-[w]_{\opp} \otimes \UU^+[u]$ to transfer 
the algebra structure on $S^+_w \bowtie S^-_u$ to an algebra structure on $\UU^-[w]_{\opp} \otimes \UU^+[u]$ which will 
be denoted by $\UU^-[w]_{\opp} \bowtie \UU^-[u]$. \thref{isom} implies that 
\begin{equation}
\label{vp-wu}
\vp_{w,u}:= \vp_w^+ \otimes \vp_u^- : S^+_w \bowtie S^-_u \cong S^+_w \otimes S^-_u \to \UU^-[w]_{\opp} \bowtie \UU^+[u]
\end{equation}
is an algebra isomorphism and that $\UU^-[w]_{\opp}, \UU^+[u] \hra \UU^-[w]_{\opp} \bowtie \UU^+[u]$ are algebra embeddings.
By way of definition $\vp_{w,u}|_{S^+_w} = \vp_w^+$ and $\vp_{w,u}|_{S^-_u} = \vp_u^-$. By \thref{SwSu3} we have the 
algebra embeddings
$$
\UU^-[w]_{\opp} \bowtie \UU^+[u] \hra \UU^-_{\opp} \bowtie \UU^+ = \UU^-[w_0]_{\opp} \bowtie \UU^+[w_0]
$$
induced by the canonical embeddings $\UU^-[w] \hra \UU^-$ and $\UU^+[u] \hra \UU^+$. 
This gives the following explicit description of the algebras $\UU^-[w]_{\opp} \bowtie \UU^+[u] \cong S^+_w \bowtie S^-_u $.

\bth{explicitUwu} {\rm{(i)}} For any simple Lie algebra $\g$, the algebra $\UU^-_{\opp} \bowtie \UU^+$
is isomorphic to the tensor product $\UU^-_{\opp} \otimes \UU^+$ in which $\UU^-_{\opp}$ and $\UU^+$ sit as subalgebras 
and commute as follows:
\begin{equation}
\label{X+-}
X_i^+ X_j^- = q^{\lcor \al_i, \al_j \rcor} X_j^- X_i^+ + \de_{i,j} (q_i^{-1} - q_i)^{-1}, \quad
\forall\,  i, j \in [1,r].
\end{equation}
In other words, $\UU^-_{\opp} \bowtie \UU^+$ is the algebra with generators $X^\pm_i$ subject to the quantum Serre relations for $X^+_i$, the opposite quantum
Serre relations for $X^-_i$ and the above mixed relations. 

{\rm{(ii)}} For all Weyl group elements $w$, $u$, the algebra $\UU^-[w]_{\opp} \bowtie \UU^+[u]$ is a subalgebra of $\UU^-_{\opp} \bowtie \UU^+$ equal to the subspace $\UU^-[w]_{\opp} \otimes \UU^+[u]$. In other words, it is the subalgebra generated by $X^-_{\be_1}, \ldots X^-_{\be_N}$, 
$X^+_{\be'_1}, \ldots, X^+_{\be'_M}$ in the notation of \S {\rm\ref{4.1}}.
\eth

\bre{iso-to-thers}
The scalars in the right hand side of \eqref{X+-} can be rescaled to any collection in $\kx$ by rescaling $X^\pm_1, \ldots, X^\pm_r$, which does
not change the Serre relations. The map $\tau \otimes \id :  \UU^-_{\opp} \bowtie \UU^+ \to \UU^- \bowtie \UU^+$ is an algebra isomorphism, 
where the target algebra is defined as the factor algebra of $\UU^- \otimes \UU^+$ by the relations \eqref{X+-}. 

{Because of these facts, the algebra $\UU^-_{\opp} \bowtie \UU^+$ is isomorphic to Kashiwara's bosonic algebra ${\mathscr{B}}_q(\g)$, \cite[\S 3.3]{Ka}, 
which plays a key role in his construction of global bases. In the simply laced case, the algebra $\UU^-_{\opp} \bowtie \UU^+$
is a subalgebra of the Hernandez--Leclerc algebra \cite[Theorem 7.3]{HL}. 
The algebra $\UU^-_{\opp} \bowtie \UU^+$ is also a subalgebra of the Heisenberg double of 
the algebra $\UU_q(\b_+)$, and plays a role in Bridgeland's realization \cite[Lemma 4.5]{B} of quantized universal enveloping algebras 
via Hall algebras  and the Berenstein--Greenstein work on double canonical bases \cite[Eq. (1.2)]{BG}. }
\ere
\begin{proof}[Proof of \thref{explicitUwu}] For given $i, j \in [1,r]$ choose two reduced expression of $w_0$ that start with $s_i$ and $s_j$ respectively. Then in the 
notation from \S \ref{4.1} for $w = u = w_0$, $X^-_{\be_1} = X^-_j$ and $X^+_{\be'_1} =  X^+_i$. \thref{SwSu} implies that 
$$
X_i^+ X_j^- = q^{\lcor \al_i, \al_j \rcor} X_j^- X_i^+ + t_{i,j} 
$$
for some $t_{i,j} \in \KK$. For degree reasons $t_{i,j} =0$ for $i \neq j$. The scalars $t_{i,i}$ are 
evaluated directly by applying \eqref{commRR}. The second part of the theorem follows from \thref{SwSu3}.
\end{proof}

Using \thref{isom} and a direct calculation, one shows that the elements $p_i$ from \eqref{p_idef} are given by
\begin{equation}
\label{pp1}
p_i = \tt \vp_{w,u}^{-1} \left( \big( \De_{u \vpi_i, w \vpi_i } \otimes \id \big)  
\left( {}^{S \tau \times 1} \RR^u_{\opp} \, {}^{\tau \times 1} \RR^w
\right)
\right)  
\end{equation}
for some $t \in \kx$, recall the notation \eqref{Rw2} and \eqref{dp}. \thref{explicitUwu} 
gives an explicit model for the quantized coordinate rings of all reduced double Bruhat cells (cf. \S \ref{3.4}),
$$
R_q[G^{w,u}/H] \cong (\UU^-[w]_{\opp} \bowtie \UU^+[u])\big[ \big((\De_{u \vpi_i, w \vpi_i } \otimes \id)
({}^{S \tau \times 1} \RR^u_{\opp} \, {}^{\tau \times 1} \RR^w)\big)^{-1}, i \in \SS(w,u)\big].
$$
\bre{comm} Assume that $\sqrt{q} \in \KK$. The following rescaling of the generators of the algebras in \thref{explicitUwu} will 
play an important role in the construction of toric frames, see \eqref{x-new}:
$$
x_i^{\pm} = q_i^{1/2}(q_i^{-1} - q_i) X_i^\pm.
$$
In terms of those, the commutation relations \eqref{X+-} simplify to
$$
x_i^+ x_i^- = q_i^2 x_i^- x_i^+ + (1 - q_i^2), \quad
x_i^+ x_j^- = q^{\lcor \al_i, \al_j \rcor} x_j^- x_i^+, \quad
i \neq j \in [1,r].
$$
\ere
\sectionnew{Classification of the homogeneous prime elements of $S_w^+ \bowtie S_u^-$}
\label{prime}
In this section we classify the homogeneous prime elements of 
the algebras $S_w^+ \bowtie S_u^-$. This determines 
the ranks of these algebras as CGL extensions. As other corollaries of the main result we obtain 
classifications of all normal elements of the algebras $S_w^+ \bowtie S_u^-$ and $R_q[G^{w,u}]$.

\subsection{Statements of the main result and its corollaries}
\label{5.1}
As in the previous sections, all of the following results are for an arbitrary base field 
$\KK$ and a non-root of unity $q \in \kx$.

Recall the definition of the normal elements $p_i \in (S_w^+ \bowtie S_u^-)_{(w-u)\vpi_i}$
from \S \ref{3.4}. 

\bth{primes} For all finite dimensional simple Lie algebras $\g$ and $w, u \in W$, the homogeneous prime elements 
of $S_w^+ \bowtie S_u^-$ are the nonzero scalar multiples of the elements 
$p_i \in (S_w^+ \bowtie S_u^-)_{(w-u)\vpi_i}$, $i \in \SS(w,u)$, defined in \eqref{p_idef},
and those elements are not associates of each other.
In particular, the rank of the algebra $S_w^+ \bowtie S_u^-$ equals $|\SS(w,u)|$.
\eth

The theorem is proved in \S \ref{5.4} and \S\S \ref{5.2}--\ref{5.3} contain preparatory material.

The elements $p_i$ quasicommute and the monomials in those elements are homogeneous normal elements of 
$S_w^+ \bowtie S_u^-$. By \cite[Theorem 1.2]{GY1} for every CGL extension $R$, the monomials in 
the homogeneous prime elements of $R$ with different collections of 
exponents have different weights with respect to the $\xh$-grading.
The normalizing automorphisms of the $p_i$ are determined from 
\begin{equation}
\label{p-norm}
p_i a = q^{ - \lcor (w+u)\vpi_i, \ga \rcor} a p_i, \quad
\forall\,  a \in (S_w^+ \bowtie S_u^-)_\ga, \; \ga \in \QQ
\end{equation}
which follows from \eqref{n1}--\eqref{n2} and \eqref{psi-pi}.

\bco{1} The homogeneous normal elements of $S_w^+ \bowtie S_u^-$ are precisely 
the monomials in the elements $p_i$, $i \in \SS(w,u)$.

The normal elements of $S_w^+ \bowtie S_u^-$ are the sums
of monomials in the $p_i$ that have the same normalizing automorphism
{\rm(}which is explicitly computed from \eqref{p-norm}{\rm)}.
\eco

The corollary follows from \thref{primes} and \cite[Proposition 2.6]{GY1} which 
states that the normal elements of an $\HH$-UFD $R$ (for a torus $\HH$) 
are all sums of monomials in the homogeneous prime elements of $R$ that 
have the same normalizing automorphism. 

The elements $\De_{w\vpi_i, \vpi_i}$, $i \in \SS(w,u)$, and 
$\De_{-u \vpi_i, - \vpi_i}$, $i \in [1,r]$, are normal elements 
of $R_q[G^{w,u}]$ and their normalizing automorphisms 
are given by \eqref{n1}--\eqref{n2}. By the above discussion and Eq. \eqref{psi-pi} 
the Laurent polynomials in those elements with different sets of exponents
are linearly independent. 
On the other hand, for all $i \notin \SS(w,u)$, 
$\De_{w\vpi_i, \vpi_i} = \De_{-u\vpi_i, -\vpi_i}^{-1}$ 
(as elements of $R_q[G^{w,u}]$).

\bco{2} The normal elements of $R_q[G^{w,u}]$ are the sums of Laurent 
monomials in $\De_{w\vpi_i, \vpi_i}$, $i \in \SS(w,u)$, and 
$\De_{-u \vpi_i, - \vpi_i}$, $i \in [1,r]$, that have the same normalizing
automorphism. 

The homogeneous normal elements of $R_q[G^{w,u}]$ with respect to the
$\PP \times \PP$-grading are the Laurent 
monomials in $\De_{w\vpi_i, \vpi_i}$, $i \in \SS(w,u)$, and 
$\De_{-u \vpi_i, - \vpi_i}$, $i \in [1,r]$
\eco

\begin{proof} 
Each element $x$ of $R_q[G^{w,u}]$ has the form $\sum_j x_j \# m_j$ where $x_j \in S^+_w \bowtie S^-_u$ 
and $m_j$ are linearly independent Laurent monomials in $\De_{-u \vpi_i, - \vpi_i}$, $i \in [1,N]$.
Such an element is homogeneous with respect to the $\PP \times \PP$-grading 
if and only if the sum has only one term and $x_1$ is a homogeneous 
element of $S^+_w \bowtie S^-_u$. Under those conditions, if $x$ 
is a normal element of $R_q[G^{w,u}]$, then $x_1$ is a 
normal element of $S^+_w \bowtie S^-_u$.
Applying \coref{1} and Eq. \eqref{psi-pi} completes the proof of the second part.

The first part follows from the fact that the normal elements of a $\Zset^l$-graded ring $R$ 
are the linear combinations of homogeneous normal elements of $R$ that have the same 
normalizing automorphism.
\end{proof}

\coref{2} strengthens \cite[Theorem 1.1]{Y-sqg} which classified the centers 
of $R_q[G^{w,u}]$ by investigating certain sets of $q$-normal elements of $R_q[G^{w,u}]$.
The results in Corollaries \ref{c1} and \ref{c2} are of independent interest due 
to the relationship of those sets of normal elements to prime spectra of quantum groups and 
related algebras \cite{J}.

\subsection{Factors of the elements $p_i$}
\label{5.2}
\ble{primes} If $p$ is a homogeneous prime element of $S_w^+ \bowtie S_u^-$, 
then $p | p_i$ for some $i \in \SS(w,u)$.
\ele
\begin{proof} Assume that this is not the case. It follows from \eqref{dBisom} that
$p$ is a homogeneous prime element of $R_q[G^{w,u}] \cong (S_w^+ \bowtie S_u^-) \# \LL_u$.
This contradicts the fact that $R_q[G^{w,u}]$ has no nontrivial homogeneous prime ideals \cite[Theorem 10.3.4]{J}. 
\end{proof} 

\subsection{Maximal terms of elements of $S_w^+ \bowtie S_u^-$}
\label{5.3}
The $\QQ$-grading of the algebras $S_w^+$ and $S_u^-$ is supported in $- \QQ^+$ and $\QQ^+$, respectively. 
Consider the vector space decomposition 
$$
S^+_w \bowtie S^-_u = \bigoplus_{\nu \in - \QQ^+, \, \nu' \in \QQ^+} (S^+_w)_\nu \otimes (S^-_u)_{\nu'}.
$$
For an element $a \in S^+_w \bowtie S^-_u$ denote by $[a]_{\nu, \nu'}$ its component 
in $(S^+_w)_{\nu} \otimes (S^-_u)_{\nu'}$. Denote its support by
$$
\supp (a) := \{ (\nu, \nu') \in (-\QQ^+) \times \QQ^+ \mid [a]_{\nu, \nu'} \neq 0 \}.
$$
Define the partial order $\geq$ on $(-\QQ^+) \times \QQ^+$ by 
\begin{multline*}
(\nu_1, \nu'_1) \geq (\nu_2, \nu'_2) \quad 
\mbox{if and only if there exist $\ga, \ga' \in \QQ^+$}
\\
\mbox{such that} \; \; \nu_1 = \nu_2 -\ga \;\; \text{and} \;\; \nu'_1 = \nu'_2 + \ga'.
\end{multline*}
Denote by $\msupp(a)$ the maximal elements in the support of an element 
$a \in S^+_w \bowtie S^-_u$. For $a \in S_w^+ \bowtie S_u^-$, $a \neq 0$, 
$$
a \in \kx \Leftrightarrow \msupp(a) = \{ (0,0) \}.
$$
Call the components $[a]_{\nu, \nu'}$ for $(\nu, \nu') \in \msupp(a)$
the \emph{maximal terms} of $a$. By \eqref{commRR},
\begin{equation}
\label{prod-max}
\msupp(xy) = \{ \text{maximal elements of} \; \msupp(x) + \msupp(y) \}, \quad \forall x,y \in S^+_w \bowtie S^-_u .
\end{equation}

For a subset $I \subseteq [1,r]$, set 
$$
\PP^+_I = \oplus_{i \in I} \Znn \vpi_i.
$$
Define the elements
$$
d^\pm_{w, \mu} := \De_{\pm \mu, \pm \mu} \De_{\pm w \mu, \pm \mu}^{-1} \in (S^\pm_w)_{\pm(w-1)\mu}, 
\; \; \forall \mu \in \PP^+_{\SS(w)}.
$$
(For the rest of this section the previous convention for right/left denominators for $S^\pm_w$ 
will not be needed.)
Their scalar multiples exhaust all homogeneous normal elements of $S^\pm_w$ 
by \cite[Theorem 6.1 (i)]{Y-sqg} and 
$$
d^\pm_{w, \mu} a = q^{ \lcor (w+1) \mu ,  \ga \rcor } a d^\pm_w, \quad
\forall\,  a \in (S^\pm_w)_\ga, \; \ga \in \QQ,
$$ 
cf. \cite[Eq. (3.30)]{Y-sqg}.

\ble{maxterm} The maximal terms of each homogeneous normal element of $S^+_w \bowtie S^-_u$ 
are of the form
$$
\tt d^+_{w, \mu} d^-_{u,\mu'}
$$
for some $\mu \in \PP^+_{\SS(w)}$, $\mu' \in \PP^+_{\SS(u)}$ and $\tt \in \kx$.
\ele
\begin{proof} Each maximal term of a normal element $c \in S^+_w \bowtie S^-_u$ 
can be written in the form 
$$
a_1 b_1 + \cdots + a_m b_m 
$$
where $a_1, \ldots, a_m \in (S^+_w)_\nu$ and $b_1, \ldots, b_m$ are linearly independent 
elements of $(S^-_u)_{\nu'}$ for some $\nu \in - \QQ^+$, $\nu' \in \QQ^+$. The normality 
of $c$, linear independence of $b_1, \ldots, b_m$ and Eq. \eqref{commRR}
easily imply that $a_1, \ldots, a_m$ 
are normal elements. By the above mentioned \cite[Theorem 6.1 (i)]{Y-sqg} 
all of them must be scalar multiples of $d^+_{w, \mu}$ where 
$\mu \in \PP^+_{\SS(w)}$ is the unique element such that $\nu = (w-1)\mu$
(if $\nu $ can be written in this form).
Reversing the roles of $S^+_w$ and $S^-_u$ leads to the desired result.
\end{proof}

\subsection{Proof of the main result}
\label{5.4}
The definition of the element $p_i$ implies that it has one maximal term which equals 
$$
\ol{p}_i := 
\begin{cases}
\tt d^+_{w, \vpi_i} d^-_{u,\vpi_i}, 
& \mbox{if} \; \; i \in \SS(w) \cap \SS(u)
\\
\tt d^+_{w, \vpi_i}, 
& \mbox{if} \; \;  i \in \SS(w) \backslash \SS(u)
\\
\tt d^-_{u, \vpi_i}, 
& \mbox{if} \; \;  
i \in \SS(u) \backslash \SS(w)
\end{cases}
$$
for some $\tt \in \kx$.

\begin{proof}[Proof of Theorem {\rm\ref{tprimes}}]
Because of \leref{primes} and 
\cite[Proposition 2.2 (b)]{GY1}, all we need to prove is that
$p_i$ are irreducible elements of $S_w^+ \bowtie S_u^-$ for $i \in \SS(w,u)$.
Let $p_i = a b$ for some $a, b \in S_w^+ \bowtie S_u^-$. The fact that 
$p_i$ has only one maximal term and Eqs. \eqref{commRR} and \eqref{prod-max} imply that 
$a$ and $b$ also have one maximal term and those maximal terms (to be denoted by $\ol{a}, \ol{b}$)
satisfy 
$$
\ol{a} \ol{b} = \ol{p}_i.
$$ 

Case 1: $i \in \SS(w) \backslash \SS(u)$. The tensor product decomposition of $S_w^+ \bowtie S_u^-$ 
into $S_w^+$ and $S_u^-$, and the fact that $d_{w,\vpi_i}^+$ are prime elements of $S_w^+$ 
(\cite[Theorem 6.1 (ii)]{Y-sqg})
imply that one of the factors $\ol{a}, \ol{b}$ is a scalar. This means that either 
$a \in \kx$ or $b \in \kx$, so $p_i$ is irreducible.
 
Case 2: $i \in \SS(u) \backslash \SS(w)$. This case is analogous to the previous one.

Case 3: $i \in \SS(w) \cap \SS(u)$. In this case 
$$
\ol{a} \ol{b} = \tt d^+_{w, \vpi_i} d^-_{u,\vpi_i}, \quad \tt \in \kx
$$
If $\ol{a}$ or $\ol{b}$ is a scalar, we complete the proof as in case 1. 
Otherwise, the primality of $d_{w, \vpi_i}^+ \in S_w^+$ and 
$d_{u, \vpi_i}^- \in S_u^-$ and \eqref{commRR} imply that either 
$$
\ol{a} = \tt_1 d_{w, \vpi_i}^+, \; \ol{b} = \tt_2 d_{u, \vpi_i}^- \quad
\mbox{or} \quad
\ol{a} = \tt_1 d_{u, \vpi_i}^-, \; \ol{b} = \tt_2 d_{w, \vpi_i}^+
$$
for some $\tt_1, \tt_2 \in \kx$.
We will restrict to the first case, leaving the analogous second case to the reader.
The definition of maximal support implies that $a \in S_w^+$ and $b \in S_u^-$. Furthermore the 
tensor product decomposition of $S_w^+ \bowtie S_u^-$ into $S_w^+$ and $S_u^-$ 
and the homogeneity of $p_i$ imply that $a$ and $b$ also need to be homogeneous. Thus
$a= \ol{a} = \tt_1 d_{w, \vpi_i}^+$, $b= \ol{b} = \tt_2 d_{u, \vpi_i}^-$ and
\begin{equation}
\label{p-dd}
p_i = \tt d_{w, \vpi_i}^+ d_{u, \vpi_i}^-.
\end{equation}
Since $i \in \SS(w) \cap \SS(u)$, by \cite[Corollary 4.4.6]{J}, $c^+_{\xi_{s_i\varpi_i}} = \De_{s_i \vpi_i, \vpi_i} \notin I^+_w$
and $c^-_{\xi_{-s_i\varpi_i}} = \De_{- s_i \vpi_i, - \vpi_i} \notin I^-_u$. Therefore, as elements of 
$S_w^+$ and $S_u^-$,  
$$
c^+_{\xi_{s_i\varpi_i}} \De_{w \vpi_i, \vpi_i}^{-1} \neq 0 \quad \mbox{and} \quad 
\De^{-1}_{- u \vpi_i, - \vpi_i} c^-_{\xi_{-s_i\varpi_i}} \neq 0. 
$$
The definition of $p_i$ implies that
$[p_i]_{(w-s_i)\vpi_i, (s_i - u) \vpi_i} \neq 0$ which contradicts 
\eqref{p-dd}. This completes the proof.
\end{proof}

\sectionnew{Quantum cluster algebra structure on $S^+_w \bowtie S^-_u$ and reduced double Bruhat cells}
\label{qS+S-}
The section describes the interval prime elements of the CGL extensions $S^+_w \bowtie S^-_u$. 
We construct quantum cluster algebra structures on $S^+_w \bowtie S^-_u$ and on the reduced 
double Bruhat cells $R_q[G^{w,u}/H]$, exhibit a large number of explicit seeds for them, and prove that they coincide 
with the corresponding upper quantum cluster algebras. 

\subsection{Classification of the prime elements of the interval subalgebras of $S^+_w \bowtie S^-_u$}
\label{6.1} As in Section \ref{qdoubB} we fix  
reduced expressions \eqref{reduced} of the Weyl group elements $w$ and $u$.
Recall the notation $k \mt |k|$ from \eqref{abs}.
Define the function $\eta : [1,N+M] \to \Zset$ by
\begin{equation}
\label{etawu}
\eta(k) := 
\begin{cases}
i_{| k |}, & \mbox{if} \; \; k \in [1,N]
\\
i'_{| k |}, & \mbox{if} \; \; k \in [N+1, N+M].
\end{cases}
\end{equation}
Consider the following rescaled generators of the CGL extension presentation of 
$S^+_w \bowtie S^-_u$ from \thref{SwSu}:
\begin{equation}
\label{x-new}
x_k := 
\begin{cases}
q_{i_{ |k| }}^{1/2}(q_{i_{ |k| }}^{-1} - q_{i_{ |k| }}) C^+_{\be_{ |k| }}, & \mbox{if} \; \; k \in [1,N] \\
{q_{i'_{ |k| }}^{-1/2}(q_{i'_{ |k| }}^{-1} - q_{i'_{ |k| }}) q^{2 \lcor \be'_{ |k| }, \rho \rcor }}
C^-_{\be'_{ |k| }}, & \mbox{if} \; \; k \in [N+1,N+M],
\end{cases}
\end{equation}
recall \eqref{Celem}. {The intrinsic reason for the difference in the two normalizations is that 
the isomorphism $\om : \UU^-[u] \to \UU^+[u]$ transforms 
the generators of these CGL extensions as follows
\begin{equation}
\label{om-prop}
\om(X^-_{\be'_k}) = (-q_{i'_k}^{-1}) (-1)^{\lcor\be'_k, \rho\spcheck \rcor} 
q^{\lcor \be'_k, \rho \rcor} X^+_{\be'_k}, 
\end{equation}
where $\rho\spcheck$ is the sum of the fundamental coweights of $\g$,
see \cite[Eq. (2.10)]{GeY}. The first scalar in the formula is precisely 
the modification of the scalars in \eqref{x-new} from the first to the second cases.}

\bth{seq-primes} Let $\g$ be an arbitrary finite dimensional simple Lie algebra 
and $(w,u)$ a pair of Weyl 
group elements. Let $\KK$ be a base field of arbitrary characteristic 
and $q \in \KK^*$ a non-root of unity such that 
$\sqrt{q} \in \KK$. Consider the CGL extension presentation \eqref{SwSu1} 
of $S_w^+ \bowtie S_u^-$ with the rescaled generators \eqref{x-new}. Then the following hold:

{\rm{(a)}} The function $\eta$ from Theorem {\rm\ref{t1}} can be chosen as the function \eqref{etawu}.

{\rm{(b)}} The normalized interval prime elements \eqref{ybarismi} are given as follows. 

\hspace{0.4cm} {\rm{1)}} If $1 \leq j < k \leq N$ and $i_{|j|} = i_{|k|}(=: i)$, then 
\begin{equation}
\label{int1}
\ol{y}_{[j,k] } = \sqrt{q}^{ \, \| ( w_{\leq |j|} - w_{< |k|}) \vpi_i \|^2/2}
(\vp^+_w)^{-1} \left( ( \De_{w_{< |k|}\vpi_i, w_{\leq |j|} \vpi_i } \otimes \id) \, {}^{\tau \times 1}\RR^w \right).  
\end{equation}

\hspace{0.4cm} {\rm{2)}} If $N < j < k \leq N+M$ and $i'_{|j|} = i'_{|k|}(=:i)$, then 
\begin{multline}
\label{int2}
\ol{y}_{[j,k] } = \sqrt{q}^{ \, \| ( u_{\leq |k|} - u_{< |j|}) \vpi_i \|^2/2} \times
\\
(\vp^-_u)^{-1} \left( ( \De_{\vpi_i, \vpi_i} \otimes \id) \, \left( 
(T_{u_{\leq |k|}}^{-1} \otimes 1) \cdot  {}^{S \tau \times 1}\RR^u_{\opp} \, \cdot ( T_{u_{< |j|}} \otimes 1) \right) \right).  
\end{multline}

\hspace{0.4cm} {\rm{3)}} If $1 \leq j \leq N < k \leq N+M$ and $i_{|j|} = i'_{|k|}(=:i)$, then 
\begin{multline}
\label{int3}
\ol{y}_{[j,k] } = \sqrt{q}^{ \, \| ( u_{\leq |k|} - w_{\leq |j|}) \vpi_i \|^2/2} \times
\\
\vp_{w,u}^{-1} 
\left( \big( \De_{\vpi_i, w_{\leq |j|} \vpi_i } \otimes \id \big)  
\left( (T_{u_{\leq |k|}}^{-1} \otimes 1) \, \cdot {}^{S \tau \times 1} \RR^u_{\opp} \, {}^{\tau \times 1} \RR^w
\right)
\right)  
\end{multline}
in terms of the algebra isomorphism \eqref{vp-wu}.

{\rm{(c)}} The condition \eqref{cond} is satisfied.

Furthermore, all formulas remain valid if $\RR^w$ and $\RR^u_{\opp}$ are 
replaced by $\RR$ and $\RR_{\opp}$, respectively, cf. \eqref{Rw2}.
\eth
Note that the three braid group terms in \eqref{int2}-\eqref{int3} are the opposite to the ones defining the quantum minors \eqref{q-minor}.
They specialize to the Fomin--Zelevinsky elements $\ol{\ol{w}}$. 

\thref{seq-primes} is proved in \S \ref{6.3}.

\subsection{Quantum cluster algebra structures on $S_w^+ \bowtie S_u^-$ and $R_q[G^{w,u}/H]$}
\label{6.2} Recall the setting of \S \ref{thm8.2GYbig} and the definitions \eqref{tau}--\eqref{GaN} of the 
sets $\Xi_N$ and $\Ga_N$. We replace everywhere $N$ with $N+M$ to account for the GK dimension 
of $S_w^+ \bowtie S_u^-$.

Let $\eta$ be the function given by \eqref{etawu} and $\ex := \{ k \in [1,N+M] \mid s(k) \neq + \infty\}$, where $s$ is the successor function for $\eta$. 
Denote by $w^\ci_N$ the longest element of $S_N$, viewed as an element of $S_{N+M}$.
Let $p', s' : [1, N+M] \to [1,N+M] \sqcup \{ \pm \infty\}$ be the predecessor and successor 
functions for the function $\eta w^\ci_N \colon [1,N+M] \to \Zset$. Define 
$$
\ep : [1,N+M] \to \{ \pm 1 \}, \; \; \ep(k)= 1 \; \; \mbox{for} \; \; k>N, \; \; 
\ep(k)= -1 \; \; \mbox{for} \; \; k \leq N.
$$
We will use the notation $\ex_\sig$, $Z$ and $Z_\sig$ for $\sig \in \Xi_{N+M}$ 
from \S \ref{thm8.2GYbig}.
The Berenstein--Fomin--Zelevinsky exchange matrix $\ol{B}$ for the double word 
$-i_1, \ldots, -i_N, i'_1, \ldots, i'_M$ is given by \cite[Definition 2.3]{BFZ},
\begin{equation}
\label{BFZ-matr}
\ol{b}_{jk} :=
\begin{cases}
- \ep(k), & \mbox{if} \; \; j = p'(k)
\\
-\ep(k) c_{i_j, i_k}, & \mbox{if} \; \; j < k < s'(j) < s'(k), \; \ep(k) = \ep(s'(j))
\\
& \mbox{or} \; \; j < k < s'(k) < s'(j), \; \ep(k) = - \ep(s'(k))
\\
\ep(j) c_{i_j, i_k}, & \mbox{if} \; \; k < j < s'(k) < s'(j), \; \ep(j) = \ep(s'(k))
\\
& \mbox{or} \; \; k < j < s'(j) < s'(k), \; \ep(j) = - \ep(s(j))
\\
\ep(j), & \mbox{if} \; \; j = s'(k)
\\
0, & \mbox{otherwise}, 
\end{cases}
\end{equation}
where $j \in [1,N+M]$, $k \in \ex_{w^\ci_N}$. Denote its columns by $\ol{b}^k$. 
Define the exchange matrix $\wt{B}$ of size $(N+M) \times \ex$ with columns $b^l$
by 
$$
b^l = Z^{-1}  Z_{w^\ci_N} \cdot
\begin{cases} 
\ol{b}^l, & \mbox{if} \; \; 
l > N
\\
\ol{b}^{N+1-s(l)}, & \mbox{if} \; \; l\leq N \; \; \mbox{and} \; \; 
s(l) \leq N
\\
-\ol{b}^{N+1- p(l)} - \cdots - \ol{b}^{N+1 - p^m(l)}, &
 \mbox{if} \; \; l\leq N \; \; \mbox{and} \; \; 
s(l) > N, 
\end{cases}
$$
where $m:= \max \{j \mid p^j(l) \neq - \infty \}$ in the last case.
For $\sig \in \Xi_{N+M}$, define the exchange matrix $\wt{B}_\sig$ of size $(N+M) \times \ex_\sig$ with columns $b^l_\sig$ using \eqref{b-sig} in
\thref{Btau}. By way of definition, $\wt{B}_{w^\ci_N} = \ol{B}$.

There are toric frames $M_\sig$ of $\Fract(S_w^+ \bowtie S_u^-)$ for $\sig \in \Xi_{N+M}$ defined as in \S \ref{thm8.2GYbig}. In particular,
\begin{equation}
\label{Mfr}
M_\sig(e_k) =
\begin{cases} 
\ol{y}_{[p^m(\sig(k)), \sig(k)]},  & \mbox{if} \; \; \sig(1) \leq \sig(k), \; \mbox{where} \; m: = \max \{l \mid p^l(\sig(k)) \in \sig([1,k]) \}
\\ 
\ol{y}_{[\sig(k), s^m(\sig(k))]}, & \mbox{if} \; \; \sig(1) > \sig(k), \; \mbox{where} \; m: = \max \{l \mid s^l(\sig(k)) \in \sig([1,k]) \}
\end{cases}
\end{equation}
for $k \in [1,N+M]$.
The matrix $\rb_\sig$ of $M_\sig$ is given by
$$
(\rb_\sig)_{kj} = \Om_{\nub_\sig}( \ol{e}_k , \ol{e}_j ), \quad \forall \; k,j \in [1, N+M],
$$
where $\nub_\sig$ is the multiplicatively skew-symmetric matrix with entries
\begin{equation}
\label{nub}
(\nub_\sig)_{kj} = \begin{cases}
\sqrt{q}^{ \, - \langle \be_{|\sig(k)|}, \be_{|\sig(j)|} \rangle }, &\text{if} \; \; 1 \le \sig(j) < \sig(k) \le N  \\
\sqrt{q}^{ \, - \langle \be'_{|\sig(k)|}, \be'_{|\sig(j)|} \rangle }, &\text{if} \; \; N <\sig(j) < \sig(k) \le N+M  \\
\sqrt{q}^{ \, \langle \be'_{|\sig(k)|}, \be_{|\sig(j)|} \rangle }, & \text{if} \; \; 1 \le \sig(j) \le N < \sig(k) \le N+M
\end{cases}
\end{equation}
(compare \eqref{lakjbowtie}).

In the case of
the longest element of $S_N$, $\sig = w^\ci_N$, the pair $(M_{w^\ci_N}, \rb_{w^\ci_N})$ has a particularly simple form. 
We will denote it by $(\ol{M}, \ol{\rb})$; it will play a special role in the rest of the paper. For $k \leq N$ and $k>N$, $\ol{M}(e_k)$ 
is given by
\begin{equation}
\label{olM1}
\ol{M}(e_k) = \sqrt{q}^{ \, \| ( w_{\leq k} - 1 ) \vpi_{i_k} \|^2/2}
(\vp^+_w)^{-1} \left( ( \De_{\vpi_{i_k}, w_{\leq k} \vpi_{i_k} } \otimes \id) \, {}^{\tau \times 1}\RR^w \right) \; \; \mbox{and}  
\end{equation}
\begin{multline}
\label{olM2}
\ol{M}(e_k) = \sqrt{q}^{ \, \| ( u_{\leq k -N } - w) \vpi_{i'_{k-N}} \|^2/2} \times
\\
\vp_{w,u}^{-1} 
\left(  \big( \De_{\vpi_{i'_{k-N}}, w \vpi_{i'_{k-N}} } \otimes \id \big) 
\left( (T_{u_{\leq  k-N}}^{-1} \otimes 1) \cdot \, {}^{S \tau \times 1} \RR^u_{\opp} \, {}^{\tau \times 1} \RR^w
\right)
\right), 
\end{multline}
respectively. 
Set $\ol{\ex}:= \ex_{w_N^\ci}$.

\bth{SwSu-qcluster}  Let $\g$ be a complex simple Lie algebra, $w, u \in W$, $\KK$ an arbitrary 
base field and $q \in \KK$ a non-root of unity such that $\sqrt{q} \in \KK$.

{\rm(a)} For all $\sig \in \Xi_{N+M}$, the pair $(M_\sig, \wt{B}_\sig)$ is a 
quantum seed for $\Fract(S^+_w \bowtie S^-_u)$. The principal part of $\wt{B}_\sig$ 
is skew-symmetrizable via the integers $d^\sig_k$, $k \in \ex_\sig$, where 
$d^\sig_k = \| \al_{i_{| \sig(k) |}} \|^2/2$ if $\sig(k) \leq N$ and $d^\sig_k = \| \al_{i'_{| \sig(k) |}} \|^2/2$
if $\sig(k) >N$.

{\rm(b)} These quantum seeds are mutation-equivalent to each other up to the $S_{N+M}$-action and 
are linked by the following mutations. Let $\sig, \sig' \in \Xi_{N+M}$ be such that
$$
\sig' = ( \sig(k), \sig(k+1)) \sig = \sig (k, k+1)
$$
for some $k \in [1,N+M-1]$.
If $\eta(\sig(k)) \neq \eta (\sig(k+1))$, then $M_{\sig'} = M_\sig \cdot (k,k+1)$
in terms of the action \eqref{r-act}.
If $\eta(\sig(k)) = \eta (\sig(k+1))$, 
then $M_{\sig'} = \mu_k (M_\sig) \cdot (k,k+1)$.

{\rm(c)} We have the following equality between the algebra $S^+_w \bowtie S^-_u$
and the quantum cluster and upper cluster algebras associated to the seed $(\ol{M}, \ol{B})$
with no inverted indices {\rm(}i.e., $\inv:=\varnothing${\rm)}:
$$
S^+_w \bowtie S^-_u = \AA(\ol{M}, \ol{B}, \varnothing) = \UU(\ol{M}, \ol{B}, \varnothing).
$$
Each generator $x_k$ of $S^+_w \bowtie S^-_u$ is a cluster variable in one of the seeds 
corresponding to $\sig \in \Ga_{N+M}$. In particular, the cluster variables in those seeds 
generate $S^+_w \bowtie S^-_u$.

{\rm(d)} The quantized coordinate ring $R_q[G^{w,u}/H]$ of the reduced double Bruhat cell
$G^{w,u}/H$ equals the corresponding quantum and upper quantum cluster algebras with all 
frozen variables inverted:
$$
R_q[G^{w,u}/H] = \AA(\ol{M}, \ol{B}, [1,N+M] \backslash \ex) = \UU(\ol{M}, \ol{B}, [1,N+M] \backslash \ex).
$$
\eth
\bex{sl2} Consider the case $G=SL_2$, $w=u=s_1$. Using the isomorphism 
\[
\vp_{s_1}^+ \otimes \vp_{s_1}^- : S^+_{s_1} \bowtie S^-_{s_1} \stackrel{\cong}{\lra}  
\UU^-[s_1]_{\opp} \otimes \UU^+[s_1]
\]
and \reref{comm}, we see that the algebra $S^+_{s_1} \bowtie S^-_{s_1}$ is isomorphic to the $\KK$-algebra 
with generators $x^+_1$ and $x^-_1$ and the relation 
$$
x^+_1 x^-_1 = q^2 x^-_1 x^+_1 + (1 - q^2).
$$
Here 
$$
x^{\pm}_1 := q^{1/2}(q^{-1} - q) X_1^\pm.
$$
In this case \thref{SwSu-qcluster} states that $S^+_{s_1} \bowtie S^-_{s_1}$ is isomorphic 
to the quantum cluster algebra of type $A_1$ having two quantum seeds with cluster variables
\[
(x^-_1, q(x^-_1 x^+_1 -1)) \quad \mbox{and} \quad (x^+_1, q(x^-_1 x^+_1 -1)),
\] 
respectively. Their mutation quivers are  $1 \rightarrow 2$ and $1 \leftarrow 2$, respectively, with $2$ being a frozen 
variable. All normalizations above match those made in \eqref{olM1}--\eqref{olM2} and \thref{SwSu-qcluster} (a).
The seeds come from the two elements of $\Xi_2 =S_2$.

In general the quantum cluster algebra on $S^+_{w} \bowtie S^-_{u}$ is not of finite type.
\eex
\bex{sl3f}
{The algebra $R$ from \exref{sl3a} is isomorphic to the algebra $\UU^-_{\opp} \bowtie \UU^+$ for $G = SL_3$, where the generators $X_i^\pm$ of 
the first algebra are sent to the generators of the second denoted in the same way. The elements $X_1^+, X_2^+ \in R$ generate a copy of $\UU^+$ 
inside of $R$ and the element $X_{12}^+$ corresponds to the root vector $X^+_{\al_1 + \al_2} = T_1 X_2^+$ used in \thref{SwSu} for $w = s_1 s_2 s_1$. 
The elements $X_1^-, X_2^- \in R$ generate a copy of $\UU^-$ inside of $R$ which is isomorphic to $\UU^-_{\opp} \subset \UU^-_{\opp} \bowtie \UU^+$
via the antiisomorphism $\tau$. Under the map $\tau$, the element $X_{12}^- \in R$ corresponds to the root vector $X^-_{\al_1 + \al_2} = T_1 X_2^-$ used 
in \thref{SwSu} for $u = s_1 s_2 s_1$. In this way the symmetric CGL extension presentation of $R$ from \exref{sl3a} matches the one 
of  $\UU^-_{\opp} \bowtie \UU^+$ (for $G = SL_3$) from \thref{SwSu}. The quantum cluster algebra structure on 
 $\UU^-_{\opp} \bowtie \UU^+$ constructed in \thref{SwSu-qcluster} is precisely the one described in \exref{sl3e}. 
 }
\eex
\begin{proof}[Proof of \thref{SwSu-qcluster}] \thref{SwSu} implies that $S^+_w \bowtie S^-_u$ is a symmetric CGL extension 
with respect to the presentation \eqref{SwSu1}. By \coref{la-k} the conditions (A) and (B) in \S \ref{normalizations} are satisfied. 
It follows from \thref{seq-primes} that the condition \eqref{cond} holds for the sequence of rescaled generators $x_1, \ldots, x_{N+M}$ from \eqref{x-new}.
Thus, we can apply \thref{maint}.  

Consider the matrix $\ol{\Psi}:= (\ol{\psi}_{kj}) \in M_{M+N}(\Zset)$ such that $(\rb_{w^\circ_N})_{kj} = \sqrt{q}^{\, \ol{\psi}_{kj}}$.
Its entries are computed from the exponents in \eqref{nub}. 
Showing that $(\ol{M}, \ol{B})$ is a quantum seed amounts to showing that
$$
(\ol{\Psi} \,  \ol{B})_{kj} = - \de_{kj} d_k 
$$
where $d_k = \| \al_{i_k} \|^2/2$ if $k \leq N$ and $d_k = \| \al_{i'_{k-N}} \|^2/2$ if $k>N$.  A compatibility of this type was 
proved by Berenstein and Zelevinsky in \cite[Theorem 8.3]{BZ} where the matrix $\ol{B}$ was replaced with one that has $r$ more 
rows and the matrix in place of $\ol{\Psi}$ is slightly different. However, if one subtracts the entries in the above compatibility 
relation from the ones in \cite{BZ} then the result is easily verified to be 0 as a consequence of the fact that
the Berenstein--Fomin--Zelevinsky cluster variables
in \cite[Theorem 2.10 and Lemma 2.12]{BFZ} are homogeneous with respect to the $\PP \times \PP$ grading of $\Cset[G]$.

The graded nature of the one-step mutations in the proof of \cite[Theorem 2.10]{BFZ} also imply that 
$\chi_{\ol{M}(\ol{b}^l)} = 1$ for all $l$ (recall that $\ol{b}^l$ denotes the $l$-th column of $\ol{B}$).

The above compatibility and \thref{seq-primes} imply that the quantum seed of $S^+_w \bowtie S^-_u$ obtained by applying \thref{maint} (b)
to this algebra and $\sig = w^\ci_N$ is precisely $(\ol{M}, \ol{B})$. Now parts (a) and (b) of the theorem follow from \thref{Btau}, 
\thref{seq-primes}, and \thref{maint} (b)-(c). 

Part (c) of the theorem follows from \thref{maint} (d).
Part (d) follows from the isomorphism \eqref{redDB}, the fact that the prime elements $p_i$, $i \in \SS(w,u)$ are precisely 
the frozen variables of the toric frame $\ol{M}$ (up to non-zero scalars), and \thref{maint} (e).
\end{proof}

\subsection{Proof of \thref{seq-primes}}
\label{6.3}
The validity of \eqref{int1} and the fact that the restriction to $[1,N]$ 
of the function $\eta$ from \eqref{etawu} satisfies the conditions of \thref{1} 
was established in \cite[Theorem 10.1]{GYbig}. {Eq. \eqref{int2}
follows from \eqref{int1} by applying the property \eqref{om-prop} of the 
automorphism $\om$ of $\UU_q(\g)$ and \prref{omS}(i) below.} 

We proceed with \eqref{int3}. It follows from \eqref{pp1} that
\eqref{int3} is valid up to a scalar. Denote by $z$ the $\vp_{w,u}^{-1}$ term of \eqref{int3}. {By \cite[Eq (3.9)]{GeY}, 
$$
T_{(w_{\leq k-1})^{-1}}^{-1}v_{\vpi_i} = \tau \big( X_{\be_k}^+ \big)^{\de_{i i_k}} T_{(w_{\leq k})^{-1}}^{-1} v_{\vpi_i},
$$
and thus
$$
\De_{\vpi_i, \vpi_i} \big(  T_{(w_{\leq k-1})^{-1}}  \tau \big( X_{\be_k}^+ \big)^{\de_{i i_k}} T_{(w_{\leq k})^{-1}}^{-1} v_{\vpi_i} \big) =1.
$$
By applying \prref{omS}(i) below and interchanging the roles of $w$ and $u$, we obtain
$$
\De_{\vpi_i, \vpi_i} \big(T^{-1}_{u_{\leq k}} \big( S \om \tau \big( X_{\be'_k}^+ \big) \big)^{\de_{i i_k}} T_{u_{\leq {k-1}}} \big) =1.
$$
Since the extremal weight spaces of the module $V(\vpi_i)$ are one dimensional, 
$$
\big( S \om \tau \big( X_{\be'_k}^+ \big) \big)^{\de_{i i_k}} T_{u_{\leq {k-1}}} v_{\vpi_i} = T^{-1}_{u_{\leq k}} \vpi_i.
$$
These identities and \eqref{om-prop} imply that}
\begin{equation}
\label{lead-z}
\mbox{the leading term of $z$ is} \; \; 
\sqrt{q}^{\, (n+1)\| \al_i \|^2/2} x_k x_{p(k)} \ldots x_j
\end{equation}
where $n \in \Zset_{>0}$ is such that $s^n(j)=k$. The above is an ordered product of generators 
of $S^-_u$ and then generators of $S^+_w$. Lemma 10.3(a) of \cite{GYbig} and \prref{int} imply
$$
\ol{y}_{[j,k] } = \prod_{l \leq |j|, \, m \leq |k|, \, i_l= i'_m=i} \sqrt{q}^{\; \,   - \lcor \be_l, \be'_m \rcor}
\sqrt{q}^{\, (\| (1 - w_{\leq |j|}) \vpi_i \|^2 + \|( 1- u_{\leq |k|}) \vpi_i \|^2)/2} z.
$$
The first product equals $\sqrt{q}^{\; \, - \lcor (1 - w_{\leq |j|}) \vpi_i,  ( 1- u_{\leq |k|}) \vpi_i \rcor }$ 
since
$$
\sum_{l \leq |j|, \, i_l=i} \be_l = (1 - w_{\leq |j|}) \vpi_i \quad
\mbox{and} \quad
\sum_{m \leq |k|, \, i'_m=i} \be'_m = ( 1- u_{\leq |k|}) \vpi_i.
$$
This proves \eqref{int3}. 
The rest of part (a) of the theorem concerning the restriction 
of $\eta$ to $[N+1, N+M]$ follows from the combination of \eqref{lead-z} and \thref{1}.

Part (c) of \thref{seq-primes}:
It follows from \cite[Lemma 10.7]{GYbig} that the condition \eqref{cond}
is satisfied when $1 \leq j < s(j) \leq N$. The fact that the condition \eqref{cond}
is satisfied when $N < j < s(j) \leq N+M$ follows from this fact and \prref{omS}(i) below 
by applying once again the property \eqref{om-prop} of $\om$. Last, we consider the case 
$1 \leq j \leq N < s(j) \leq N+M$. The third statement in \thref{SwSu2} 
implies that it is sufficient to verify the condition \eqref{cond} when 
$j=1$, $k= N+M$ and $i_1=i'_M$, $i_l \neq i_1$ for $l \in [2,N]$, $i'_m \neq i_M$ 
for $m \in [1,M-1]$. This is done analogously to the proof of 
\cite[Lemma 10.7]{GYbig}. \qed

\sectionnew{Modified Berenstein--Zelevinsky quantum seeds}
\label{App}
In this section we recall the definition of the Berenstein--Zelevinsky seeds for quantum double Bruhat cells, study the automorphism
$S \om$ on $R_q[G]$, and describe a modification of the BZ seeds via $S \om$. Throughout the section we assume that $\sqrt{q} \in \KK$.
\subsection{The Berenstein--Zelevinsky seeds}
\label{6a.1}
Fix two Weyl group elements $w$ and $u$ and reduced expressions of them as in \eqref{reduced}. Consider the reversed reduced expressions
\begin{equation}
\label{red-expr}
w^{-1} = s_{i_N} \ldots s_{i_1} \quad \mbox{and} \quad
u^{-1} = s_{i'_M} \ldots s_{i'_1},
\end{equation} 
and set
\begin{equation}
\label{reversed-sub}
w_{< k}^{-1} := s_{i_N} \ldots s_{i_{k+1}}, \; \; \forall k \in [1,N].  
\end{equation}
Consider the reduced double word
\begin{equation}
\label{double-word}
1, \ldots, r, i_1, \ldots, i_N, - i'_1, \ldots, - i'_M 
\end{equation}
in the terminology of \cite{BZ}. The last two parts form a reduced word for the Weyl group element $(u,w) \in W \times W$ where 
the generators of the first copy of $W$ are assigned negative indices. 

{\em{The Berenstein--Zelevinsky toric frame}} \cite{BZ} of $\Fract(R_q[G^{u,w}])$
associated to \eqref{double-word} has {\em{cluster variables}}
\begin{align}
M_{BZ}(e_i) &:=  \De_{\vpi_i, w^{-1} \vpi_i}, &&i \in [1, r], 
\nn
\\
\label{BZ-clvar}
M_{BZ}(e_{r+k}) &:= \De_{\vpi_{i_k}, w^{-1}_{<k} \vpi_{i_k}}, &&k \in [1, N], 
\\
M_{BZ}(e_{r+N+k}) &:= \De_{u_{\leq k} \vpi_{i'_k}, \vpi_{i'_k}}, &&k \in [1, M]
\nn
\end{align}
in the notation \eqref{reversed-sub}. 
The {\em{matrix of the toric frame}} is $\rb_{BZ}:=(q^{\mu_{jk}/2})$ where
\begin{equation}
\label{BZ-matrix}
\mu_{jk} := \lcor \ga_j, \ga_k \rcor - \lcor \de_j, \de_k \rcor \; \; \mbox{for $j>k$}, \; \; 
\mu_{jk} :=  - \mu_{kj} \; \mbox{for $j<k$}, \; \; \mu_{jj} =0,
\end{equation}
and $\ga_k, \de_k$ are the first and second weights of the quantum minor defining $M_{BZ}(e_k)$, cf. \cite[Eqs. (8.5), (10.13)]{BZ},
We interchange the roles of $u$ and $w$ ($R_q[G^{u,w}]$ vs. $R_q[G^{w,u}]$)
to match the Berenstein--Zelevinsky picture to ours in Sect. \ref{qdoubB}. Consider the function $\eta \colon [1, r+N+M] \to \Zset$
defined by taking the absolute values of the numbers in the list \eqref{double-word}. Denote the corresponding predecessor/successor functions 
$p, s : [1,r+ N +M] \to [1,r+ N +M] \cup \{ \pm \infty \}$ by \eqref{pred.succ}. Define the set of {\em{exchangeable indices}} to
be 
\begin{equation}
\label{BZ-exch}
\ex:= \{ k \in [r+1, r+N+M] \mid s(k) \neq + \infty\}.
\end{equation}
The {\em{BZ quantum seed}} \cite[Eq. (8.7)]{BZ} is the pair $(M_{BZ}, \wt{B}_{BZ})$ for the {\em{exchange matrix}} $\wt{B}_{BZ}$ with entries $b_{jk}$ given by
\begin{equation}
\label{BZ-exmatr}
b_{jk}:= 
\begin{cases}
- \epsilon (k), & \mbox{if} \; \; j = p(k) \\
- \epsilon (k) c_{\eta(j) \eta(k)}, & \mbox{if} \; \; j < k < s(j) < s(k) \; \mbox{and} \; \ep(k) = \ep(s(j))  \\
                                                   & \mbox{or} \; \; j < k \leq r+N < s(k) < s(j)         \\
 \epsilon (j) c_{\eta(j) \eta(k)}, & \mbox{if} \; \; k < j < s(k) < s(j) \; \mbox{and} \; \ep(j) = \ep(s(k))  \\
                                                   & \mbox{or} \; \; k < j \leq r+N < s(j) < s(k)         \\
\epsilon (j), & \mbox{if} \; \; j = s(k) \\      
0, & \mbox{otherwise}                                                                                              
\end{cases}
\end{equation}
for $j \in [1, r+N+M]$, $k \in \ex$, where $\ep(k) := +1$ for $k \leq r+N$ and $\ep(k) := -1$ otherwise.

\subsection{The action of $S \om$ on quantum minors}
\label{6a.2}
Recall the definition of the automorphism $\om$  of $\UU_q(\g)$ from \S \ref{3.5}. It is a coalgebra antiautomorphism, 
so, $S^* \om^* = (\om S)^* \in \Aut R_q[G]$. For simplicity we will denote $\om^* : R_q[G] \to R_q[G]$ by $\om$, just like the antipode 
$S^*$ of $R_q[G]$ is denoted by $S$. The maps $S$ and $\om$ do not commute.

{\em{Throughout the rest of the paper we will use the antiautomorphism
$S \om$ of $\UU_q(\g)$ and the automorphism $\om S$ of $R_q[G]$.}} 

\bpr{omS} {\rm{(i)}} For all $u, w \in W$, $\mu \in \PP^+$, and $\iota_1, \iota_2 = \pm1$,
$$
\De_{\mu, \mu} (T_w^{\iota_1} \cdot  S \om(x) \cdot T_u^{\iota_2}) = \De_{\mu, \mu} (T_{u^{-1}}^{\iota_2} \cdot x \cdot T_{w^{-1}}^{\iota_1}), 
\quad \forall\,  x \in \UU_q(\g). 
$$
{\rm{(ii)}} For $u, w \in W$ and $\mu \in \PP^+$, 
$$
\om S(\De_{u \mu, w \mu}) = (-1)^{\lcor \rho\spcheck, (w-u)\mu \rcor}
q^{\lcor \rho, (w -u) \mu \rcor} \De_{w \mu, u \mu}
$$
where $\rho$, $\rho\spcheck$ are the sums of the fundamental weights {\rm(}resp. coweights{\rm)} of $\g$.
\epr

\begin{proof} Part (i). Using $S \om(K_i) = K_i$, $S \om(X^+_i) = - X^-_i K_i$, $S \om(X^-_i) = - K_i^{-1} X^+_i $,
we obtain that for any $x = x_- x_0 x_+ \in \UU_q(\g)$ where $x_\pm \in \UU^\pm$, $x_0 \in \UU^0$,
$$
\De_{\mu, \mu} (x) = \varepsilon(x_-) \varepsilon(x_+) \De_{\mu, \mu} (x_0) 
= \De_{\mu, \mu} ( S \om(x)),
$$
because $S \om (x_0) = x_0$ ($\varepsilon$ denotes the counit of $\UU_q(\g)$). This proves the 
special case $w=u=1$.
For the general case, define the operators $S \om (T_w) \colon V(\mu) \to V(\mu)$ by 
$$
S \om (T_w) := S \om(T_{i_N}) \ldots S \om(T_{i_1})
$$
and by setting 
$$
(S \om (T_i)) v := \sum_{j,k,l \geq 0, j - k + l =\lcor \la, \al_i\spcheck\rcor} \frac{(-1)^k q_i^{k -jl}}{[j]_{q_i}! [k]_{q_i} ! [l]_{q_i} !} 
(- X^-_i K_i)^j  (- K_i^{-1} X^+_i )^k (- X^-_i K_i)^l v
$$
on $v \in V(\mu)_\la$. The latter formula is obtained by the direct application of $S \om$ to the element of $\UU_q(\g)$ defining $T_w$, cf. \cite[Eq. 8.6 (2)]{Ja}.
By commuting the $K_i$ elements 
to the right, one obtains that $(S \om(T_i))(v) = (- q_i)^{\lcor \la, \al_i\spcheck \rcor} \, {}^\om T_i (v)$ in terms of the notation 
\cite[Eq. 8.6 (2)]{Ja}. It follows from the first identity in \cite[Eq. 8.6 (7)]{Ja} that $S \om(T_i) = T_i$, and thus, that
$S \om (T_w) = T_{w^{-1}}$. Now we can apply the special case,
\begin{align*}
\De_{\mu, \mu} ( T_u  \cdot S \om (x) \cdot T_w) &= \De_{\mu, \mu} ( (S \om) T_{u^{-1}} (S \om) (x) \om S(T_{w^{-1}}))  \\
& =  \De_{\mu, \mu} (S \om ( T_{w^{-1}} \cdot x \cdot T_{u^{-1}} ) ) = \De_{\mu, \mu} (T_{w^{-1}} \cdot x \cdot T_{u^{-1}})), \quad \forall\,  x \in \UU_q(\g),
\end{align*}
the point being that in the last expression $T_{w^{-1}} \cdot x \cdot T_{u^{-1}}$ acts on $v_\mu$ by an element of $\UU_q(\g)$ and the previous 
term is obtained by applying $S \om$ to it. The other cases for $\iota_1, \iota_2 = \pm 1$ are analogous.

Part (ii). One checks that
$$
T_w v_\mu = (-1)^{\lcor \rho\spcheck - w^{-1}\rho\spcheck, \mu \rcor} 
q^{\lcor \rho - w^{-1}\rho, \mu \rcor} T_{w^{-1}}^{-1} v_\mu
$$
by iterating the special case of this formula for $\UU_q(\sl_2)$. This, \eqref{q-minor}, and part (i) give
\begin{align*}
(\om S(\De_{u \mu, w \mu}))(x) &= \De_{\mu, \mu} (T_{u^{-1}} \cdot S \om(x) \cdot T_{w^{-1}}^{-1}) = 
\\
&= \De_{\mu, \mu} (T_{w}^{-1} \cdot x \cdot T_u) = a 
\De_{\mu, \mu} (T_{w^{-1}} \cdot S \om(x) \cdot T_{u^{-1}}^{-1}) = a \De_{w \mu, u \mu} (x)
\end{align*}
where $a= (-1)^{\lcor (w^{-1} - u^{-1}) \rho\spcheck,\mu \rcor} q^{\lcor (w^{-1} - u^{-1}) \rho, \mu \rcor} 
= (-1)^{\lcor \rho\spcheck, (w-u)\mu \rcor} q^{\lcor \rho, (w -u) \mu \rcor} $.
\end{proof}

\subsection{The action of $S \om$ on quantum double Bruhat cells}
\label{6a.3}
\bpr{SRq} Let $u, w \in W$.

{\rm{(i)}} $S^{\pm1}(I_{u,w}) = I_{u^{-1},w^{-1}}$
and $S^{\pm 1}$ induce antiisomorphisms between $R_q[G^{u,w}]$ and 
$R_q[G^{u^{-1},w^{-1}}]$.

{\rm{(ii)}} $\om^{\pm 1} S^{\pm 1} (I_{u,w}) = I_{w^{-1},u^{-1}}$
and $\om^{\pm 1} S^{\pm 1}$ induce isomorphisms between
$R_q[G^{u,w}]$ and $R_q[G^{w^{-1},u^{-1}}]$.
\epr

\begin{proof} Part (i).
For $i \in [1,r]$ denote by $\BB^\pm_i$ the subalgebra of $\UU_q(\g)$ 
generated by $X_i^\pm$ and $K_1^{\pm 1}, \ldots, K_r^{\pm 1}$.
It follows from \cite[Corollary 4.4.6]{J} and \cite[Proposition 9.2]{BZ} that
the ideals $I_{u,w}$ are given by 
$$
I_{u,w} = (\BB_{i'_1}^- \ldots \BB_{i'_M}^- \BB_{i_1}^+ \ldots \BB_{i_N}^+)^\perp
= (\BB_{i_1}^+ \ldots \BB_{i_N}^+ \BB_{i'_1}^- \ldots \BB_{i'_M}^-)^\perp
$$
in terms of the reduced expressions \eqref{reduced},
where the orthogonal complements are computed with respect to the natural 
pairing between $R_q[G]$ and $\UU_q(\g)$. The first part of the proposition 
follows from this fact by using that $S(\BB_i^\pm) = \BB_i^\pm$ for all $i \in [1,r]$. 
Therefore $S^{\pm 1}$ induce antiisomorphisms
between 
$$
R_q[G]/I_{u,w} \quad \mbox{and} \quad
R_q[G]/I_{u^{-1},w^{-1}}.
$$ 
Denote by $M_{w,u}$ the multiplicative subset 
of $R_q[G^{w,u}]$ generated by $\kx$ and the elements $\De_{w \vpi_i, \vpi_i}^{\pm}$,
$\De_{-u\vpi_i, - \vpi_i}^{\pm}$ for $i \in [1,r]$. The elements $S^{\pm 1}(\De_{u\vpi_i, \vpi_i})$ 
and $S^{\pm 1}(\De_{- w\vpi_i, - \vpi_i})$ are normal in $R_q[G^{u^{-1},w^{-1}}]$ and
homogeneous with respect to the $\QQ \times \QQ$ grading. By \coref{2}
all of them belong to $M_{w^{-1},u^{-1}}$, so $S^{\pm 1} (M_{u,w}) \subseteq M_{u^{-1},w^{-1}} $.
Using that $S^{-1}$ is the inverse of $S$ gives the inverse inclusion, and thus
$$
S^{\pm 1}(M_{u,w}) = M_{u^{-1},w^{-1}} \quad \mbox{and} \quad 
R_q[G^{u,w}] = (R_q[G]/I_{u,w})[(S^{\pm 1}(E_{u^{-1},w^{-1}}))^{-1}]
$$ 
which implies the second part of (i). Part (ii) is proved analogously.
\end{proof}

\subsection{A modification of the Berenstein--Zelevinsky seeds via $\om S$}
\label{6a.4}
The isomorphism $\om S : R_q[G^{u,w}] \to R_q[G^{w^{-1}, u^{-1}}]$ from \prref{SRq} 
and the relations in \prref{omS} (ii) give a quantum seed of 
$\Fract(R_q[G^{w^{-1}, u^{-1}}])$ which we will call the {\em{modified Berenstein--Zelevinsky quantum seed}}. The 
{\em{cluster variables}} of the {\em{toric frame}} are given by 
\begin{align}
M_{mBZ}(e_i) &:=  \De_{w^{-1} \vpi_i, \vpi_i}, && i \in [1, r], 
\nn
\\
M_{mBZ}(e_{r+k}) &:= \De_{w^{-1}_{<k} \vpi_{i_k}, \vpi_{i_k}}, && k \in [1, N], 
\label{mod-BZ-clvar}
\\
M_{mBZ}(e_{r+N+k}) &:= \De_{\vpi_{i'_k}, u_{\leq k} \vpi_{i'_k}}, && k \in [1, M].
\nn
\end{align}
The {\em{matrix of the toric frame}} is given by \eqref{BZ-matrix} and {\em{exchangeable indices}} by \eqref{BZ-exch}. Define 
the modified BZ quantum seed to be the pair $(M_{mBZ}, \wt{B}_{BZ})$. 

Define an action of $\HH = (\kx)^r$ on $R_q[G]$ by 
$$
h \cdot c = h^{\nu' + \nu} c \quad \mbox{for} \quad c \in R_q[G]_{\nu, \nu'}, h \in \HH,
$$
recall the identification $\xh \cong \QQ$ in \eqref{Hchar}. There exists a unique $h_* \in \HH$ such that 
\begin{equation}
\label{hst}
h_*^\ga = (-1)^{ \lcor \rho\spcheck, \ga \rcor} q^{\lcor \rho, \ga \rcor}, \quad \forall\,  \ga \in \QQ,
\end{equation}
recall that $\PP\spcheck \subseteq \PP$. Taking into account that $\De_{w \mu, u \mu} \in R_q[G]_{- w\mu, u\mu}$ and 
\prref{omS} (ii) we obtain,
\bpr{BZtomBZ} For all $w, u \in W$, the $\om S$-image of the Berenstein--Zelevinsky quantum seed and the modified BZ quantum seed 
are related by an automorphism of $R_q[G^{w^{-1}, u^{-1}}]$, 
$$
M_{mBZ}(f) = h_* \cdot (\om S) M_{BZ}(f), \quad \forall\,  f \in \Zset_{r+N+M}.
$$ 
\epr

\sectionnew{A quantum twist map for double Bruhat cells}
\label{twist}
In this section we develop a quantum analog of the Fomin--Zelevinsky twist map 
\cite{FZ} for all double Bruhat cells. We use it to
relate the cluster variables in the quantum seed in \thref{SwSu-qcluster} for $\Fract(S^+_w \bowtie S^-_u)$
corresponding to the element $\sig = w^\ci_N \in \Xi_{N+M}$ to the cluster variables in the modified Berenstein--Zelevinsky quantum seed \cite{BZ} 
for $\Fract(R_q[G^{w^{-1}, u^{-1}}])$ described in the previous section. Throughout the section, $\g$ is an arbitrary simple Lie algebra,
$w, u$ are Weyl group elements with fixed reduced expressions as in \eqref{reduced}, and 
reversed expressions \eqref{red-expr} are used for the elements $w^{-1}$, $u^{-1}$. 

\subsection{First half of the construction: a twist map for quantum Schubert cells}
\label{7.1}
For $w \in W$ define the algebra antiisomorphism
$$
\zeta_{w^{-1}} := \tau S^{-1} \tau \om T_w^{-1} \colon \UU_q(\g) \to \UU_q(\g).
$$
The double appearance of $\tau$ is needed because of the nature of the isomorphisms in \thref{isom} and 
the fact that $S$ and $\tau$ do not commute. 

\bpr{zeta-1} For all simple Lie algebras $\g$ and $w \in W$, 
the map $\zeta_{w^{-1}} : \UU_q(\g) \to \UU_q(\g)$ restricts to antiisomorphisms 
$$
\zeta_{w^{-1}} \colon \UU^{\pm}[w] \to \UU^{\pm}[w^{-1}]
$$
and satisfies
$$
\zeta_{w^{-1}}( T_{i_1} \ldots T_{i_{k-1}}(X^{\pm}_{i_k})) = 
\tt^{\pm}_{w, k} T_{i_N} \ldots T_{i_{k+1}}(X^{\pm}_{i_k}), \quad 
\forall\,  k \in [1,N]
$$
for some $\tt^\pm_{w, k} \in \kx$, in terms of \eqref{reduced}, \eqref{red-expr}. Moreover, $\tt_{w,k}^+ \tt_{w,k}^- = 1$
for $k \in [1,N]$, so,
$$
(\zeta_{w^{-1}} \otimes \zeta_{w^{-1}})  \RR^w = \RR^{w^{-1}}.
$$
\epr

The map $\zeta_{w^{-1}}$ is inverse to the twist antiisomorphism 
$$
\Theta_w := T_w \om \tau S \tau : \UU^{\pm}[w^{-1}] \to \UU^{\pm}[w]
$$
studied in \cite{LeY}. The first two statements of the proposition follow from \cite[Proposition 6.1]{LeY}. The property 
that $\tt_{w,k}^+ \tt_{w,k}^- = 1$ follows from the fact that 
$$
T_w \om (x) = \tt_\ga \om T_w (x), \quad \forall\,  x \in \UU_q(\g)_\ga, \; \ga \in \QQ
$$
for certain $\tt_\ga \in \kx$ satisfying $\tt_\ga \tt_{-\ga} =1$ (see \cite[Eq. 8.18 (5)]{Ja}) and
from the fact that $\tau$ and $\om$ commute.

We will denote by the same notation the antiisimorphisms
\begin{equation}
\label{zetaw-2}
\zeta_{w^{-1}} \colon S^w_{\pm} \to S^{w^{-1}}_{\pm}
\end{equation}
obtained from $\zeta_w$ by composing it with the (anti-)isomorphisms $\vp^\pm_w : S^w_\pm \to \UU^\mp[w]$ from \thref{isom}
and their inverses. 

\ble{two-zetaw} For all $w \in W$, the twist antiismorphism $\zeta_{w^{-1}} \colon S^w_+ \to S^{w^{-1}}_+$ is given by 
$$
\zeta_{w^{-1}} \big( c_\xi^+ \De_{w \mu, \mu}^{-1} \big) = (\vp_{w^{-1}}^+)^{-1} \Big(
(c_\xi^+ \otimes \id)
\big( (T_{w^{-1}}^{-1} \otimes 1) \cdot {}^{\om S \tau \times 1} \RR^{w^{-1}} \big) \Big)
$$
for $\mu \in \PP^+$, $\xi \in V(\mu)^*$ {\rm{(}}recall \eqref{dp}{\rm{)}}, and the equality holds with  
$\RR^{w^{-1}}$ replaced by $\RR$.
\ele

\begin{proof} Using Propositions \ref{pomS} (i) and \ref{pzeta-1}, and \eqref{iota-ident}, we obtain
\begin{align*}
\vp_{w^{-1}}^+ \big( \zeta_{w^{-1}} ( c_\xi^+ \De_{w \mu, \mu}^{-1} ) \big) 
&= ( c_{\xi, v_{w \mu}} \otimes \zeta_{w^{-1}}) \; {}^{\tau \times 1} \RR^w  \\
&= ( c_{\xi, v_{w \mu}} \otimes \id) \; {}^{\tau \zeta_{w^{-1}}^{-1} \times 1} \RR^{w^{-1}} 
= ( c_{\xi, v_{w \mu}} \otimes \id)
\; {}^{T_{w^{-1}}^{-1} \om S \tau \times 1} \RR^{w^{-1}}  \\
&= (c_\xi^+ \otimes \id) \big( (T_{w^{-1}}^{-1} \otimes 1) \cdot \; {}^{\om S \tau \times 1} \RR^{w^{-1}} \big).
\end{align*}  
At the end one can replace $\RR^{w^{-1}}$ by $\RR$ because $c_\xi^+ = 0$ in $R_q[G^{w,u}]$ for $\xi \in (V(\mu)_{w \mu}^+)^\perp$.
\end{proof}

\subsection{Quantum Bruhat cells}
\label{7.2}
In the setting of \S \ref{3.3}, $I^-_{w_0} = 0$ and $I_{w,w_0} = I_w^+ R^-$. The quantized coordinate ring of the Bruhat cell $B^+ w B^+$ is defined to be
$$
R_q[B^+ w B^+] := (R_q[G]/I_w^+R^-) [(E_w^+)^{-1}].
$$
Since $I_w^+ R^- \subseteq I_{w,u}$ and $E_w^+ \subseteq E_{w,u}$ for $u \in W$, we have the canonical projection
\begin{equation}
\label{proj-wu}
\kappa : R_q[B^+ w B^+] \to R_q[G^{w,u}]
\end{equation}
Keeping in mind \eqref{prod}, define
\begin{equation}
\label{loc-minors-}
\De_{w \mu, \mu} \in R_q[B^+ w B^+] \quad \mbox{for} \; \;  \mu \in \PP
\end{equation}
by $\De_{w \mu, \mu} := \De_{w \mu_1, \mu_1} \De_{w \mu_2, \mu_2}^{-1}$
for any $\mu_i \in \PP^+$ such that $\mu = \mu_1 - \mu_2$. 

\bpr{olS}
The set
$$
\ol{S}^{\, +}_w := \{ c_{\xi_{w \mu}, v} \De^{-1}_{w \nu, \nu} \mid v \in V(\mu)_\nu, \; \mu \in \PP^+, \; \nu \in \PP \}
$$
is a subalgebra of $R_q[B^+ w B^+]$ and the map $\psi_w^+ : \ol{S}^{\, +}_w \to \UU^-$ given by 
$$
\psi_w^+ \big( c_{\xi_{w \mu}, v} \De^{-1}_{w \nu, \nu} \big) := (c_{\xi_\mu, v} \otimes \id) ( {}^{\tau \times 1} \RR )
$$
is a well defined algebra isomorphism.
\epr

\begin{proof} The fact that the set $\ol{S}_w^{\, +}$ is a subalgebra of $R_q[B^+ w B^+]$ follows from the identity 
$(T^{-1}_{w^{-1}} v_\la) \otimes (T^{-1}_{w^{-1}} v_\mu) = T^{-1}_{w^{-1}} (v_\la \otimes v_\mu)$
for $\la, \mu \in \PP^+$ (see e.g. \cite[(2.19)-(2.20)]{Y-sqg}). The second isomorphism is completely analogous to \thref{isom}
proved in \cite[Theorem 2.6]{Y-sqg}.
\end{proof}

\subsection{Definition of the quantum twist map for double Bruhat cells}
\label{7.3} Similarly to the previous subsection, define
$$
\De_{-u \mu, - \mu} \in R_q[G^{w, u}], \; \;  \De_{w^{-1} \mu, \mu} \in R_q[G^{w^{-1},u^{-1}}] \quad \mbox{for}
\; \;  \mu \in \PP
$$
by 
\begin{equation}
\label{loc-minors}
\De_{- u \mu, - \mu} = \De_{- u \mu_1, - \mu_1} \De_{- u \mu_2, - \mu_2}^{-1}, \quad
\De_{w^{-1} \mu, \mu} = \De_{w^{-1} \mu_1, \mu_1} \De_{w^{-1} \mu_2, \mu_2}^{-1}
\end{equation}
for any $\mu_i \in \PP^+$ such that $\mu = \mu_1 - \mu_2$.
Recall the algebra isomorphisms \eqref{vp-wu}. We define a linear map 
$$
\zeta_{w^{-1}, u^{-1}} : (S^+_w \bowtie S^-_u) \# \LL_u^- \to R_q[G^{w^{-1}, u^{-1}}] \cong 
\big( (S^+_{w^{-1}} \bowtie S^-_{u^{-1}})[p_i^{-1}, \; i \in \SS(w,u)] \big) \# \LL_{w^{-1}}^+
$$
by the following formulas on the three components of the left hand side identified (as a vector space) with the tensor product
$S^+_w \otimes S^-_u \otimes \LL_u^-$. 
Firstly, 
\begin{equation}
\label{zeta-wu-1}
\zeta_{w^{-1},u^{-1}}|_{S^+_w} := \zeta_{w^{-1}} \,.
\end{equation}
Secondly, we define $\zeta_{w^{-1},u^{-1}}|_{S^-_u}$ as the composition
\begin{equation}
\label{zeta_wu-2}
S^-_u \,  \stackrel{\vp^-_u}{\longrightarrow}  \, \UU^+[u] \,   \stackrel{\om}{\longrightarrow}   \, \UU^-[u] \hra \UU^- \,
\xrightarrow{(\psi^+_{w^{-1}})^{-1}} \,  \ol{S}^{\, +}_{w^{-1}} \,   \stackrel{\kappa}{\longrightarrow}  \, R_q[G^{w^{-1},u^{-1}}].
\end{equation}
The map $\zeta_{w^{-1},u^{-1}}|_{S^-_u}$ is an antihomomorphism because $\psi^+_{w^{-1}}$ is an antiisomorphism, 
$\vp^-_u$ and $\om$ are isomorphisms and $\kappa$ is a homomorphism.

Finally, we define
$\zeta_{w^{-1},u^{-1}}|_{\LL_u^-}$ by
\begin{equation}
\label{zeta-De1}
\zeta_{w^{-1},u^{-1}} ( \De_{- u \mu, - \mu}) := \De_{w^{-1}( u \mu), u \mu}, \quad \forall\,  \mu \in \PP  
\end{equation}
where the right hand side uses the notation \eqref{loc-minors} for the weight $u \mu \in \PP$. By \eqref{prod}, $\zeta_{w^{-1},u^{-1}}|_{\LL_u^-}$ is an (anti)homomorphism of commutative algebras. It is straightforward to check that
\begin{equation}
\label{zeta-weight}
\zeta_{w^{-1},u^{-1}} ((S^+_w \bowtie S^-_u) \# \LL_u^-)_{\mu, \nu}) 
\subseteq R_q[G^{w^{-1}, u^{-1}}]_{- w^{-1} \mu, - u \nu}, \quad
\mu, \nu \in \PP.
\end{equation} 

\bth{quantum-twist} For all simple Lie algebras $\g$ and $w, u \in W$, the following hold:

{\rm{(i)}} The map $\zeta_{w^{-1},u^{-1}}$ extends to an antiisomorphism $R_q[G^{w, u}] \to R_q[G^{w^{-1},u^{-1}}]$
whose inverse is $\zeta_{w,u}$. 

{\rm{(ii)}} The antiisomorphism $\zeta_{w^{-1},u^{-1}}$ links the modified BZ cluster variables to the cluster variables in the 
toric frame $\ol{M}$ in Theorem {\rm\ref{tSwSu-qcluster}} as follows:
\begin{align}
\label{zeta-De2}
\zeta_{w^{-1},u^{-1}}\big( \ol{M}(e_k) \big) &=
\sqrt{q}^{ \; - \| ( w_{\leq k} - 1) \vpi_{i_k} \|^2/2} 
\De_{w^{-1}_{<k} \vpi_{i_k}, \vpi_{i_k}}
\De_{ w^{-1} \vpi_{i_k}, \vpi_{i_k}}^{-1},
\\
\label{zeta-De3}
\zeta_{w^{-1},u^{-1}} \big( \ol{M}(e_{N+j}) \big) &= 
\sqrt{q}^{ \; - \| ( u_{\leq j } - w) \vpi_{i'_j} \|^2/2}
\De_{\vpi_{i'_j}, u_{\leq j} \vpi_{i'_j}} 
\De_{ w^{-1}( u_{\leq  j} \vpi_{i'_j}), u_{ \leq j} \vpi_{i'_j}}^{-1},
\end{align}
for all $k \in [1,N]$, $j \in [1,M]$ in the notation \eqref{reversed-sub} and \eqref{loc-minors}.
\eth

We will call $\zeta_{w^{-1}, u^{-1}} : R_q[G^{w, u}] \to R_q[G^{w^{-1},u^{-1}}]$ the {\em{quantum twist map}}.
In \S \ref{7.4} we describe its relation to the classical twist map of Fomin--Zelevinsky \cite{FZ}. In \S \ref{7.5} 
we prove that the map $\zeta_{w^{-1},u^{-1}}$ is an antihomomorphism. In \S \ref{7.6} and \ref{7.7} we prove parts (ii) and (i) of 
\thref{quantum-twist}, respectively.

\subsection{Motivation for the quantum twist map}
\label{7.4}
Recall that $G$ denotes a connected, simply connected complex simple {algebraic} group with Lie algebra $\g$,
$B^\pm$ a pair of opposite Borel subgroups with unipotent radicals $U^\pm$, and $H= B^+ \cap B^-$ the 
corresponding maximal torus of $G$. The set $G_0 = U_- H U_+$ is a Zariski open subset of $G$ and we have the canonical projections 
$$
G_0 \to U^-, \; H, \; U^+, \quad x \mt [x]_-, \; [x]_0, \;  [x]_+
$$ 
where $x = [x]_- [x]_0 [x]_+$ and $x \in G_0$, $[x]_- \in U^-$, $[x]_0  \in H$, $ [x]_+ \in U^+$. The elements $w \in W$ have two representatives 
$\ol{w}$, $\ol{\ol{w}}$ in the normalizer of $H$ in $G$, defined in \cite[\S 1.4]{FZ}. It is easy to show that their action on integrable $\g$-modules are the 
specializations of $T_{w^{-1}}^{-1}$ and $T_w$, respectively. The involution $\om$ of $U_q(\g)$ specializes to the involution 
of $\g$ given by $\om(h_i) = - h_i$, $\om(e_i) = f_i$, $\om(f_i) = e_i$ in terms of the Chevalley generators of $\g$. The latter involution integrates to an involution 
of $G$ to be denoted by the same letter.  The Fomin--Zelevinsky twist map \cite{FZ} is the algebraic isomorphism 
$$
\zeta^{w,u} : G^{w,u} \to G^{w^{-1}, u^{-1}}, \quad 
\zeta^{w,u}(x) := \om \Big( \left[ {\ol{w}}^{\, -1} x \right]_-^{-1} {\ol{w}}^{\, -1} x \, {\ol{\ol{u}}}^{\, -1} \left[ x {\ol{\ol{u}}}^{\, -1} \right]_+^{-1}  \Big).
$$ 
Denote the canonical projections $\upsilon_\pm : G \to G/B^\pm$ and their restrictions $\upsilon_+ : G^{w,u} \to (B^+ w B^+)/B^+$ 
and $\upsilon_- : G^{w,u} \to (B^- u B^-)/B^-$. Following the specialization arguments in \cite[Sect. 4]{Y-cm}, one can show that the pullback algebras 
$$
A^+_w := \upsilon_+^*\big( \Cset[(B^+ w B^+)/B^+] \big) \subseteq \Cset[G^{w,u}], \quad 
A^-_u := \upsilon_-^*\big( \Cset[(B^- w B^-)/B^-] \big) \subseteq \Cset[G^{w,u}]
$$
are specializations of $S^+_w, S^-_u \subseteq R_q[G^{w,u}]$ with respect to the nonrestricted integral form of $R_q[G]$. Each element 
$x \in G^{w,u}$ has a unique decomposition of the form
$$
x = n_+ \ol{w} m_+ h = m_- h' \, \ol{\ol{u}} \, n_-
$$
where $n_+ \in U^+ \cap w U^- w^{-1}$, $n_- \in U^- \cap u^{-1} U^+ u$, $m_\pm \in U^\pm$ and $h, h' \in H$. 

\ble{zeta-class} In the above setting, the following hold for $x \in G^{w,u}$:
\begin{align}
(\zeta^{w,u})^*(f)(x) = f(\om(\ol{w}^{\, -1} n_+^{-1})), \quad \forall f \in A^+_{w^{-1}}, 
\label{geom-twist1}
\\
(\zeta^{w,u})^*(g)(x) = g(\om(m_+ \ol{u}^{\, -1})), \quad \forall g \in A^-_{u^{-1}}.
\label{geom-twist2}
\end{align}
\ele

\begin{proof} Using that $\left[ {\ol{w}}^{\, -1} x \right]_- = \ol{w}^{\, -1} n_+ {\ol{w}}$ and $\left[ x {\ol{\ol{u}}}^{\, -1} \right]_+ = \ol{\ol{u}} n_- \ol{\ol{u}}^{\, -1}$, 
and the invariance properties of the functions in $A^+_w$ and $A^-_u$, we obtain
$$
(\zeta^{w,u})^*(f)(x) = f(\om(\ol{w}^{\, -1} n_+^{-1} m_- h'))= f(\om(\ol{w}^{\, -1} n_+^{-1}))
$$ 
for $f \in A^+_{w^{-1}}$ and
$$
(\zeta^{w,u})^*(g)(x) = g(\om(m_+ h n_-^{-1} \ol{\ol{u}}^{\, -1})) 
= g(\om(m_+ \ol{u}^{\, -1}))
$$
for $g \in A^-_{u^{-1}}$.
\end{proof}

The defining formulas \eqref{zeta-wu-1}--\eqref{zeta-De1} for the quantum twist map were set up in such a way that 
they would quantize the geometric maps behind \eqref{geom-twist1}--\eqref{geom-twist2}. In this way $\zeta_{w^{-1},u^{-1}} : R_q[G^{w,u}] \to R_q[G^{w^{-1},u^{-1}}]$ 
is a quantization of the pullback of the Fomin-Zelevinsky map $\big(\zeta^{w^{-1},u^{-1}}\big)^* : \Cset[G^{w,u}] \to \Cset[G^{w^{-1},u^{-1}}]$. 

\subsection{$\zeta_{w^{-1},u^{-1}}$ is an antihomomorphism}
\label{7.5}
It was shown in \S\S \ref{7.1}, \ref{7.3} that the restrictions of $\zeta_{w^{-1},u^{-1}}$ to the three algebras $S^+_w$, $S^-_u$ and $\LL_u^-$
are antihomomorphisms. The algebra $\LL_u^-$ is spanned by the localized minors $\De_{- u \mu, - \mu}$, $\mu \in \PP$ which are 
normal elements of $R_q[G^{w^{-1},u^{-1}}]$. The corresponding commutation relations with the elements of $R_q[G^{w^{-1},u^{-1}}]$ 
are given by \eqref{n2}. It follows from \eqref{zeta-weight} and \eqref{n1} that $\zeta_{w^{-1},u^{-1}}$ anti-preserves those relations. It remains to prove that 
$\zeta_{w^{-1},u^{-1}}$ anti-preserves the commutation relation \eqref{commRR} between $S^+_w$ and $S^-_u$.

We will prove a stronger statement. It follows from the definitions of $S_{w^{-1}}^+$ and $R_q[B^+w^{-1}B^+]$ that the first algebra is canonically 
a subalgebra of the second. By \prref{olS},
 $\ol{S}^{\, +}_{w^{-1}}$ is also a subalgebra of $R_q[B^+ w^{-1} B^+]$.
 
\ble{zeta-le} The map $S^+_w \bowtie S^-_u \to R_q[B^+ w^{-1} B^+]$ given by $\zeta_{w^{-1}}$ on the first term and by the 
composition 
\begin{equation}
\label{comp}
S^-_u \, \stackrel{\vp^-_u}{\longrightarrow} \, \UU^+[u]  \, \stackrel{\om}{\longrightarrow} \,  \UU^-[u] \hra \UU^- 
\, \xrightarrow{(\psi^+_{w^{-1}})^{-1}}  \,   \; \ol{S}^{\, +}_{w^{-1}} \subset R_q[B^+ w^{-1} B^+].
\end{equation}
on the second is an antihomomorphism.
\ele

The fact that $\zeta_{w^{-1},u^{-1}}$ is an antihomomorphism follows from the lemma and the fact that 
its restriction to $S^+_w \bowtie S^-_u$
coincides with the composition of the map in the lemma and the 
 homomorphism $\kappa : R_q[B^+ w^{-1} B^+] \to R_q[G^{w^{-1},u^{-1}}]$.

\begin{proof} By \thref{SwSu3}, we have the embedding $S^+_w \bowtie S^-_u \hra S^+_w \bowtie S^-_{w_\ci}$, so it is sufficient to 
prove the lemma for $u = w_\ci$. Furthermore, $S^-_{w_\ci}$ is generated by $(\vp^-_{w_\ci})^{-1}(X_i^+)$ for $i \in [1,r]$, and 
using the embeddings in \thref{SwSu3} one more time, we see that it is sufficient to prove the lemma for $u = s_i$, $i \in [1,r]$.
In this case 
$$
(\vp^-_{s_i})^{-1}(X_i^+) = (q_i^{-1} - q_i)^{-1} \De_{-s_i \vpi_i, - \vpi_i}^{-1} \De_{- \vpi_i, - \vpi_i}.
$$
and the commutation relation \eqref{commRR} reduces to 
\begin{multline*}
\big[ \De_{-s_i \vpi_i, - \vpi_i}^{-1} \De_{- \vpi_i, - \vpi_i} \big] 
\big[ c^+_\xi \De_{w^{-1} \mu, \mu}^{-1} \big]
\\
= q^{ \lcor \al_i, \nu - w^{-1} \mu  \rcor} 
\big[ c^+_\xi \De_{w^{-1} \mu, \mu}^{-1} \big]
\big[ \De_{-s_i \vpi_i, - \vpi_i}^{-1} \De_{- \vpi_i, - \vpi_i} \big] 
+ (q_i^{-1} - q_i) 
\big[ c^+_{S^{-1}(X_i^+) \xi} \De_{w^{-1} \mu, \mu}^{-1} \big] 
\end{multline*}
for $\xi \in (V(\mu)^*)_{- \nu}$, $\mu \in \PP^+$, $\nu \in \PP$.
It is straightforward to check that the composition map in the lemma anti-preserves this commutation relation by using 
\leref{two-zetaw} and computing directly the image of $\De_{-s_i \vpi_i, - \vpi_i}^{-1} \De_{- \vpi_i, - \vpi_i}$. The lemma now follows
from the fact that  the maps $\zeta_{w^{-1}}$ and \eqref{comp} are antihomomorphisms.
\end{proof}

\subsection{Proof of \thref{quantum-twist} (ii)}
\label{7.6} 
First we prove \eqref{zeta-De2}. Similarly to the proof of \leref{two-zetaw} (using Propositions \ref{pomS} (i) and \ref{pzeta-1} and \eqref{iota-ident}),
we obtain
\begin{multline*}
\sqrt{q}^{ \; \| ( w_{\leq k} - 1) \vpi_{i_k} \|^2/2}  \vp_{w^{-1}}^+ \zeta_{w^{-1},u^{-1}} \big( \ol{M}(e_k) \big) = (\De_{\vpi_{i_k}, w_{\leq k} \vpi_{i_k}} \otimes\, \id) ( {}^{\tau \zeta_{w^{-1}}^{-1} \times 1} \RR^{w^{-1}} ) \\
\begin{aligned}
&= (\De_{\vpi_{i_k}, \vpi_{i_k}} \otimes\, \id) \big((T_{w^{-1}}^{-1} \otimes 1) \cdot \; 
{}^{\om S \tau \times 1} \RR^{w^{-1}}  \cdot (T_{w^{-1}_{< k}} \otimes 1) \big) \\
&= ( \De_{\vpi_{i_k}, \vpi_{i_k}} \otimes\, \id) \big((T_{w_{<k}} \otimes 1) \cdot \; 
{}^{\tau \times 1} \RR^{w^{-1}}  \cdot (T_w^{-1} \otimes 1) \big)  
= \De_{w^{-1}_{<k} \vpi_{i_k}, \vpi_{i_k}} \De_{ w^{-1} \vpi_{i_k}, \vpi_{i_k}}^{-1}.
\end{aligned}
\end{multline*}
Next we proceed with \eqref{zeta-De3}. Because of the embeddings 
$S^+_w \bowtie S^-_{u_{\leq j}} \hra S^+_w \bowtie S^-_u$ and the definition of $\zeta_{w^{-1},u^{-1}}$, it is 
sufficient to prove \eqref{zeta-De3} for $j = M$. Set for brevity $i:= i'_M$. We need to prove that
\begin{equation}
\label{pf-zeta-De3}
\sqrt{q}^{ \; \| ( u - w) \vpi_i \|^2/2}
\zeta_{w^{-1},u^{-1}} \big( \ol{M}(e_{N+M}) \big) \cdot \big( \De_{\vpi_i, u \vpi_i} 
\De_{ w^{-1}( u \vpi_i), u \vpi_i}^{-1} \big)^{-1} = 1.
\end{equation}
The frozen cluster variable 
$\ol{M}(e_{N+M})$ is a normal element 
of $R_q[G^{w,u}]$ by \eqref{n1}--\eqref{n2} and  
$\De_{\vpi_i, u \vpi_i} \De_{ w^{-1}( u \vpi_i), u \vpi_i^{-1}}^{-1}$
is a normal element of $R_q[G^{w^{-1},u^{-1}}]$ by the same properties and Propositions \ref{pomS} (ii) and \ref{pSRq} (i).
Using that $\zeta_{w^{-1},u^{-1}}$ is an antiisomorphism gives that the element in the left hand side of \eqref{pf-zeta-De3} is in the 
center of $R_q[G^{w^{-1},u^{-1}}]$. It easily follows from \eqref{zeta-weight} that this element is also in $R_q[G^{w^{-1},u^{-1}}]_{0,0}$.
By \cite[Theorem 4.14 (3)]{HLT}, $Z(R_q[G^{w^{-1},u^{-1}}])_{0,0} = \KK$, thus
\begin{equation}
\label{tt}
\sqrt{q}^{ \; \| ( u - w) \vpi_i \|^2/2}
\zeta_{w^{-1},u^{-1}} \big( \ol{M}(e_{N+M}) \big) \cdot \big( \De_{\vpi_i, u \vpi_i} 
\De_{ w^{-1}( u \vpi_i), u \vpi_i}^{-1} \big)^{-1} = t
\end{equation}
for some $t \in \KK$. For a set of elements $r_1, \ldots r_n$ in an algebra $R$ that generate a quantum torus, we will denote by 
$\TT(r_1, \ldots, r_n)$ this quantum torus.
It follows from \eqref{p_idef} that 
\begin{multline*}
\sqrt{q}^{ \; \| ( u - w) \vpi_i \|^2/2} \, \ol{M}(e_{N+M}) 
- \big( \De_{- u \vpi_i, - \vpi_i}^{-1} \De_{- \vpi_i, - \vpi_i} \big) 
\big( \De_{ \vpi_i, \vpi_i} \De_{w \vpi_i, \vpi_i}^{-1} \big) \\
\in \TT( \{\ol{M}(e_k) \mid k \in [1,N+M-1]\} ).
\end{multline*}
Applying $\zeta_{w^{-1}, u^{-1}}$ to both sides and using \eqref{zeta-De2} and \eqref{zeta-De3} for $j <M$, leads to
\begin{multline}
\label{diff-2}
\sqrt{q}^{ \; \| ( u - w) \vpi_i \|^2/2} \zeta_{w^{-1}, u^{-1}} \big( \ol{M}(e_{N+M}) \big)  -  \\
\big( \De_{ \vpi_i, \vpi_i} \De_{w^{-1} \vpi_i, \vpi_i}^{-1} \big) 
\big( \De_{w^{-1} \vpi_i, u \vpi_i} 
\De_{w^{-1}( u \vpi_i), u \vpi_i}^{-1} \big)  
\in \TT( \{M_{mBZ}(e_k) \mid k \in [1,r+N+M-1]\} ).
\end{multline}
Here we also use that $\zeta_{w^{-1}, u^{-1}}(\De_{ \vpi_i, \vpi_i} \De_{w \vpi_i, \vpi_i}^{-1} ) = 
\De_{ \vpi_i, \vpi_i} \De_{w^{-1} \vpi_i, \vpi_i}^{-1}$, obtained from \eqref{zeta-De2} for $k=N$, and 
$\zeta_{w^{-1}, u^{-1}}(\De_{- u \vpi_i, - \vpi_i}^{-1} \De_{- \vpi_i, - \vpi_i} ) = \De_{w^{-1} \vpi_i, u \vpi_i} 
\De_{w^{-1}( u \vpi_i), u \vpi_i}^{-1}$ which is easily deduced from the definition of $\kappa$ and \prref{omS} (i). 
Hence, \eqref{tt} and \eqref{diff-2} give
$$
t \De_{\vpi_i, u \vpi_i} - \De_{ \vpi_i, \vpi_i} \De_{w^{-1} \vpi_i, \vpi_i}^{-1} \De_{w^{-1} \vpi_i, u \vpi_i} 
\in \TT( \{M_{mBZ}(e_k) \mid k \in [1,r+N+M-1]\} ).
$$
From here one gets that $t=1$ by using that
$$
S( S^-_u ) := \{ c_{\xi_\mu, v} \De_{\mu, u \mu}^{-1} \mid \mu \in \PP^+, v \in V(\mu) \} \subset R_q[G^{w^{-1}, u^{-1}}]
$$
(obtained from \prref{SRq}) and applying the antiisomorphism $\vp^-_{u} S^{-1} : S( S^-_u) \to \UU^+[u]$, cf. \thref{isom}.
\subsection{Proof of \thref{quantum-twist} (i)} 
\label{7.7}
It follows from \eqref{psi-pi} and \eqref{zeta-De3} that 
\begin{equation}
\label{zeta-De1'}
\zeta_{w^{-1},u^{-1}} (\De_{w \mu, \mu}) = \De_{u^{-1} (u \mu), - u \mu}, \quad \mu \in \PP^+.
\end{equation}
Thus, $\zeta_{w^{-1},u^{-1}}$ extends to an antihomomorphism 
$R_q[G^{w, u}] \to R_q[G^{w^{-1},u^{-1}}]$. Analogously to \S \ref{7.6} 
one shows that 
\begin{align}
\label{zeta-De2'}
&\zeta_{w,u}\big( \De_{w^{-1}_{<k} \vpi_{i_k}, \vpi_{i_k}} \big) =
\sqrt{q}^{ \; - \| ( w_{\leq k} - 1) \vpi_{i_k} \|^2/2} \, \ol{M}(e_k) \De_{ - u( u^{-1} \vpi_{i_k}), - u^{-1} \vpi_{i_k}},
\\
\label{zeta-De3'}
&\zeta_{w,u} \big( \De_{\vpi_{i'_j}, u_{\leq j} \vpi_{i'_j}} \big) = 
\sqrt{q}^{ \; - \| ( u_{\leq j } - w) \vpi_{i'_j} \|^2/2} \, \ol{M}(e_{N+j}) \De_{ - u( u^{-1}_{< j} \vpi_{i'_j}), - u^{-1}_{<j} \vpi_{i'_j}}
\end{align}
for all $k \in [1,N]$, $j \in [1,M]$.
Below we prove the first identity, leaving the second to the reader. \prref{omS} (i) and \leref{two-zetaw} imply that
\begin{multline*}
\vp_w^+ \big( \zeta_{w,u} ( \De_{w^{-1}_{<k} \vpi_{i_k}, \vpi_{i_k}} \De_{w^{-1} \vpi_{i_k}, \vpi_k}^{-1} ) \big) = 
(\De_{w^{-1}_{<k} \vpi_{i_k}, \vpi_{i_k}}  \otimes \id) \big( (T_w^{-1} \otimes 1) \cdot \; {}^{\om S \tau \times 1} \RR^w \big)  \\
= (\De_{\vpi_{i_k}, \vpi_{i_k}}  \otimes \id) \big( (T_{w_{\leq k}}^{-1} \otimes 1) \cdot \; {}^{\om S \tau \times 1} \RR^w \big)  
= (\De_{\vpi_{i_k}, \vpi_{i_k}}  \otimes \id) \big(  {}^{\tau \times 1} \RR^w \cdot (T_{(w_{\leq k})^{-1}}^{-1} \otimes 1) \big).
\end{multline*}
This proves \eqref{zeta-De2'} because $X_{\be_l}^+ v_{w_{\leq k} \vpi_{i_k} }=0$ for $l >k$, 
which truncates ${}^{\tau \times 1}\RR^w$ to ${}^{\tau \times 1}\RR^{w_{\leq k}}$.

It follows from \eqref{zeta-De2}--\eqref{zeta-De3} and \eqref{zeta-De2'}--\eqref{zeta-De3'} that 
$\zeta_{w,u}$ and $\zeta_{w^{-1},u^{-1}}$ are inverse homomorphisms because the set $\{\ol{M}(e_k) \mid k \in [1,N+M]\}$ 
generates $\Fract(S^+_w \bowtie S^-_u)$ and the set $\{M_{mBZ}(e_k) \mid k \in [1,N+M+r]\}$ generates 
$\Fract(R_q[G^{w^{-1},u^{-1}}])$ (by the definition of toric frame). This completes the proof of \thref{quantum-twist} (i).

\sectionnew{Quantum cluster algebra structures on $R_q[G^{w,u}]$, the Berenstein--Zelevinsky conjecture}
\label{QcldoubleB}
In this section we relate the quantum seeds for the reduced double Bruhat cells from \thref{SwSu-qcluster} to the BZ quantum seeds 
via the composition of a reduction map, the quantum twist map (Section \ref{twist}) and the map $\om S$ (Section \ref{App}), 
and complete the proof of the Berenstein--Zelevinsky conjecture.

\subsection{Reduction of graded quantum cluster algebras} 
\label{8.1}
We will need an easy general fact on such reductions. We will work in the setting of \S \ref{q-cl}. Let $\PP$ be a free abelian group.

Assume that $R$ is a $\PP$-graded Ore domain.
Let $(M, \wt{B})$ be a $\PP$-graded seed in $\Fract(R)$ by which we mean that every cluster variable in a mutation 
equivalent seed is homogeneous. This is equivalent to saying that the $M(e_k)$ are homogeneous for $k \in [1,N]$ 
and $\deg M(b^k) =0$ for all columns $b^k$ of $\wt{B}$ (i.e., $k \in \ex$).
We will make the following assumption:

(i) The set $\inv$ of inverted frozen indices contains $[1,n]$ for some $n \in [1,N]$ and the function $\phi : \Zset^n \to \PP$, $\phi(f) := \deg M(f)$ 
is an isomorphism.
 
For each quantum seed $(M', \wt{B}')$ mutation equivalent to $(M,  \wt{B})$, we define a quantum seed 
$(M', \wt{B}')_{\redd}:= (M'_{\redd}, \wt{B}'_{\redd})$ of $\Fract (R_0)$. The cluster variables in the latter will be indexed by $[n+1,N]$, the exchangeable 
ones by $\ex$ and the inverted frozen ones by $\inv \backslash [1,n]$. Denote by $\Zset^{[n+1,N]}$ the set of 
integer vectors whose components are indexed by $[n+1,N]$. The matrix $\wt{B}'_{\redd}$ 
is obtained from $\wt{B}'$ by removing the first $n$ rows. The toric frame $M'_{\redd}$ is defined by
$$
M'_{\redd}(f):= M'(f - \phi^{-1}(\nu)) \; \; \mbox{where} \; \; \nu: = \deg M'(f), 
\quad \forall\,  f \in \Zset^{[n+1,N]}.
$$
The matrix $\rb'_{\redd}$ of $M'_{\redd}$ is given by 
$$
(\rb'_{\redd})_{kj} = \Om_{\rb}( e_k - \phi^{-1}(\nu_k), \, e_j - \phi^{-1} (\nu_j)), \quad k,j \in [n+1,N]
$$
(recall \eqref{Om}), where $\nu_k:= \deg M'(e_k)$ and $\rb'$ is the matrix of $M'$. The graded condition 
on $(M, \wt{B})$ yields at once that $(M'_{\redd}, \wt{B}'_{\redd})$ is a quantum seed of $\Fract(R_0)$. 
The definition of reduced toric frames and \eqref{MOm} directly imply that, 
\begin{equation}
\label{MtoMredd}
{\mbox{if $(M'', \wt{B}'') = \mu_k (M', \wt{B}')$, $k \in \ex$, then 
$(M''_{\redd}, \wt{B}''_{\redd}):= \mu_k(M'_{\redd}, \wt{B}'_{\redd})$.}} 
\end{equation}
An illustration of the reduction construction is given in \exref{sl3g} in the next subsection.
The following theorem is inverse to \cite[Theorem 4.6]{GL}.  

\bth{gr} Assume that $(M, \wt{B})$ is a graded quantum seed of the $\PP$-graded Ore domain $R$ that satisfies the condition {\rm{(i)}}. 
Then the quantum cluster algebra associated to the reduced seed $(M_{\redd}, \wt{B}_{\redd})$ satisfies
\begin{align*}
&\AA(M, \wt{B}, \inv) = \AA(M_{\redd}, \wt{B}_{\redd}, \inv \backslash [1,n]) [M(e_k)^{\pm 1}, \, 1 \leq k \leq n], \\
&\UU(M, \wt{B}, \inv) = \UU(M_{\redd}, \wt{B}_{\redd}, \inv \backslash [1,n]) [M(e_k)^{\pm 1}, \, 1 \leq k \leq n].
\end{align*}
In particular, $\AA(M, \wt{B}, \inv) = \UU(M, \wt{B}, \inv)$ if and only if $\AA(M_{\redd}, \wt{B}_{\redd}, \inv \backslash [1,n])=
\UU(M_{\redd}, \wt{B}_{\redd}, \inv \backslash [1,n])$.
\eth

The theorem follows from the fact that every graded subalgebra $A$ of $R$ that contains the frozen variables $M(e_k)^{\pm 1}$, $k \in [1,n]$ satisfies
$A = A_0[M(e_k)^{\pm 1}, \, 1 \leq k \leq n]$ and that
$$
\TT(M', \wt{B}') =  \TT(M'_{\redd}, \wt{B}'_{\redd}) [M(e_k)^{\pm 1}, \, 1 \leq k \leq n]
$$
for every seed $(M', \wt{B}')$ mutation equivalent to $(M, \wt{B})$, recall \S \ref{q-cl}.

One can define an analogous reduction procedure for classical (not quantum) cluster algebras.
\subsection{Cluster structures on quantum double Bruhat cells} 
\label{8.2}
Let $\zeta : R \to A$ be an antiisomorphism of Ore domains. To a quantum seed $(M: \Zset^N \to R, \wt{B})$ of $R$, we associate 
a quantum seed  $\zeta(M, \wt{B}):=(\zeta M, -\wt{B})$ of $A$ as follows. Set $\zeta M(e_k):= \zeta \circ M(e_k)$, for $k \in [1,N]$ and let the matrix of $\zeta M$ be 
$(r_{kj}^{-1})$ where $\rb = (r_{kj})$ is the matrix of $M$. It is obvious that $(\zeta M, -\wt{B})$  is a quantum seed of $A$. We carry over the 
set $\ex$ of exchangeable variables and the set $\inv$ of invertible frozen variables from the first to the second.

For the Berenstein--Zelevinsky quantum seed $(M_{BZ}, \wt{B}_{BZ})$ of $R_q[G^{u,w}]$ (recall \S \ref{6a.1}), we define the 
set of exchangeable indices \eqref{BZ-exch} to be 
$$
\ex_{BZ} := \{ k \in [r+1, r+N+M] \mid s (k) \ne + \infty\}
$$ 
for the successor function $s$ in \S \ref{6a.1} and declare all frozen variables to be invertible, 
$$
\inv_{BZ}:=  [1,r+N+M] \backslash \ex_{BZ}.
$$
For the modified Berenstein--Zelevinsky quantum seeds $(M_{mBZ}, \wt{B}_{BZ})$ of $R_q[G^{w^{-1}, u^{-1}}]$ (see \S \ref{6a.4})
set $\ex_{mBZ}:=\ex_{BZ}$ and $\inv_{mBZ} := \inv_{BZ}$.

We apply the graded reduction from \thref{gr} to the seeds $(M_{BZ}, \wt{B}_{BZ})$ of $R_q[G^{u,w}]$  
and $(M_{mBZ}, \wt{B}_{BZ})$ of $R_q[G^{w^{-1}, u^{-1}}]$. In the first case we choose the grading 
$R_q[G^{u,w}]_{\nu, -}$, $\nu \in \PP$ and in the second $R_q[G^{w^{-1}, u^{-1}}]_{-, \nu}$, $\nu \in \PP$. The frozen variables used for the reduction 
in the two respective cases are 
$$
M_{BZ}(e_i)^{\pm 1} :=  \De^{\pm 1}_{\vpi_i, w^{-1} \vpi_i}, 
\quad
M_{mBZ}(e_i)^{\pm 1} :=  \De^{\pm 1}_{w^{-1} \vpi_i, \vpi_i}, 
\quad i \in [1, r].
$$
\bex{sl3g} 
{
The Berenstein-Zelevinsky quantum seed from \S \ref{6a.1} for $G = SL_3$, 
$w =u = s_1 s_2 s_1$ has cluster variables
\[
\De^q_{1,3}, \; \; \De^q_{12,23}, \quad \; \; \De^q_{1,2}, \; \; \De^q_{12,12}, \; \; \De^q_{1,1}, \quad \; \; 
\De^q_{2,1}, \; \; \De^q_{23,12}, \; \; \De^q_{3,1}
\]
in terms of the standard notation for quantum minors of $R_q[SL_3]$. The mutation quiver of this seed is
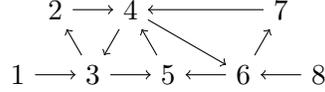
\[
\begin{tikzpicture}
\node (a) at (0,0) {1};
\node (b) at (0.5,1.732/2) {2};
\node (1) at (1,0) {3}; 
\node (3) at (2,0) {5}; 
\node (4) at (3,0) {6}; 
\node (6) at (4,0) {8};
\node (2) at (1.5,1.732/2) {4};
\node (5) at (3.5,1.732/2) {7};
\draw[<-] (1) -- (a);
\draw[<-] (2) -- (b);
\draw[<-] (b) -- (1);
\draw[<-] (3) -- (1);
\draw[<-] (3) -- (4);
\draw[<-] (4) -- (6);
\draw[<-] (2) -- (5);
\draw[<-] (1) -- (2);
\draw[<-] (2) -- (3);
\draw[<-] (4) -- (2);
\draw[<-] (5) -- (4);
\end{tikzpicture}
\]
{This quiver coincides with the quiver of the seed in \cite[Example 2.5]{BFZ} and the quiver of the quantum seed in 
\cite[Example 3.2]{BZ} after a reenumeration of vertices.}
}

{
The corresponding reduced quantum seed has cluster variables
\begin{align*}
&q^{1/2} \De^q_{1,2} (\De^q_{1,3})^{-1}, \; \; q^{1/2} \De^q_{12,12} (\De^q_{12,23})^{-1}, \; \; q^{1/2} \De^q_{1,1} (\De^q_{1,3})^{-1},
\\
&q^{1/2} \De^q_{2,1} \De^q_{1,3} (\De^q_{12,23})^{-1}, \; \; \De^q_{23,12} \De^q_{1,3}, \; \; \De^q_{3,1}  (\De^q_{12,23})^{-1}.
\end{align*}
Its mutation quiver is obtained from the above quiver by removing the vertices 1, 2 and the adjacent edges, and then renaming the other vertices according to 
$k \mt k -2$. However, this is precisely the quiver of the second quantum seed of the algebra $R$ in \exref{sl3e}. Comparing the matrices of the related toric frames, 
one sees that the reduced Berenstein--Zelevinsky quantum cluster algebra for $R_q[SL_3^{s_1 s_2 s_1, s_1 s_2 s_1}]$ is isomorphic to the 
algebra $R$ in \exref{sl3a}. 
}
\eex

\bre{comparison-FWZ}  
{
In \cite[\S 4.1-4.2]{FWZ} Fomin, Williams and Zelevinsky define the notion of a cluster subalgebra of a cluster algebra.
\exref{sl3g} illustrates the similarities and differences between this notion and the reduction construction for graded (quantum) cluster algebras.
}

{
On the one hand, in both cases one removes some rows and the corresponding columns of an exchange matrix (cf. \cite[Definition 4.1.1]{FWZ}),
though in our case we only remove rows corresponding to a collection of frozen variables.
}

{
On the other hand, the cluster variables in a cluster subalgebra are a subset of the original cluster variables and there is a vanishing condition 
on certain entries of an exchange matrix (cf. \cite[Definition 4.2.3]{FWZ}), while in our case the cluster 
variables in the reduced algebra are cluster monomials for the original (quantum) cluster algebra.   
}
\ere

The next facts describe the relations between the quantum seeds of reduced double Bruhat cells and of double Bruhat cells
from Sections \ref{qS+S-} and \ref{App}, respectively.

\bpr{connections} For all simple Lie algebras $\g$ and $w,u \in W$, the following hold:

{\rm(i)} The quantum seed $(\ol{M}, \ol{B})$ of $R_q[G^{w,u}/H]$ in Theorem {\rm\ref{tSwSu-qcluster}  (c)} 
is a reduction of the modified Berenstein--Zelevinsky seed of $R_q[G^{w^{-1},u^{-1}}]$ in the sense that 
$$
(\ol{M}, \ol{B}) = \zeta_{w,u} \big( (M_{mBZ}, \wt{B}_{BZ})_{\redd} \big)
$$
where $\zeta_{w,u} : R_q[G^{w^{-1},u^{-1}}] \to R_q[G^{w,u}]$ is the quantum twist map {\rm(}antiisomorphism{\rm)} from Section {\rm\ref{twist}}. 
 
{\rm(ii)} The quantum seed $(\ol{M}, \ol{B})$ of $R_q[G^{w,u}/H]$ is a reduction of the Berenstein--Zelevinsky seed of $R_q[G^{u,w}]$,
$$
(\ol{M}, \ol{B}) = \zeta_{w,u} \circ (h_* \cdot) \circ (\om S) \big( (M_{BZ}, \wt{B}_{BZ})_{\redd} \big),
$$
recall the isomorphism $\om S : R_q[G^{u,w}] \to R_q[G^{w^{-1}, u^{-1}}]$ and the definition \eqref{hst} of $h_* \in \HH$.

Moreover, the sets $\ex$ and $\inv$ in Theorem {\rm\ref{tSwSu-qcluster} (d)} 
match those for $R_q[G^{w^{-1}, u^{-1}}]$ by converting the index sets from $[1, N+M]$ to $[r+1,r+N+M]$.
\epr

The key part of (i) is that $M(e_k) = \zeta_{w,u} (M_{mBZ})_{\redd}(e_{r+k})$ for $k \in [1,M+N]$. This is exactly 
what is proved in \thref{quantum-twist} (ii) (though written in an explicit form there). The facts that the matrices $\rb$ of the 
toric frames and the exchange matrices match are verified directly. Part (ii) follows from \prref{BZtomBZ} and the first part.

\bth{BZc} {\rm{[Berenstein--Zelevinsky Conjecture]}} For all simple Lie algebras $\g$, Weyl group elements $w, u \in W$, arbitrary base fields $\KK$, 
and non-roots of unity $q \in \kx$ such that $\sqrt{q} \in \KK$, the quantum coordinate ring $R_q[G^{u,w}]$ equals the quantum cluster 
algebra $\AA(M_{BZ}, \wt{B}_{BZ}, \inv_{BZ})$ and the upper quantum cluster algebra $\UU(M_{BZ}, \wt{B}_{BZ}, \inv_{BZ})$.
\eth

The theorem follows from \thref{SwSu-qcluster} (d) via the reduction procedure of \thref{gr} and the connections in \prref{connections}.
Theorems \ref{tSwSu-qcluster} (a)--(b) and \ref{tgr} also give large sets of explicit quantum seeds for the quantum cluster algebra structure
on $R_q[G^{u,w}]$. These seeds are parametrized by $\sig \in \Xi_{N+M}$ (recall \eqref{tau}) 
and the one in the theorem corresponds to $\sig= w_N^\ci$. The cluster variables in these seeds have the form 
$$
q^{m} \big( \zeta_{w,u} (h_* \cdot) \om S ( \wt{M}_\sig(e_k)) \big) \De_{\mu, w^{-1} \mu}, \quad k \in [1, N+M],
$$
see \eqref{Mfr}. Here $m \in \Zset$ and $\mu \in \PP$ depend on $\sig$ and $k$. The exact values of $m$ and $\mu$ can be determined by 
applying the sequence of mutations in \thref{SwSu-qcluster} (b) to get from the seed $(\ol{M}, \ol{B})$ to the seed 
$(M_\sig, \wt{B}_\sig)$ and then using \eqref{MtoMredd} at each step. The result is more technical and we leave the details to
the reader. The principal parts of the exchange matrices of the seeds are the principal parts of the matrices $\wt{B}_\sig$.


\end{document}